\documentclass[11pt]{amsart}

\usepackage{fullpage, amsfonts, amsmath, amscd, amsthm,stmaryrd}
\usepackage[all]{xy}
\usepackage{epsf, graphicx, amssymb, color}

\newtheorem{thm}{Theorem}[section]
\newtheorem{conj}[thm]{Conjecture}
\newtheorem{cor}[thm]{Corollary}
\newtheorem{lem}[thm]{Lemma}
\newtheorem{prop}[thm]{Proposition}

\theoremstyle{definition}
\newtheorem{defin}[thm]{Definition}

\def\c{{\mathcal C}}
\def\f{{\mathcal F}}
\def\h{{\mathcal H}}
\def\s{{\mathfrak s}}
\def\t{{\mathfrak t}}

\def\C{{\mathbb C}}
\def\D{{\mathcal D}}
\def\F{{\mathcal F}}
\def\G{{\mathcal G}}

\def\L{{\mathbb L}}
\def\M{{\mathcal M}}

\def\P{{\mathcal P}}
\def\Q{{\mathbb Q}}

\def\T{{\mathbb T}}
\def\Z{{\mathbb Z}}

\def\del{{\partial}}

\def\mod{{\textup{mod} \;}}

\begin{document}

\title{A spanning tree model for the Heegaard Floer \\ homology of a branched double-cover}

\author{Joshua Greene}
\thanks {Partially supported by an NSF Graduate Fellowship.}
\address {Department of Mathematics, Princeton University\\ Princeton, NJ 08544}
\email {jegreene@math.princeton.edu}

\maketitle

\begin{abstract}
Given a diagram of a link $K$ in $S^3$, we write down a Heegaard diagram for the branched-double cover $\Sigma(K)$.  The generators of the associated Heegaard Floer chain complex correspond to Kauffman states of the link diagram.  Using this model we make some computations of the homology $\widehat{HF}(\Sigma(K))$ as a graded group.  We also conjecture the existence of a $\delta$-grading on $\widehat{HF}(\Sigma(K))$ analogous to the $\delta$-grading on knot Floer and Khovanov homology.
\end{abstract}




\section{Introduction.}\label{S Intro}

Given a link $K \subset S^3$, let $\Sigma(K)$ denote the double-cover of $S^3$ branched along $K$.  This is a closed, oriented $3$-manifold, to which we can associate its Heegaard Floer homology group $\widehat{HF}(\Sigma(K))$, defined by Ozsv\'ath and Szab\'o \cite{OS3mfld1, OS3mfld2}.  When the determinant of $K$ is non-zero, $\Sigma(K)$ is a rational homology sphere, and this invariant takes the form of an abelian group, graded by rational numbers, which decomposes according to $spin^c$ structures on $\Sigma(K)$.  Closely related to this construction are two other invariants: the knot Floer homology $\widehat{HFK}(K)$, defined by Ozsv\'ath and Szab\'o \cite{OSknot} and by Rasmussen \cite{Jakethesis}; and the Khovanov homology $Kh(K)$ \cite{khovanov}.  Both of these knot homology theories are abelian groups with an integer bigrading, one of which is commonly denoted by a $\delta$ in both theories.

These invariants exhibit many similarities, and understanding their precise relationship is an area of active interest.  For example, Ozsv\'ath and Szab\'o have constructed a spectral sequence from the reduced Khovanov homology $Kh(K)$ to $\widehat{HF}(-\Sigma(K))$, using $\Z / 2$ coefficients \cite{OSdoublecover}.  Rasmussen observed an inequality of ranks $\text{rk} \; Kh_\delta (K) \geq \text{rk} \; \widehat{HFK}_\delta (K)$ for many knot types, and conjectured the existence of a spectral sequence between these groups to explain this phenomenon (\cite{Jake}, Section 5).  Such a result would be particularly interesting since $\widehat{HFK}$ is known to distinguish the unknot \cite{OSgenus}, whereas the corresponding fact is unknown for $Kh$.

Both of the knot homology theories $\widehat{HFK}$ and $Kh$ possess {\em spanning tree models} as well.  Given a planar projection $D$ of a link $K$, we can color its regions black and white in checkerboard fashion, and define a graph $B$ whose vertices correspond to the black regions and whose edges correspond to incidences between black regions at the crossings in $D$.  There is a standard construction of a doubly-pointed Heegaard diagram for $K \subset S^3$ out of $D$, and its corresponding Floer chain complex $\widehat{CFK}$ is freely generated by spanning trees of $B$ \cite{OSaltknots}.  In a different direction, by performing reductions on the chain complex appearing in the definition of Khovanov homology, Champanerkar and Kofman \cite{CK} and Wehrli \cite{W} have shown that the same is true for Khovanov homology.  These models are particularly economical in many cases.  For example, when $D$ is a connected alternating diagram, the differentials on $\widehat{CFK}$ and $CKh$ vanish.  Moreover, this model remains useful in making calculations for small non-alternating knots, and for understanding pieces of $\widehat{HFK}$ for other families of knots \cite{OSmutation}.  However, the differential on $\widehat{CFK}$ and $CKh$ in the spanning tree models remains a mystery in general.  We note that the original construction of Khovanov homology is algorithmically computable, and the same is now known for knot Floer homology as well \cite{MOS,OScube}, although the complexes involved in these constructions typically have much larger rank than their homology groups.  It remains a puzzle to describe a simple, algorithmic model for either theory $\widehat{HFK}$ or $Kh$ based on a chain complex generated by spanning trees.

\begin {figure}
\begin {center}
\includegraphics[width=5in]{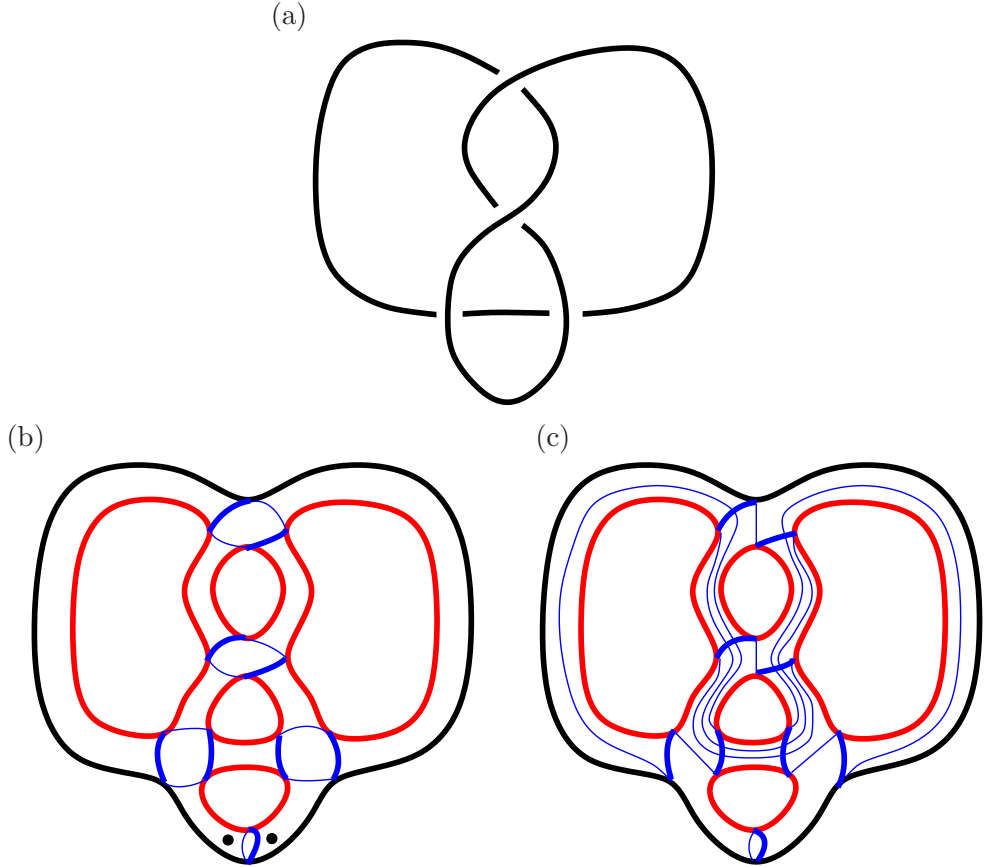}
\put(-270,320){(a)}
\put(-370,160){(b)}
\put(-170,160){(c)}
\caption{(a) A diagram $D$ for the unknot $K$.  (b) The doubly-pointed Heegaard diagram presenting $K \subset S^3$.  (c)  The Heegaard diagram presenting $\Sigma(K)$.  The Heegaard surfaces shown have genus five, and the five red $\alpha$ curves encircle the holes in the surface.  The blue $\beta$ curves appear thinner when they pass onto the bottom of the transparent Heegaard surface. } \label{f intro}
\end {center}
\end {figure}

The main result of this paper is a spanning tree model for $\widehat{HF}(\Sigma(K))$ which is similar in spirit to the spanning tree model for $\widehat{HFK}(K)$.  Figure~\ref{f intro} gives an example of the Heegaard diagram involved.  Shown there is a diagram $D$ for the unknot $K$, the associated doubly-pointed Heegaard diagram for $K \subset S^3$ (following \cite{OSaltknots}), and a Heegaard diagram for $\Sigma(K)$.  Comparing the two pictures, observe that the $\alpha$ and $\beta$ curves and their intersection points can be put into one-to-one correspondence between the two Heegaard diagrams, and that the top half of both diagrams are identical.  The first observation indicates that both of the associated Floer chain complexes have the same generating set, while the second is a bit more subtle, leading to a pair of spectral sequences which we examine in $\S$\ref{S specseq}.  This model for $\widehat{HF}(\Sigma(K))$ suffers the same drawback as does the model for $\widehat{HFK}(K)$, in that it does not lead to an algorithmic calculation of the homology group.  However, we can use it to make computations for many manifolds of interest, and it suggests the existence of a $\delta$-grading on $\widehat{HF}(\Sigma(K))$ analogous to the ones on $\widehat{HFK}$ and $Kh$.
\\

This paper is organized as follows.  In $\S$\ref{S prelim} we discuss some of the basic combinatorics of knot diagrams.  In particular, we define the Goeritz matrices and coloring matrix associated to a knot diagram, and we introduce Kauffman states to encode spanning trees of the black graph.  In $\S$\ref{S HDs} we construct two Heegaard diagrams for $\Sigma(K)$ out of a planar diagram $D$ of $K$.  The first of these, $\h_1(D)$, arises out of the standard cut-and-paste description of $\Sigma(K)$.  The second, $\h_2(D)$, is the one alluded to already.  It is constructed by resolving $D$ into a single unknotted curve with a collection of resolving arcs.  The branched double-cover of the unknotted curve is $S^3$, and the preimages of the resolving arcs gives rise to a simple framed link $\L \subset S^3$.  Surgery on $\L$ produces $\Sigma(K)$, and by writing down a Heegaard diagram subordinate to $\L$ we obtain $\h_2(D)$.  This procedure is depicted in Figure~\ref{f HDderivation}.  Using this diagram, we give a quick proof that the branched double-cover of a non-split alternating link is an L-space.

In $\S$\ref{S absgradings} we write down a simple formula for the grading and $spin^c$ structure of a Kauffman state $x$ (Theorem \ref{t grading}), interpreted as a generator of $\widehat{CF}(\h_2(D))$.  One of the terms in the grading formula is a quantity $\delta(x)$ akin to the $\delta$-grading on $\widehat{HFK}$ and $Kh$.  After proving these formulas, we show in $\S$\ref{ss alternating} a simple way to compute the correction terms for $\Sigma(K)$ when $K$ is a non-split alternating link (cf. \cite{OSdoublecover}, Theorem 3.4).  In $\S$\ref{S domain} we use the Goeritz matrices, coloring matrix, and the results of $\S$\ref{ss regions} to write down the domain of a homotopy class connecting a pair of generators (Theorem \ref{t domain}).

In $\S$\ref{S specseq} we revisit the similarity between the diagram $\h_2(D)$ and the standard doubly-pointed Heegaard diagram for $K \subset S^3$ gotten from $D$.  We write down a pair of spectral sequences whose $E_0$ terms are freely generated by Kauffman states, and which converge to $\widehat{HF}(\Sigma(K))$ and $\widehat{HFK}(K)$.  Moreover, the $d_0$ differential in both spectral sequences counts holomorphic disks whose domain is supported on top of the Heegaard diagram, so the $E_1$ term in both sequences are the same.  Theorem \ref{t E1} asserts that this group is freely generated by so-called {\em solitary} Kauffman states of $D$ (Definition \ref{d solitary}).  We note that by construction the $d_0$ differential preserves the Alexander grading in $\widehat{CFK}$, and respects the decomposition of $\widehat{CF}(\Sigma(K))$ into $spin^c$ structures.  In $\S$\ref{S calculations} we apply the spanning tree model of $\widehat{HF}(\Sigma(K))$ to make calculations for some knots with $\leq 10$ crossings.  Specifically, we apply it to those knots which are neither alternating nor Montesinos knots, since the invariant for these knots can be computed by other means (as in $\S$\ref{ss alternating} and \cite{OSdoublecover, OSplumbed}).

We conclude in $\S$\ref{S speculation} with some speculation.  We focus on the term $\delta(x)$ appearing in the grading formula Theorem \ref{t grading} and conjecture the existence of a natural $\delta$-grading on $\widehat{HF}(\Sigma(K))$, analogous to the one on $\widehat{HFK}(K)$ and $Kh(K)$.  In support of this conjecture are some small calculations of $\widehat{HF}(Y)$, where $Y$ is a rational homology sphere arising as the boundary of a negative definite plumbing on a tree.  In these examples, the differential in the spanning tree model behaves surprisingly nicely, and suggests an algorithmically computable model for $\widehat{HF}(Y)$.  We also point out a strong similarity with work by N\'emethi \cite{nemethi}, who has proposed an alternative algorithmically computable model for $HF^+(Y)$, which comes equipped with an additional integer grading reminiscent of our $\delta$.
\\

We point out that Heegaard diagrams presenting $\Sigma(K)$, and more general $m$-fold cyclic branched covers $\Sigma_m(K)$, have appeared in the work of Grigsby \cite{grigsby} and Levine \cite{levine}.  In their work, the approach is to begin with a doubly- or multiply-pointed Heegaard diagram presenting $K \subset S^3$, form the appropriate cyclic branched cover of the Heegaard surface, and thereby obtain a Heegaard diagram presenting the preimage of the link $\widetilde{K} \subset \Sigma_m(K)$.  This approach is particularly well-suited to the calculation of the knot Floer homology $\widehat{HFK}(\widetilde{K} \subset \Sigma_m(K))$.  By contrast, our constructions are specific to the case $m=2$, and do not immediately present the preimage $\widetilde{K} \subset \Sigma(K)$.

Another Heegaard diagram for $\Sigma(K)$ was described by Manolescu (\cite{Mnil}, Section 7).  In that setting, the link $K$ is first presented as a plat closure.  Manolescu then shows how the set of generators used by Bigelow \cite{bigelow} to compute the Jones polynomial $V(K)$ can also be used as a generating set for both the Seidel-Smith ``symplectic Khovanov cohomology'' $Kh_{symp}(K)$ \cite{SS} and the Floer homology $\widehat{HF}(\Sigma(K) \# (S^1 \times S^2))$.  The triple use of the Bigelow generators here bears a surface similarity to that of the spanning trees described above.


\subsection*{Acknowledgment.}  It is a pleasure to thank my advisor, Zolt\'an Szab\'o, for his patience and guidance during this project.  I am also grateful to John Baldwin, Eli Grigsby, Adam Levine, and Zhong Tao Wu for conversations about this work.  Thanks lastly to Jake Rasmussen for drawing my attention to the work \cite{Mnil} in his excellent course at Princeton, Spring 2008.


\vspace{.5in}

\section{Preliminaries on knot diagrams.}\label{S prelim}

\subsection{Diagrams, graphs, and the Goeritz form.}\label{ss prelim_1}

Consider a connected diagram $D$ of a link $K$ with a marked point $p$ on one of its edges.  The diagram $D$ splits the plane into connected regions, which we color white and black in checkerboard fashion.

\begin{figure}
\centering
\includegraphics[width=3in]{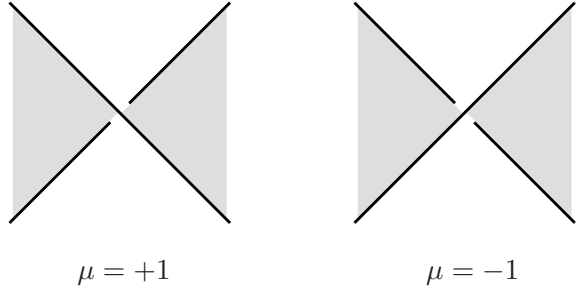}
\put(-190,-20){$\mu = +1$}
\put(-58,-20){$\mu = -1 $}
\caption{The incidence number of a crossing.}  \label{f crossingincidence}
\end{figure}

\begin{figure}
\centering
\includegraphics[width=6in]{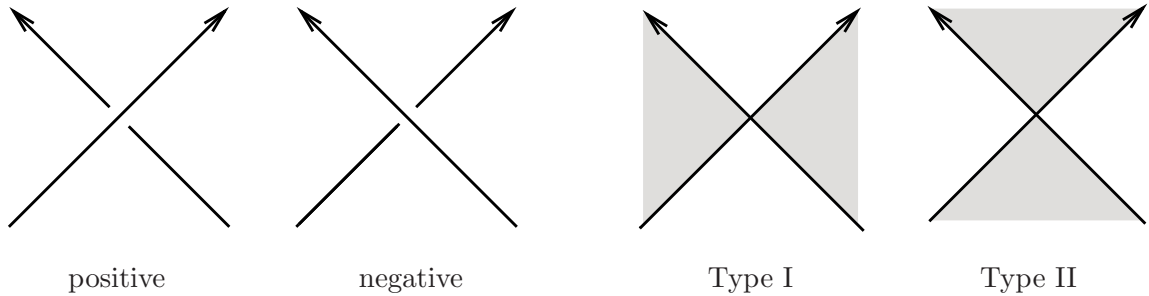}
\put(-410,-20){positive}
\put(-300,-20){negative}
\put(-168,-20){Type I}
\put(-65,-20){Type II}
\caption{The sign and type of a crossing.}  \label{f crossingsignandtype}
\end{figure}

With respect to this coloration, each crossing $c$ in $D$ has an incidence number $\mu(c) = \pm 1$ given as in Figure \ref{f crossingincidence}.  When the link $K$ is oriented, each crossing additionally has a sign and a type (I or II) given as in Figure \ref{f crossingsignandtype}.

We form a planar graph $B$ by drawing a vertex in every black region and an edge for every crossing that joins two black regions.  Dually, we obtain a graph $W$ on the white regions.  In both $B$ and $W$, we associate the label $\mu(e) := \mu(c)$ to the edge $e$ corresponding to the crossing $c$, and mark the vertex in the region adjacent to the marked point $p$ in $D$. We call these vertices {\em roots} and denote them by $r_B$ and $r_W$.  We refer to these decorated plane drawings $B$ and $W$ as the \emph{black} and \emph{white} graphs corresponding to $D$.  Deleting the roots and their incident edges result in the {\em reduced} black and white graphs $\widetilde{B}$ and $\widetilde{W}$.

\begin{figure}
\begin{center}
\includegraphics[width=6in]{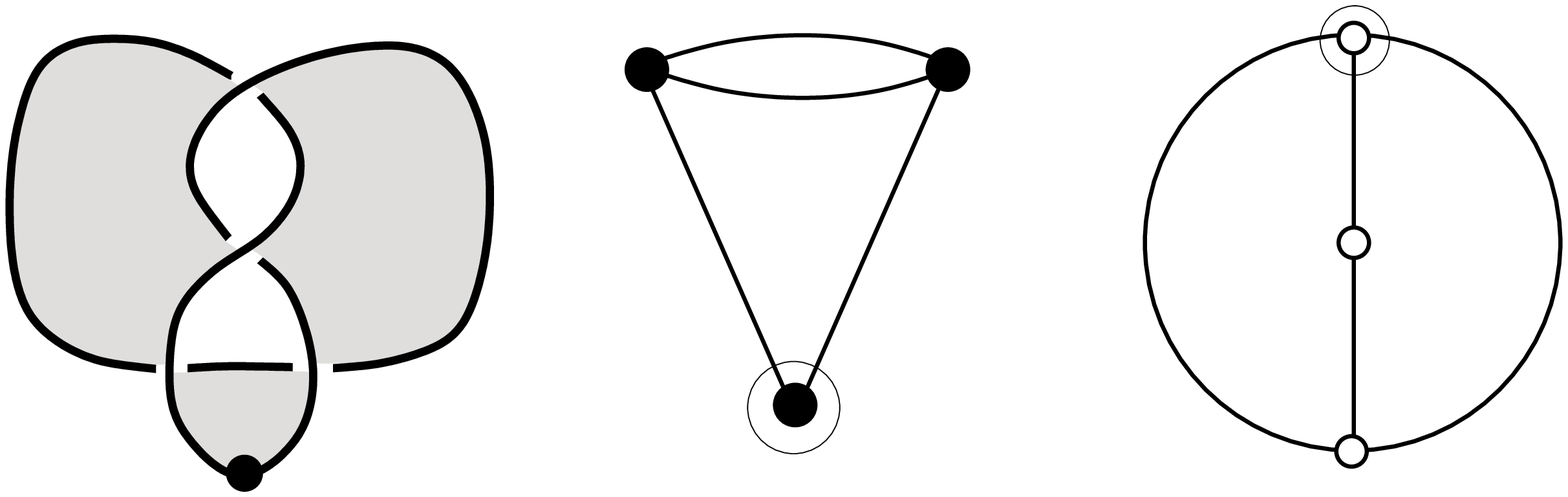}
\put(-355,-5){$p$}
\put(-450,120){$D$}
\put(-280,120){$B$}
\put(-115,120){$W$}
\put(-195,5){$r_B$}
\put(-45,130){$r_W$}
\put(-220,130){$-1$}
\put(-220,95){$-1$}
\put(-258,65){$+1$}
\put(-185,65){$-1$}
\put(-75,95){$-1$}
\put(-75,40){$-1$}
\put(-137,65){$+1$}
\put(5,65){$-1$}
\caption{A diagram $D$ for the unknot $K$ with a checkerboard coloring of its regions, along with its corresponding black graph $B$ and white graph $W$.} \label{f blackwhitegraph}
\end{center}
\end{figure}

The {\em Goeritz matrix} $G_B = (g_{ij})$ of the black graph $B$ is defined as follows.  Enumerate the vertices of $\widetilde{B}$ by $v_1,\dots,v_n$, and set $$g_{ij} = \sum_{e \text { joining }  v_i \text { and } v_j} \mu(e), \quad i \ne j, \quad \quad \text{and} \quad \quad g_{ii} = -\sum_{e \text { incident } v_i} \mu(e).$$  For example, if we label the top left vertex of the graph $B$ in Figure \ref{f blackwhitegraph} by $v_1$ and the top right vertex $v_2$, then $$G_B = \left( \begin{matrix} 1 & -2 \\ -2 & 3 \end{matrix} \right).$$  The Goeritz matrix $G_W$ of the white graph $W$ is defined analogously.

We mention a few important properties of the Goeritz matrix $G = G_B$.  It is a symmetric $n \times n$ matrix, and $|\det(G)| = \det(K)$.  Orienting $K$ somehow, we can define $$\mu(D) = \sum_{c \text { of type II}} \mu(c);$$ then Gordon and Litherland \cite{GL} have shown that the signature of $K$ can be computed as $$\sigma(K) = \sigma(G) - \mu(D).$$  (This conforms with the somewhat backwards convention that a positive link has {\em negative} signature.)  The Goeritz matrix is also a presentation matrix for $H_1(\Sigma(K);\Z) \cong H^2(\Sigma(K);\Z)$, a fact that we will deduce in $\S$\ref{ss h_1 1}.

\subsection{The coloring matrix.}\label{ss prelim_2}

When a crossing of $D$ is marked, we also get an associated {\em coloring matrix} $A = (a_{ij})$.  Enumerate the unmarked crossings of $D$ by $c_1,\dots,c_m$ and the unmarked arcs by $\alpha_1,\dots,\alpha_m$.  Then we set

$$ a_{ij} = \left\{
\begin{array}{cl}
\mu(c_i) & , \quad \mbox{if $\alpha_j$ is an understrand at $c_i$} \cr
& \cr
-2 \mu(c_i) & , \quad \mbox{ if $\alpha_j$ is an overstrand at $c_i$} \cr
& \cr
0 & , \quad \mbox{otherwise}.
\end{array}
\right.
$$  For example, orient the diagram $D$ in Figure \ref{f blackwhitegraph} out of the point $p$ and to the left.  Following this orientation out of $p$, let $c_1,c_2,c_3$ denote the first three crossings of $D$, encountered in that order; and $\alpha_1, \alpha_2,\alpha_3$ the three arcs of $D$, traversed in that order, after the marked one.  Marking the remaining unused crossing, we obtain the coloring matrix $$A = \left( \begin{matrix} 1 & 1 & 0 \\ 0 & 0 & -1 \\ -1 & 0 & 2 \end{matrix} \right).$$

The coloring matrix gets its name because of its relationship with $n$-colorings of the knot diagram (not to be confused with the checkerboard coloring!).  Recall that an {\em $n$-coloring} is a mapping $\{ \mbox{arcs of D} \} \to \Z / n \Z$ with the property that at every crossing, the identity $2b \equiv a + c \, (\mod n)$ holds, where $b$ denotes the value on the overstrand, and $a$ and $c$ the values on the two understrands.  The $n$-colorings of the diagram which take the value $0$ on the marked arc are precisely the elements in $\ker(A) \; (\mod n)$.  The matrix $A$ is typically {\em not} symmetric, but like the Goeritz matrix, it has the property that $|\det(A)| = \det(K)$, and it gives a presentation for $H_1(\Sigma(K);\Z)$.

\subsection{Kauffman states.}\label{ss prelim_3}

\begin {figure}
\begin {center}
\includegraphics[width=6in]{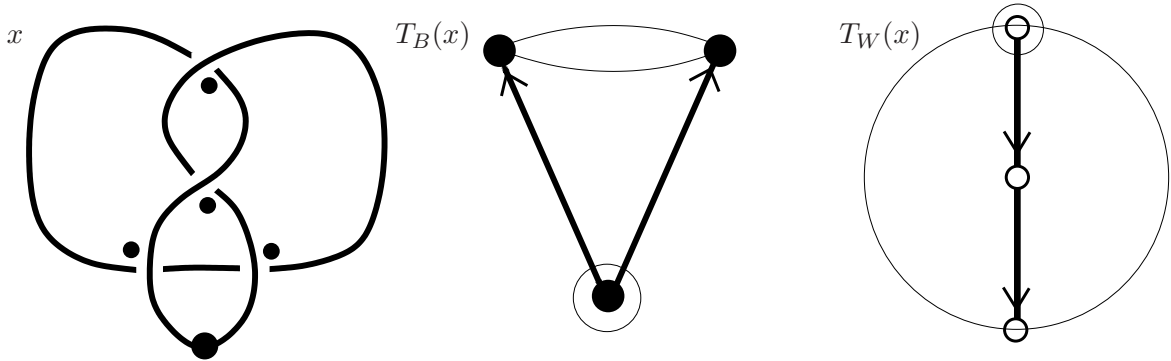}
\put(-440,120){$x$}
\put(-293,120){$T_B(x)$}
\put(-125,120){$T_W(x)$}
\caption{A Kauffman state $x$ of the diagram $D$ pictured in Figure~\ref{f blackwhitegraph}, along with the oriented trees $T_B(x)$ and $T_W(x)$.}\label{f kauffman}
\end {center}
\end {figure}

Every crossing in $D$ is incident a corner of four (not necessarily distinct) regions.  A {\em Kauffman state} $x$ of $D$ is a matching between its crossings and corners of {\em unmarked} regions, so that each crossing gets paired with one of its incident corners.  The Kauffman state $x$ can be visualized by placing a small marker nearby each crossing $c$ in the corner that $c$ gets paired with.  See Figure~\ref{f kauffman}.  We write $x(c)$ to indicate the corner of the region with which $c$ gets paired and $x(R)$ to denote the crossing with which the region $R$ gets paired.

There is a 1-1 correspondence between Kauffman states of $D$ and spanning trees of $B$ (and $W$).  Given a Kauffman state $x$, consider the edge $e \in E(B)$ for each crossing $c$ for which $x(c)$ is black; the collection of these edges is a spanning tree $T_B(x)$ of the black graph.  Moreover, $T_B(x)$ inherits a natural orientation from $x$, gotten by directing each edge $e$ to point towards the region containing $x(c)$.  This is the same as orienting $T_B(x)$ out of its root $r_B$ -- that is, directing every edge of $T_B(x)$ to point towards its endpoint which is further away from $r_B$ inside $T_B(x)$.  In the same way, $x$ gives rise to the spanning tree $T_W(x)$ in the white graph, which is planar dual to $T_B(x)$.  Conversely, starting with a spanning tree $T$ of $B$, we may construct its planar dual $T^*$, and orient each of these trees out of their respective roots.  We obtain a Kauffman state $x(T)$ as follows: at each crossing $c$, there is a corresponding edge $e \in E(T) \cup E(T^*)$, and we take $x(c)$ to lie in the corner of the region to which $e$ is directed.

\vspace{.5in}

\section{Heegaard diagrams for the branched double-cover.}\label{S HDs}

In this section we describe two Heegaard diagrams for the branched double-cover of a link $K \subset S^3$, beginning with a connected planar projection $D$.  In $\S$\ref{ss h_1 1} we define the Heegaard diagram $\h_1(D)$ based on a surgery diagram for $\Sigma(K)$ written down by Ozsv\'ath and Szab\'o \cite{OSdoublecover}.  Using it we confirm that the Goeritz matrix for $D$ presents $H_1(\Sigma(K))$, and we obtain a simple presentation for $\pi_1(\Sigma(K))$.  In $\S$\ref{ss h_1 2} we recover $\h_1(D)$ from the standard cut-and-paste description of $\Sigma(K)$, and fit it into a Heegaard triple presenting the double-cover of $D^4$, branched along a spanning surface for $K$.  In $\S$\ref{ss h_2} we turn to our main construction of the Heegaard diagram $\h_2(D,T)$.  We show that $\widehat{CF}(\h_2(D,T))$ is generated by Kauffman states of $D$, and examine what happens when $D$ is alternating.  In $\S$\ref{ss regions} we examine the regions of $\h_2(D,T)$, work which is needed when we determine the domain of a homotopy class in $\S$\ref{S domain}.  In $\S$\ref{ss h_2 indep} we show that $\h_2(D,T)$ is independent of the choice of spanning tree $T$ used in its definition.  Finally, in $\S$\ref{ss heegaardmoves} we examine a sequence of Heegaard moves relating the diagrams $\h_1(D)$ and $\h_2(D)$.  This relationship enables us to identify a subset of the generators of $\widehat{CF}(\h_1(D))$ with Kauffman states of $D$, and use them in the determination of the grading on $\widehat{HF}(\Sigma(K))$ in $\S$\ref{S absgradings}.

\subsection{The Heegaard diagram $\h_1(D)$.}\label{ss h_1 1}

Ozsv\'ath and Szab\'o \cite{OSdoublecover} explain how to use the graph $\widetilde{W}$ to produce a surgery diagram for the branched double-cover $\Sigma(K)$ (they actually use the black graph $\widetilde{B}$, but this difference is cosmetic). Specifically, label each vertex $w_i \in V(\widetilde{W})$ with the value $g_{ii}$ and center a round unknot at it.  For every edge $e$ between two vertices, we add a right-/left-handed clasp between the corresponding unknots according as $\mu(e) = \pm 1$.  Finally, we frame each unknotted component by the label on its vertex.  Denote the resulting framed link by $\L$; then $\Sigma(K) = S^3(\L)$.

\begin{figure}[htb!]
\centering
\includegraphics[height=6cm]{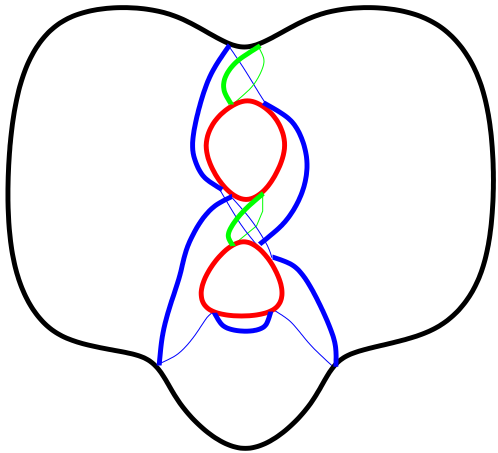}
\caption{The Heegaard triple $\h_1$ for the diagram $D$ in Figure \ref{f blackwhitegraph}.  Omitting the green $\beta$ curves produces a Heegaard diagram presenting $\Sigma(K)$.}  \label{f HD1}
\end{figure}

This dual description neatly leads to a Heegaard triple subordinate to $\L$ (in the sense of \cite{OS4mfld}, Section 4), and in particular to a Heegaard diagram for $\Sigma(K)$.  First, let $\Sigma$ denote the boundary of a closed regular neighborhood $V$ of $W$, let $\alpha = \{ \alpha_i \}$ denote the collection of curves obtained on intersecting $\Sigma$ with the unmarked faces of $W$, and let $\beta' = \{ \beta_i' \}$ denote a collection of curves on $\Sigma$ chosen so that each $\beta'_i$ meets $\alpha_i$ geometrically once and avoids all other $\alpha_j$.  Thus the $\alpha$-handlebody $U_\alpha$ describes the exterior of $V$, while $U_\beta$ describes $V$ itself.  Next, let $\gamma'$ denote the collection of curves on $\Sigma$ obtained by pushing the $\alpha$ curves slightly into the upper half-space determined by the plane of the diagram.  Finally, choose a meridian on $\Sigma$ around each edge $e$ of $W$, and perform a $\mu(e)$-Dehn twist along it. Let $\beta$ and $\gamma$ denote the images of $\beta'$ and $\gamma'$ following these Dehn twists.  Then $(\Sigma,\alpha,\beta,\gamma)$ is a Heegaard triple subordinate to $\L$, and $Y_{\alpha \gamma} \cong \Sigma(K)$.  Verification of this fact is straightforward; the key points are that $V$ is equal to a regular neighborhood of a bouquet of the link $\L$, and that the collection of curves $\gamma$ is an isotopic copy of the framed link $\L$ drawn on the surface $\Sigma$.  We denote this Heegaard diagram for $\Sigma(K)$ so derived from $D$ by $\h_1 = \h_1(D)$ (Figure~\ref{f HD1}).

The diagram $\h_1$ permits an easy calculation of the first homology and fundamental group of $\Sigma(K)$.  We begin by reviewing the general procedure of how this works.  Given a Heegaard diagram $\h = (\Sigma,\alpha,\beta)$ for $M$, enumerate its $\alpha$ curves $\alpha_1,\dots,\alpha_n$ and $\beta$ curves $\beta_1,\dots,\beta_n$.  Orient these curves somehow, and fix an orientation on the surface $\Sigma$.

\begin{defin}\label{d intersectionmatrix}
Given these orientations, the matrix of intersection numbers $(\#(\alpha_i \cap \beta_j))$ is called the {\em intersection matrix of $\h$}, denoted $I(\h)$.
\end{defin}  \noindent The importance of the definition is that $I(\h)$ is a presentation matrix for $H_1(M)$.  A presentation for $\pi_1(M)$ is given as follows.  We associate a generator $x_i$ to each $\alpha_i$ and a relation $r_j$ to each $\beta_j$.  The relation $r_j$ is gotten by traversing one full circuit around the curve $\beta_j$: each time $\beta_j$ encounters a curve $\alpha_i$, we write down a generator $x_i^{\pm 1}$, the sign of the exponent chosen according to the intersection number between $\alpha_i$ and $\beta_j$ at that point.  The product of these generators in order is the word $r_j$. Of course, this word depends on the initial point from which we traverse $\beta_j$, but two different points give rise to conjugate words, and we will blur this small indeterminacy.  In total, we obtain the presentation $\pi_1(M) \cong \langle x_1, \dots, x_n \; | \; r_1, \dots, r_n \rangle$.

Returning to the case at hand of $\h_1$, observe that $\#(\alpha_i \cap \beta_j) = g_{ij}$ for all $i,j$. It follows that $I(\h_1)$ agrees with the Goeritz matrix $G$, and we recover the well-known fact that $G$ is a presentation matrix for $H_1(\Sigma(K))$.  For the fundamental group, we write down a generator $x_i$ for each unmarked vertex $v_i \in B$, and for every edge $e$ incident the vertex $v_j$, we write down the word $(x_j^{-1}x_k)^{\mu(e)}$, where $v_k$ denotes the other endpoint of $e$.  When the neighbor $v_k$ is the marked vertex, we understand this word by taking $1$ in place of $x_k$.  Traverse a small clockwise circuit around the vertex $v_j$ and multiply these words together in the order the edges incident $v_j$ are encountered.  The result is a relation $r_j$, and $\h_1$ induces the presentation $\pi_1(\Sigma(K)) \cong \left< x_1,\dots,x_n \; | \; r_1,\dots,r_n \right>$.  In the example of Figure \ref{f blackwhitegraph}, we obtain the presentation $\langle x_1, x_2 \; | \; (x_1 x_2^{-1})^{-1}(x_1 x_2^{-1})^{-1}(x_1)^{+1}, (x_2 x_1^{-1})^{-1}(x_2 x_1^{-1})^{-1}(x_2)^{-1} \rangle$ for $\pi_1(\Sigma(K))$. It is easy to check that this is a presentation of the trivial group, and this is consistent with the fact that $\Sigma(K) \cong S^3$.

\vspace{.2in}

\subsection{A second pass at $\h_1(D)$.}\label{ss h_1 2}

There is a conceptually clearer derivation of $\h_1$ which follows the standard cut-and-paste description of $\Sigma(K)$, which we presently recall (see also \cite{licky}, p. 85 and \cite{GL}, p. 56).  From the black regions of the diagram $D$, we obtain a spanning surface $F := F_B$ for the link $K$: it consists of the planar black regions away from the crossings, along with a half-twisted band nearby each crossing in $D$.  Form a closed regular neighborhood $U$ of $int(F)$ with spine equal to the link $K$.  Its boundary is the orientable double-cover of $F$ branched along $K$, and so comes equipped with an involution $i$.  We take two copies of $S^3 - int(U)$, identified by some homeomorphism $j$, and glue them together along their boundaries by means of the homeomorphism $j \circ i| \del (S^3 - int(U))$, and the result is $\Sigma(K)$.

\begin{figure}[htb!]
\centering
\includegraphics[height=6cm]{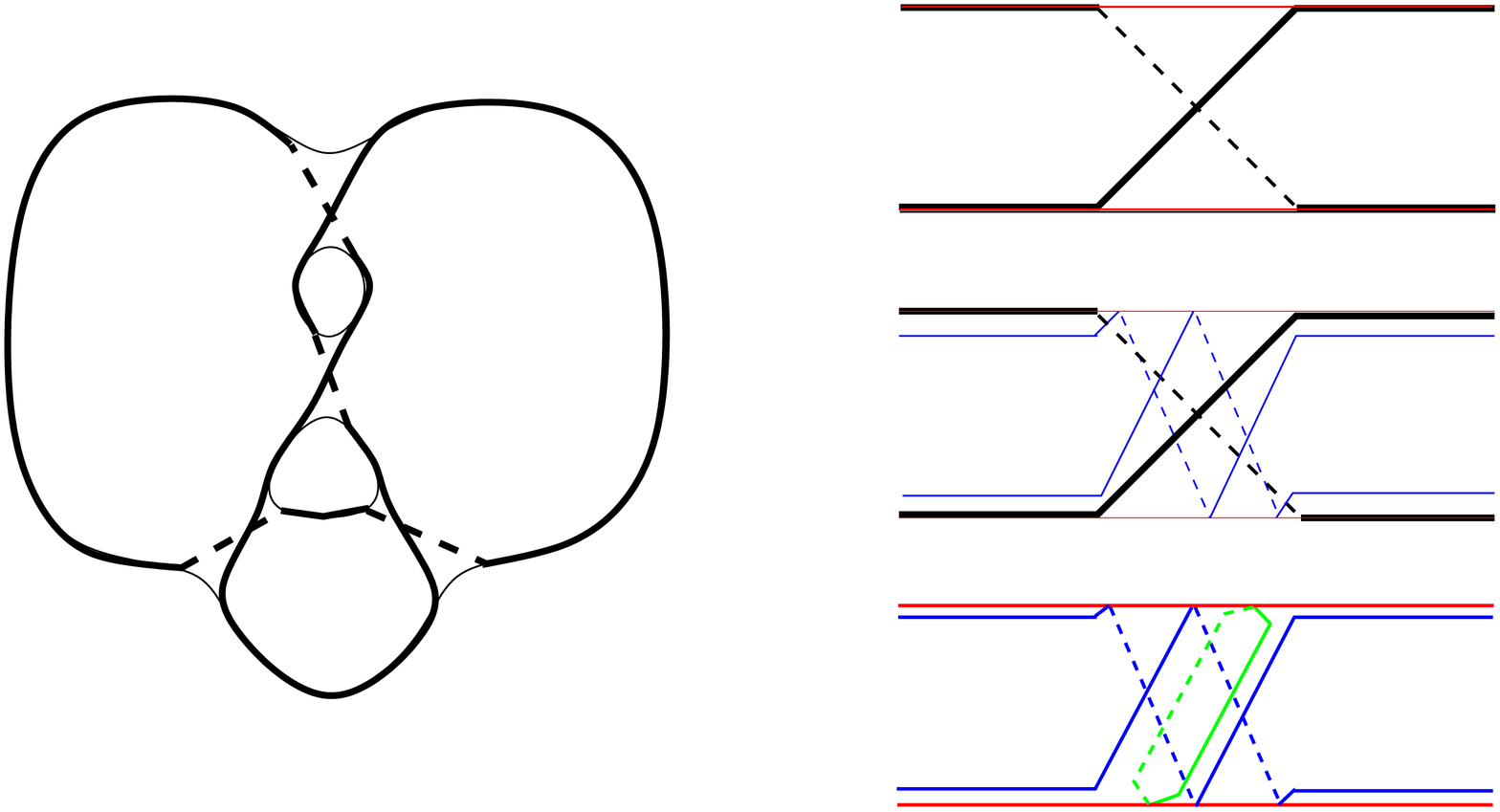}
\put(-350,150){(a)}
\put(-150,170){(b)}
\put(-150,110){(c)}
\put(-150,50){(d)}
\caption{(a) The link $K$ drawn on $\Sigma$. (b) Nearby a crossing.  The red arcs are pieces of $\alpha$ curves. (c) The effect of the involution $i$ on the $\alpha$ curves, perturbed to give (blue) $\gamma$ curves. (d)  The local picture of the Heegaard triple nearby the crossing.}  \label{f H1 derivation}
\end{figure}

We spell out this construction further in order to clarify how it gives rise to a Heegaard diagram for $\Sigma(K)$.  Resolve every crossing of $D$ so that its two incident black regions merge, yielding a collection $C$ of black regions in the plane.  Form the product $[-1,1] \times C$ and quotient it by collapsing the segment $[-1,1] \times p$ to the point $(0,p)$ for each point $p \in \del C$.  The result is the handlebody $U$, which is naturally embedded in $S^3$ so that $0 \times C$ is identified with $C$.  Notice that, away from its crossings, $D$ coincides with the intersection of $U$ with the plane of the diagram.  Nearby a crossing, we can push the overstrand of $D$ onto the top half of $\Sigma := \del U$ and the understrand onto its lower half.  In this way we obtain an embedding of the link $K$ in $\Sigma$ (Figure~\ref{f H1 derivation}(a)), which extends to an embedding of $F$ in $U$.  Next we describe an involution $i$ on $U$ which fixes $F$.  Away from the crossings of $D$, $i$ interchanges points $(\pm t, p)$ for $t \in [-1,1]$, $p \in U$.  Nearby a crossing, the pair $(U,F)$ is locally modelled by the pair $(D^2 \times I, [-1,1] \times I)$ (viewing $D^2 \subset \C$).  With respect to this model, $i$ is given by identifying points $(z,t)$ and $(\overline{z},t)$.

Let $\alpha$ denote the collection of curves obtained on intersecting $\Sigma$ with the unmarked white regions of the link diagram.  Thus the pair $(\Sigma, \alpha)$ describes the handlebody which is the complement of $U$ in $S^3$.  In order to form $\Sigma(K)$, we take two of these handlebodies and identify their boundaries by means of the homeomorphism $i$.  Let $\gamma$ denote the image $i(\alpha)$, perturbed so as to meet $\alpha$ transversally.  Then $(\Sigma,\alpha,\gamma)$ is a Heegaard diagram for $\Sigma(K)$.  In order to describe the $\gamma$ curves explicitly, observe that $i$ fixes the $\alpha$ curves away from the crossings of $D$.  Nearby a crossing, the effect of $i$ is as pictured in Figure~\ref{f H1 derivation}(c).  We get a collection of curves $\gamma$, well-defined up to isotopy, by pushing $i(\alpha)$ off of $\alpha$ into the top half of $\Sigma$ away from the crossings of $D$, and extending in the obvious way nearby the crossings of $D$.  One sees at once that $(\Sigma,\alpha,\gamma)$ agrees with the diagram $\h_1$ described above.

The $3$-manifold $\Sigma(K)$ naturally bounds $\Sigma(F)$, the double-cover of the four-ball $B^4$ branched along a pushed-in copy of the surface $F$, and we can adapt the preceding construction to a Heegaard triple which describes this $4$-manifold.  This triple will ultimately enable us to calculate the absolute gradings on the Floer chain complex, as well as its decomposition into $Spin^c$ structures ($\S$\ref{S absgradings}).  We begin with a careful description of $\Sigma(F)$.  Form the product $[-1,1] \times S^3$ and remove $(0,1] \times int(U)$ from $[-1,1] \times S^3$.  The result is homeomorphic to the original $[-1,1] \times S^3$, but now the new $1 \times S^3$ boundary component has been decomposed into $U$, a collar neighborhood $N$, and the complement.  We extend the involution $i$ on $U$ to one on $U \cup N$ in the obvious way.  Now form two copies of the modified $[-1,1] \times S^3$ and identify the distinguished $U \cup N$ contained within the boundaries of each by means of the homeomorphism $i$.  The result is the $4$-manifold $\Sigma(F)$ with two small balls removed.

Now we show how to extract a Heegaard triple for $\Sigma(F)$ from this description.  Decompose $S^3 = U_\beta \cup_\Sigma U_\alpha$, now writing $U_\beta$ in place of $U$.  Let $\Delta$ denote a triangle with edges $e_\alpha, e_\beta, e_\delta$ in clockwise order, and vertices $v_{\alpha \beta}, v_{\beta \delta}, v_{\delta \alpha}$ at which the corresponding edges meet.  Form the product $\Delta \times \Sigma$, and glue $I \times U_\alpha$ onto $e_\alpha \times \Sigma$ and $I \times U_\beta$ onto $e_\beta \times \Sigma$.  Smoothing its corners, this space is diffeomorphic to the modified $[-1,1] \times S^3$ from the previous paragraph.  To aid in seeing this, note that its boundary consists of $U_\alpha \cup_{v_{\alpha \beta} \times \Sigma} U_\beta$ and $U_\beta \cup_{v_{\beta \delta} \times \Sigma} (e_\delta \times \Sigma) \cup_{v_{\delta \alpha} \times \Sigma} U_\alpha$.  The first of these boundary components gets identified with $-1 \times S^3$, and the second gets identified with $1 \times S^3$ decomposed into $U$, $N$, and the complement.  Finally, take two copies of this space and identify them along the subsets $U_\beta \cup_{v_{\beta \delta} \times \Sigma} (e_\delta \times \Sigma)$ sitting within the boundary of each by means of the homeomorphism $i$.  Once again we obtain $\Sigma(F)$ minus two small balls.

Presented in this fashion, we see that $\Sigma(F) -2 B^4$ can be described by a Heegaard triple $(\Sigma,\alpha,\beta,\gamma)$, where the $\beta$ curves are chosen so that $(\Sigma,\alpha,\beta)$ specifies the Heegaard decomposition $S^3 \cong U_\alpha \cup_\Sigma U_\beta$.  There is a great deal of flexibility in choosing $\beta$ curves that serve this purpose, and we focus on one way of doing so which depends on a choice of Kauffman state $x$ in the diagram $D$.

\begin{defin}\label{d edgecurve}
At each crossing of $D$, there is a corresponding edge $e \in E(W)$, and the handlebody $U_\beta$ is locally modelled by the cylinder $D^2 \times I$ nearby it.  We let $\beta_e$ denote the {\em edge curve} $S^1 \times *$ drawn there, and $\beta(x)$ the collection of edge curves $\beta_e$.
\end{defin}

\noindent Now, the Kauffman state $x$ gives rise to a spanning tree $T$ in the white graph $W$.  Each $\beta_e$ clearly bounds a disk $D^2 \times * \subset U_\beta$, and the curves in $\beta(x)$ are homologically independent and of the right number because $T$ is a spanning tree.  We denote the resulting Heegaard triple by $\h_1(x) = \h_1(x,D)$.  Observe that $\h_1(x)$ possesses an involution which exchanges the $\alpha$ and $\gamma$ curves and preserves $\beta(x)$; thus $Y_{\beta \gamma} \cong Y_{\alpha \beta} \cong S^3$.  In total, we have shown:

\begin{prop}\label{p Htriple}

Let $D$ be a connected, marked diagram of a link $L$, and $F$ the spanning surface for $L$ induced its black graph $B$.  For any choice of Kauffman state $x$, the Heegaard triple $\h_1(x)$ presents the space $\Sigma(F)$ with two small balls deleted.

\end{prop}

\subsection{The Heegaard diagram $\h_2(D,T)$.}\label{ss h_2}

We will ultimately be interested in calculating the Heegaard Floer homology of $\Sigma(K)$, and there is a Heegaard diagram $\h_2$ presenting $\Sigma(K)$ somewhat better suited to this purpose than $\h_1$.  Its description closely resembles the standard doubly-pointed Heegaard diagram for $K \subset S^3$ gotten from a diagram of $K$ \cite{OSaltknots} (cf. also $\S$\ref{S specseq}); in particular, the resulting Floer chain complex is generated by the Kauffman states of $D$ (Proposition \ref{p  Kauffman}).  We begin with a description of $\h_2$, making a choice of spanning tree $T \subset B$ in its construction, and prove that it presents $\Sigma(K)$.  Then we prove Proposition \ref{p  Kauffman} and use this result to show that the branched double-cover of a connected alternating link is an L-space.  

Starting with the given diagram $D$, choose a spanning tree $T$ of the black graph $B$, and let $T^*$ denote the dual tree in the white graph $W$.  We will imagine the diagram $D$ drawn simultaneously with the graphs $B$ and $W$ in the plane, so that edges of $B$ and $W$ meet in pairs at the crossings of $D$, and in particular every crossing lies on a unique edge in $T \cup T^*$.  Let $U$ denote a closed regular neighborhood of $D$ and $\Sigma$ its boundary.  Just as before, let $\alpha$ denote the collection of curves obtained on intersecting $\Sigma$ with every region of $D$ except the one marked white region.  Draw a meridional curve $\beta_p$ chosen to encircle the marked edge of $D$ nearby $p$, and at every crossing of $D$, introduce a pair of $\beta$ arcs on $\Sigma$ drawn on the top half of $\Sigma$ as in Figure~\ref{f intro}.  We will complete each pair of $\beta$ arcs to a full $\beta$ curve by connecting their endpoints with an additional pair of $\beta$ arcs on the bottom half of $\Sigma$.  At a particular crossing in $D$, let $e$ denote the edge of $T \cup T^*$ which passes through it.  Project that portion of the edge $e$ which meets $U$ down onto the bottom half of $\Sigma$.  In this way, we obtain a $\beta$ arc at every crossing of the diagram $D$.  To complete each to a $\beta$ curve, observe that on deleting $e$ from the spanning tree to which it belongs, we obtain two subtrees, one of which contains a marked vertex.  The other subtree $T_e$ has a silhouette in the bottom half of $\Sigma$ which we can connect up with the ends of the $\beta$ arc at the crossing through which $e$ passes, and this completes the $\beta$ curve at that crossing.  Observe that we can so complete each $\beta$ curve simultaneously so that no two meet.  For example, if edge $f$ lies on the subtree $T_e$, then the arc we draw for $f$ will appear nested inside the one drawn for $e$.  In short, there is a unique way, up to isotopy, to extend each $\beta$ arc to a $\beta$ curve in the bottom half of $\Sigma$ in such a way that no two meet and remain disjoint from the meridional $\beta$ curve drawn at the outset.  The result $(\Sigma, \alpha,\beta)$ is a Heegaard diagram, and we denote it $\h_2(D,T)$.

\begin{prop}\label{p h2}

The Heegaard diagram $\h_2(D,T)$ presents the space $\Sigma(K)$.

\end{prop}

As an example, Figure~\ref{f intro}(c) depicts $\h_2(D,T)$ for the marked diagram $D$ pictured in Figure~\ref{f blackwhitegraph} and the spanning tree $T$ pictured in Figure~\ref{f kauffman}.

\begin{proof}

The approach is to fully resolve the diagram $D$ into a single unknotted component $C$; then $\Sigma(K)$ and $\Sigma(C) \cong S^3$ will be related by surgery along a particular framed link $\L \subset S^3$.  We write down a Heegaard triple subordinate to $\L$ and in this way produce the desired Heegaard diagram.

\begin{figure}
\begin{center}
\includegraphics[height=6in]{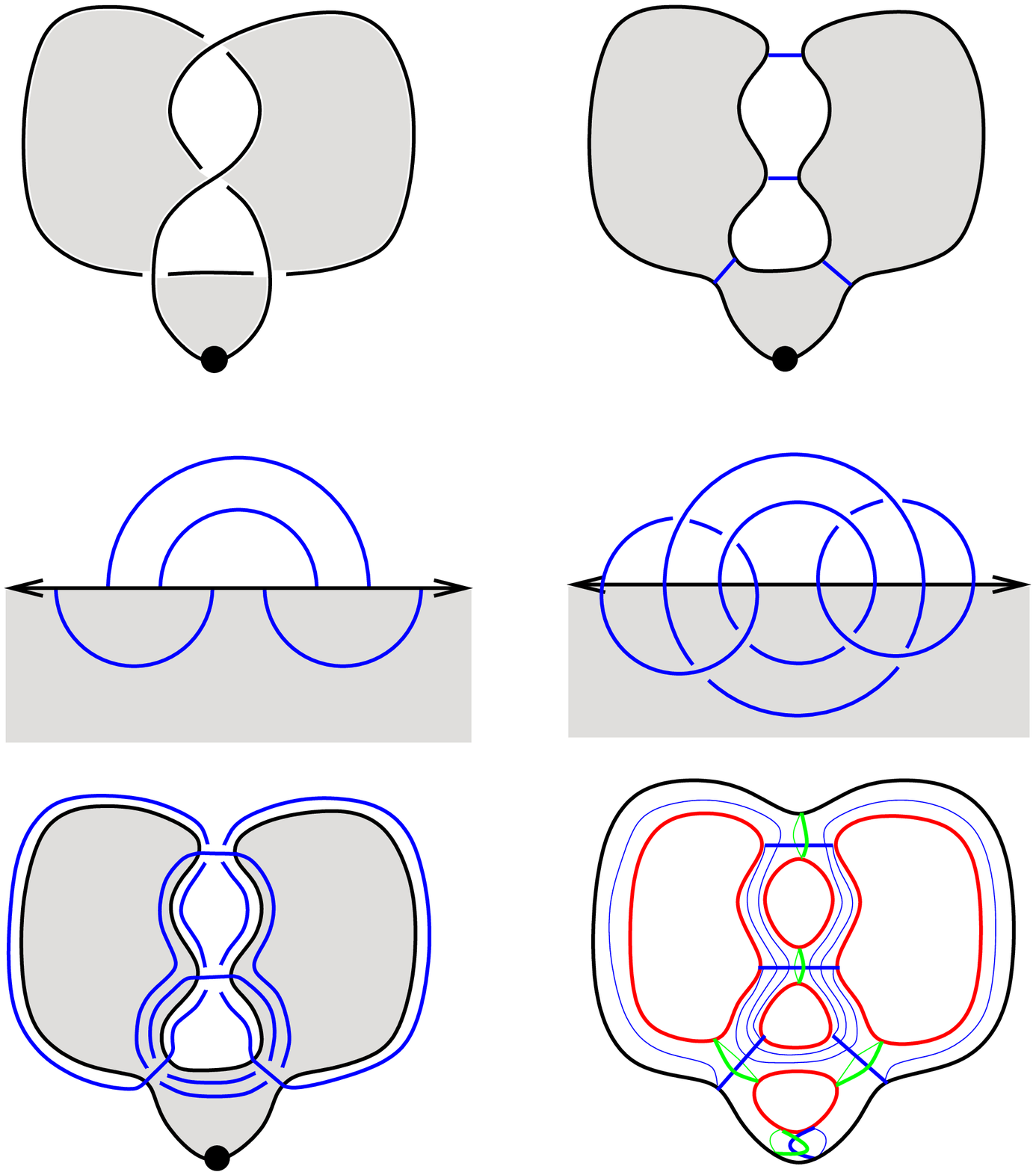} 
\put(-400,420){(a)}
\put(-200,420){(b)}
\put(-400,260){(c)}
\put(-200,260){(d)}
\put(-130,265){$+1$}
\put(-110,250){$+1$}
\put(-170,180){$+1$}
\put(-30,185){$-1$}
\put(-400,130){(e)}
\put(-200,130){(f)}
\caption{(a) The marked diagram $D$.  (b) The full resolution $C$ specified by the pair of dual trees in Figure~\ref{f kauffman}, along with the resolving arcs.  (c)  Straightening the full resolution $C$ into a line.  (d)  The framed link $\L$ obtained as the preimage of the resolving arcs in the double-cover of $S^3$ branched along $C$.  (e)  An ambient isotopy of $S^3$ which restores $C$ to its original form and carries the link $L$ with it.  (f) A Heegaard triple subordinate to the link $L$ with $\infty$-framing on each component. The blue $\gamma$ curves away from the marked point correspond to components of $L$, and the green $\beta$ curves are meridians to these components.  The $\beta$-handlebody describes the regular neighborhood of a bouquet of $L$, while the $\alpha$-handlebody describes its complement.  By performing a $(-1)^{\sigma_i}$ Dehn twist to each curve $\gamma_i$ along the meridian $\beta_i$, we obtain a Heegaard triple subordinate to the framed link $\L$.  Omitting the $\beta$ curves from this triple results in the Heegaard diagram pictured in Figure~\ref{f intro}(c).}\label{f HDderivation}
\end{center}
\end{figure}

Enumerate the crossings of $D$ by $c_1,\dots,c_n$, and center a small closed ball $B_i$ around each one.  The choice of the pair $(T,T^*)$ specifies a full resolution of the diagram $D$ into a single connected component $C$ in the plane.  Observe that all the black regions in $D$ merge to become the interior of $C$, and similarly the white regions merge to become its exterior.  At each crossing $c_i$, connect the two strands in $C \cap B_i$ by a small arc $\eta_i$ (Figure~\ref{f HDderivation}(b)). Notice that $\eta_i$ lies in the interior of $C$ precisely when the edge $e_i$ of $T \cup T^*$ which passes through $c_i$ belongs to the tree $T \subset B$; with this notation, we may take as $\eta_i$ a sub-arc of the dual edge $e_i^*$, a technical observation that will be of use soon.  Form the double-cover of $S^3$, branched along $C$, and let $L_i$ denote the preimage of $\eta_i$.  Since $C$ is unknotted, this double-cover is simply $S^3$ itself, and $L :=\cup_i L_i$ defines a link inside it.

We can describe the link $L$ concretely as follows.  Isotope $C$ to a line in the plane so that the marked point $p$ gets sent to $\infty$.  The interior and exterior of $C$ become complementary half-planes, naturally distinguished black and white, and the disjoint curves $\eta_i$ sit to either side of $C$.  Now reflect each $\eta_i$ into its complementary half-plane to obtain an immersed collection of curves in the plane.  We obtain the link $L$ by pushing the reflected image of $\eta_i$ down slightly from the plane (Figure~\ref{f HDderivation}(c)-(e)).

\begin{figure}
\begin{center}
\includegraphics[width=5in]{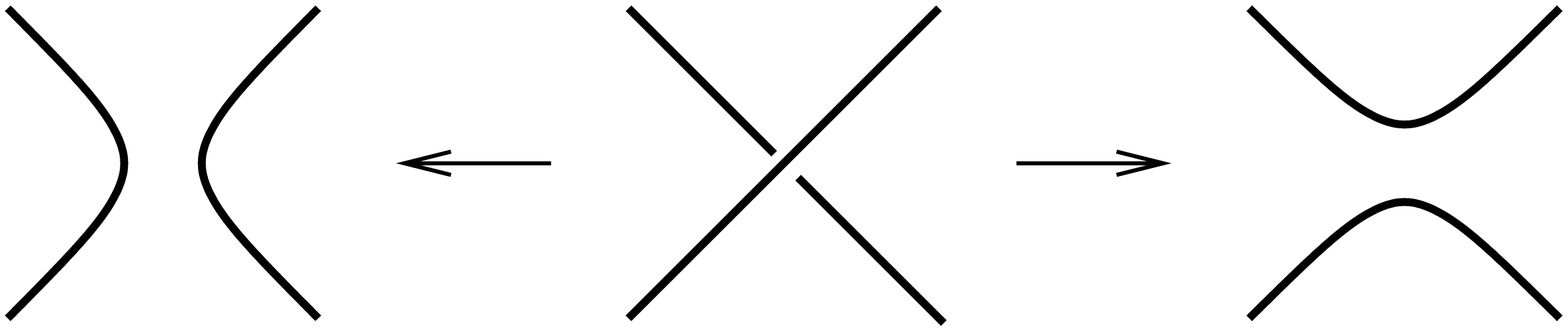} 
\put(-250,50){0}
\put(-115,50){1}
\caption{Resolutions of a crossing.}\label{f resolutions}
\end{center}
\end{figure}

Frame $L$ by placing the coefficient $(-1)^{\sigma_i}$ on component $L_i$, where $\sigma_i \in \{0,1\}$ corresponds to the type of resolution that takes place at $c_i$ (Figure~\ref{f resolutions}).  We claim that $\Sigma(K) = S^3(\L)$.  To see why this is, begin by observing that $A:= K \cap (S^3 - \cup_i B_i) = C \cap (S^3 - \cup_i B_i)$, and the double-cover of $S^3 - \cup_i B_i$ branched along $A$ is the exterior of the link $L$.  The double-cover of $B_i$ branched along either of $K \cap B_i$ or $C \cap B_i$ is a solid torus, and the preimage of $B_i$ inside $\Sigma(C) \cong S^3$ is a regular neighborhood of the link component $L_i$. Therefore, $\Sigma(K)$ is obtained from $\Sigma(C) \cong S^3$ by some surgery along $L$.  Determination of the framing follows as in \cite{OSdoublecover}, Proposition 2.1.

Next we find a Heegaard triple subordinate to $\L$, but proceed slightly differently from the construction of the Heegaard triple in $\S$\ref{ss h_1 1} (which would apply to any link which admits a diagram in which each component is an embedded, unknotted curve).  Form a regular neighborhood of the graph $C  \cup_i \eta_i$ and let $\Sigma$ denote its boundary.  We obtain the $\alpha$ curves as before, intersecting $\Sigma$ with all the unmarked faces of $C \cup_i \eta_i$, we obtain a curve $\beta_i$ chosen to encircle the arc $\eta_i$ for every $i$, and we obtain one exceptional meridian $\beta_p$ near the marked point $p$.  Next, for each $i$, push the arc $\eta_i$ to the top half of $\Sigma$, and extend its ends slightly to meet the intersection of $\Sigma$ with the dual edge $e_i^*$ (recall the definition of $\eta_i$ above).  Complete this to a curve on the surface by an arc on the bottom of $\Sigma$ chosen so as to run parallel to $C$ and avoid the exceptional $\beta$ curve.  Let $\gamma_i'$ denote the resulting curve. Observe that we can draw all these $\gamma_i'$ on $\Sigma$ simultaneously so as to remain disjoint, and the resulting collection of curves constitute an isotopic copy of the link $L$.  Lastly, let $\gamma_p$ denote a small isotopic translate of $\beta_p$, chosen to meet it in two transverse intersection points.  It is easy to see that $(\Sigma,\alpha,\beta,\gamma')$ is a Heegaard triple subordinate to the link $L$ with $\infty$-framings on each of its components (Figure~\ref{f HDderivation}(f)).  The main observation here is that each $\gamma_i'$ is a Seifert-framed longitude for the corresponding link component.  To obtain a Heegaard triple subordinate to $\L$, we perform a $(-1)^{\sigma_i}$ Dehn twist to $\gamma_i$ along the curve $\beta_i$ for every $i$.  Direct inspection shows that $(\Sigma, \alpha,\gamma)$ concurs with $\h_2(D,T)$ (Figure~\ref{f intro}(c)).

\end{proof}

Although the $\beta$ curves may appear somewhat meandering in $\h_2(D,T)$, their intersections with the $\alpha$ curves are very controlled.  Specifically, we obtain one intersection point for every incidence between an unmarked region and a crossing in the diagram $D$, as well as one special intersection point between the meridional curve $\beta_p$ and the unique $\alpha$ curve it meets.  Therefore, a collection of intersection points between these curves, one point on each curve, will consist of this special intersection point, along with one point for each unmarked region of $D$ (which lies on its associated $\alpha$ curve) and one point for each crossing of $D$ (which lies on its associated $\beta$ curve).  We have established the following result.

\begin{prop}\label{p Kauffman}

The generators of the Floer chain complex $\widehat{CF}(\h_2(D,T))$ are in one-to-one correspondence with Kauffman states of $D$.

\end{prop}

Recall that for a rational homology sphere $Y$, there is an inequality $\text{rk} \; \widehat{HF}(Y) \geq |H_1(Y)|$ (\cite{OS3mfld2}, Proposition 5.1), and $Y$ is an {\em L-space} precisely when equality holds.  Proposition \ref{p  Kauffman} leads to a quick proof of the following result, which is a special case of \cite{OSdoublecover}, Proposition 3.3.

\begin{cor}\label{c Lspace}

Let $K$ denote a non-split alternating link.  Then $\Sigma(K)$ is an L-space.

\end{cor}

\begin{proof}

The link $K$ admits a connected alternating diagram, and its number of Kauffman states equals the determinant $\det(K)$.  This value, when finite (as it is here), is well-known to equal $|H_1(\Sigma(K))|$.  The resulting sequence of inequalities $$|H_1(\Sigma(K))|  = \textup{rk}(\widehat{CF}(\h_2(D,T))) \geq \textup{rk}(\widehat{HF}(\Sigma(K))) \geq |H_1(\Sigma(K))|$$ implies the statement of the corollary.

\end{proof}

\subsection{Regions of $\h_2(D,T)$.} \label{ss regions}

Now we turn to understanding the regions of $\h_2(D,T)$.  In addition to the curve collections $\alpha$ and $\beta$ already defined, we define two other special curves in $\Sigma$.  The first $\delta_1$ is the intersection of $\Sigma$ with the marked white region of $D$.  The second $\delta_2$ is the result of pushing the boundary of a regular neighborhood of $T$ down onto the bottom of $\h_2(D,T)$.  We arrange so that $\delta$ avoids all $\beta$ curves except for $\beta_p$, which it meets in a single point.  Let $\h_2'(D,T) = (\Sigma, \alpha', \beta', \delta_1, \delta_2)$.  Although $\h_2'(D,T)$ is not an actual Heegaard diagram, we will think of it as a modified version of one. In particular, there is now an intersection point between $\alpha \cup \{ \delta_1 \}$ and $\beta$ corresponding to every incidence between a crossing and a region in $D$.  The curves in $\alpha \cup \{ \delta_1 \}$ naturally split $\h_2'(D,T)$ into its {\em top} and {\em bottom}.  Recall that for the marked diagram $D$, its {\em arcs} are the connected line segments in $D$ with each endpoint at an undercrossing or the marked point $p$.

\begin{prop}\label{p h2regions}

The regions on top of $\h_2'(D,T)$ correspond to arcs in the marked diagram $D$, and the regions on bottom correspond to regions of $D$.  All regions of $\h_2'(D,T)$ are simply-connected, and the same is true of $\h_2(D,T)$.

\end{prop}

It is instructive to verify this proposition in the example pictured in Figure~\ref{f intro}.

\begin{proof}

The assertion about the regions on top of $\h_2'(D,T)$ is immediate from inspection.  Note also the analogy to the case of the standard doubly-pointed Heegaard diagram of $K \subset S^3$ gotten from $D$ \cite{OSaltknots}; we pursue this relationship further in $\S$\ref{S specseq}.

For the second assertion, select an unmarked region $r$ of $D$.  The region $r$ abuts regions $s_1,\dots,s_k$ and crossings $c_1,\dots,c_k$ in counterclockwise order, where $r$ and the $s_i$ have opposite colors, and for every $i$, the crossing $c_i$ is incident $s_i$ and $s_{i+1}$ (subscripts $(\mod k)$).  Let $v_i$ denote the intersection point corresponding to $c_i$ and $s_i$, and $w_i$ the intersection point corresponding to $c_i$ and $s_{i+1}$.  Then the vertices $v_1,w_1,\dots,v_k,w_k$ induce a cycle $C \subset \Sigma$, whose edges alternate between $\alpha \cup \{ \delta_1 \}$ arcs (between $w_i$ and $v_{i+1}$) and $\beta$ arcs (between $v_i$ and $w_i$).  We claim that there is a region $R$ on the bottom of $\h_2'(D,T)$ corresponding to $r$ which is simply-connected and has boundary $\del R = C$.  This is apparent from the construction of $\h_2'(D,T)$.  The case of a marked region $r$ is similar.  In this case, the only change is that the region $\del R$ includes the curve $\delta_2$ and a portion of the curve $\beta_p$.  Finally, every region on the bottom of $\h_2'(D,T)$ must abut an arc of a curve in $\alpha \cup \{ \delta_1 \}$, and each such arc in turn abuts a unique region $R$ constructed in the above way; it follows that we have enumerated all the regions on the bottom of $\h_2'(D,T)$.

Next we analyze the regions of $\h_2(D,T)$.  Observe first that eliminating the curve $\delta_2$ merges the two regions of $\h_2'(D,T)$ which correspond
to the marked regions of $D$.  The resulting region is simply-connected.  Eliminating the curve $\delta_1$ causes some regions on the bottom and top of $\Sigma$ to merge.  Specifically, every region atop $\Sigma$ and incident $\delta_1$ merges with a {\em unique} region on the bottom of $\Sigma$ along the common arc of $\delta_1$ where they abut.  It follows that every region following the merger is simply-connected, establishing the claim about $\h_2(D,T)$.

\end{proof}

\subsection{Independence of $\h_2(D,T)$ on $T$.} \label{ss h_2 indep}

The union of $\alpha$ and $\beta$ curves in $\h_2(D,T)$ is an embedded graph $\Gamma \subset \Sigma$, whose vertices are the intersection points between the $\alpha$ and $\beta$ curves.  At a vertex in $\Gamma$, there is a natural cyclic ordering of the pieces of edges incident it, gotten by taking its incident arcs in counterclockwise order as they are embedded in $\Sigma$.  Thus $\Gamma$ has the structure of a {\em ribbon graph}: it is a graph with a cyclic ordering on the pieces of edges incident each vertex.  The reason for the phrase ``pieces of edges'' is that an edge may be incident a given vertex twice, and we will want to distinguish these two incidences.

\begin{lem}\label{l ribbon}

The ribbon graph $\Gamma$ so obtained is independent of the choice of spanning tree $T$.

\end{lem}

\begin{proof}

We describe the graph $\Gamma$ by constructing its vertex set, edge set, and collection of cyclic orderings without reference to $T$.  The vertex set consists of one vertex for each incidence of an unmarked region of $D$ with a crossing, as well as one for the marked point $p$.  Traverse the boundary of an unmarked region of $D$ distinct from the one bound by $\alpha_p$, and for every pair of consecutive crossings, put in an edge between the corresponding vertices.  Do the same for the region bound by $\alpha_p$, but in this case treat the marked point $p$ just like a crossing.  Next consider a crossing $c$, and enumerate its incident regions $R_1,R_2,R_3,R_4$ in counterclockwise order, so that $R_1$ and $R_2$ abut along the overstrand at $c$.  Letting $v_1,v_2,v_3,v_4$ denote the corresponding vertices, put in edges $(v_1,v_3),(v_3,v_2),(v_2,v_4),(v_4,v_1)$.  If, say, $R_4$ is the marked region, then the vertex $v_4$ does not exist, and we instead take edges $(v_1,v_3),(v_3,v_2),(v_2,v_1)$.  This collection of vertices and edges specifies the underlying graph $\Gamma$.  Orienting the $\alpha$ and $\beta$ curves induces an orientation on the edges of $\Gamma$, and additionally orienting $\Sigma$ specifies an intersection number $\pm 1$ at each vertex of $\Gamma$.  Fix a vertex $v$ in $\Gamma$, and let $a^+,a^-,b^+,b^-$ denote the pieces of edges incident $v$, so that $a$ edges correspond to $\alpha$ arcs and $b$ edges to $\beta$ arcs, and the symbols $+, -$ indicate whether the piece of edge is directed out of or into $v$.  If the intersection number at $v$ is $+1$, we take the cyclic ordering $(a^+,b^+,a^-,b^-)$, and if it is $-1$ we take $(b^+,a^+,b^-,a^-)$ instead.  It is straightforward to check that the resulting ribbon graph $\Gamma$ agrees with the one given by $\h_2(D,T)$, regardless of the choice of $T$.

\end{proof}

\begin{prop}\label{p h2indep}

The diagram $\h_2(D,T)$ is independent, up to diffeomorphism, of the choice of spanning tree $T$.

\end{prop}

Here a {\em diffeomorphism of Heegaard diagrams} $(\Sigma, \alpha, \beta)$ is a diffeomorphism $f: \Sigma \to \Sigma'$ carrying each $\alpha_i$ onto $\alpha_i'$ and $\beta_j$ onto $\beta_j'$.

\begin{proof}

A ribbon graph $\Gamma$ uniquely specifies an embedding of its underlying graph into a closed surface $S$, up to diffeomorphism of $S$, in the following standard way.  Each vertex gives rise to a $0$-handle, and each edge gives rise to a $1$-handle, which we attach to the $0$-handles according to the cyclic ordering of the edges at each vertex.  The graph $\Gamma$ is embedded in the union of these $0$- and $1$- handles in the obvious manner (the ribbon-like surface obtained so far is responsible for the terminology ``ribbon graph'').  Finally, attach a $2$-handle by the identity map to each boundary component of the ribbon surface to produce the desired embedding $\Gamma \subset S$.

Now let $\Gamma$ denote the ribbon graph obtained from the Heegaard diagram $\h_2(D,T)$, with each edge decorated by the label of the curve in $\alpha \cup \beta$ of which it is an arc.  By Lemma \ref{l ribbon}, $\Gamma$ is independent of $T$.  The embedding $\Gamma \subset S$ from the preceding paragraph agrees with the original embedding $\Gamma \subset \Sigma$, up to diffeomorphism, because every region of $\h_2(D,T)$ is simply-connected (Proposition \ref{p h2regions}).  Now from the embedding $\Gamma \subset \Sigma$ and the edge decorations we can recover the $\alpha$ and $\beta$ curves, and so reconstruct $\h_2(D,T)$ up to diffeomorphism.

\end{proof}

Following Proposition \ref{p h2indep}, we denote the Heegaard diagram $\h_2(D,T)$ by $\h_2(D)$, or just $\h_2$ when $D$ is understood.  We additionally make $\h_2$ into a pointed Heegaard diagram by placing a basespoint $z$ nearby the marked point $p$ in the knot diagram.  Thus $z$ lies to one side of the meridional curve $\beta_p$; by Proposition \ref{p h2regions}, the regions just to either side of $\beta_p$ are the same, so the region of $\h_2$ in which $z$ is chosen is well-defined.  We make $\h_1$ into a pointed diagram in an analogous manner.

\subsection{The relationship between $\h_1$ and $\h_2$.}\label{ss heegaardmoves}

According to a theorem of Singer, the two diagrams $\h_1$ and $\h_2$ are related by a sequence of Heegaard moves: isotopies, handleslides, and (de-)stabilizations.  In this section, we explicitly describe such a sequence from $\h_2$ to $\h_1$.  This description will enable us to compare the absolute gradings of the generators of $\widehat{CF}(\h_1)$ and $\widehat{CF}(\h_2)$, a key to establishing the absolute gradings formula in $\S$\ref{S absgradings}.

\begin{figure}
\begin{center}
\includegraphics[width=5.5in]{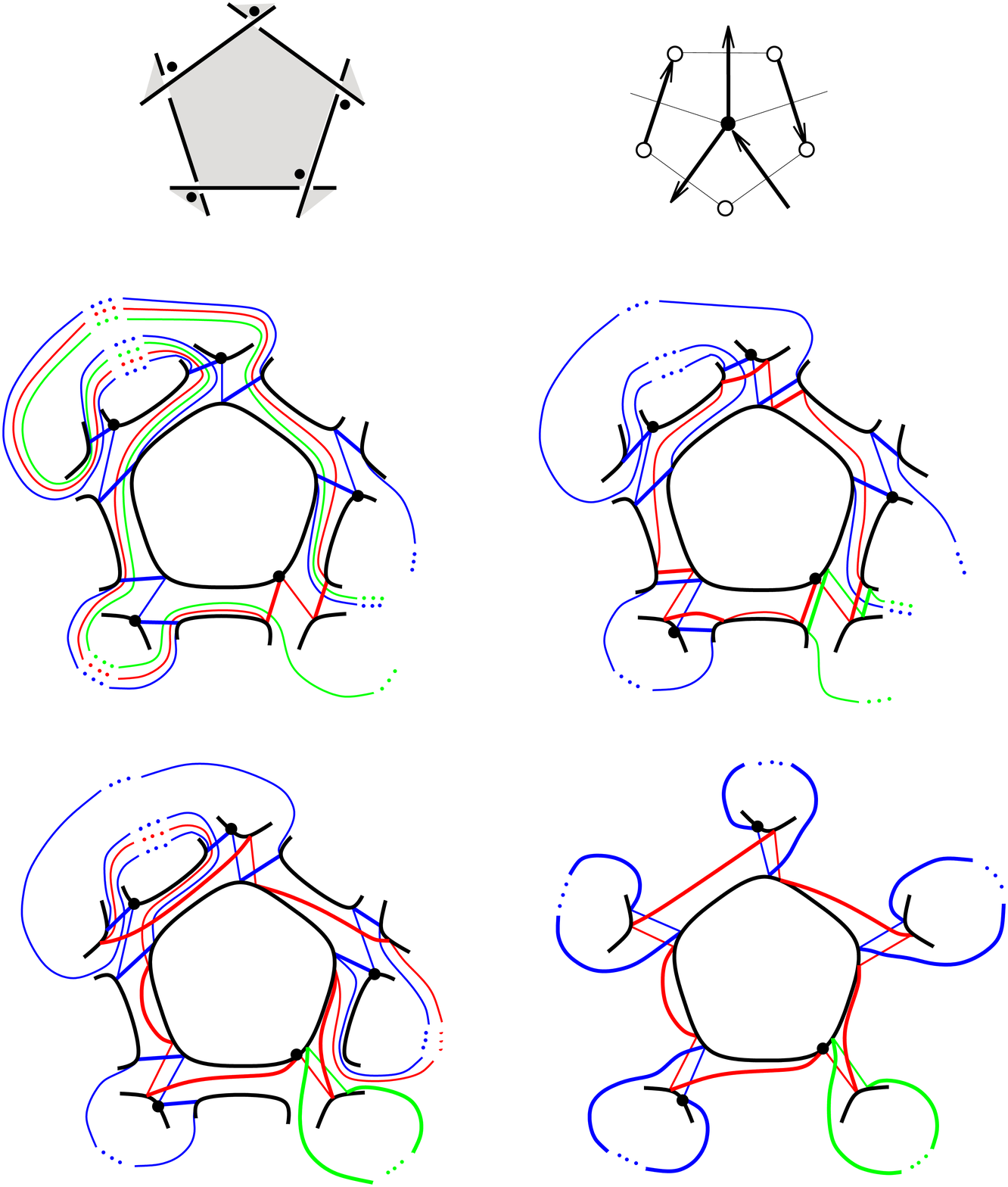}
\put(-370,450){$(a)$}
\put(-170,450){$(b)$}
\put(-400,350){$(c)$}
\put(-200,350){$(d)$}
\put(-400,150){$(e)$}
\put(-200,150){$(f)$}
\put(-123,430){$v_i$}
\put(-85,390){$e_i$}
\put(-272,212){\textcolor{red}{$\beta_i$}}
\caption{(a) A region in a knot diagram $D$, with part of a Kauffman state $x$.  (b) The portion of the black and white graphs nearby the pictured region.  The edges in $T = T_B(x)$ and $T_W(x)$ are thickened and directed accordingly. (c) The corresponding local picture of $\h_2(D)$.  All the $\alpha$ curves pictured are black.  The curve $\beta_i$ is shown in red, and the curves $\beta_j$ where $v_j$ neighbors $v_i$ are shown in blue.  The arc shown in green represents a number of nested parallel arcs of the $\beta_k$, where $e_k \ne e_i$ is an edge along the path from $r$ to $v_i$ in $T$. (d) The effect of handlesliding all the curves $\beta_j$ with $j \leq i$.  Now the green curve corresponds to a unique $\beta_j$.  (e) The effect of another sequence of handleslides.  (f)  Following the destabilizations, we obtain the local picture of $\h_1(D)$ shown.  The Kauffman generator $x \in \widehat{CF}(\h_2(D))$ gets carried onto the corresponding generator in $\widehat{CF}(\h_1(D))$ through this sequence of moves.  Local pieces of both are pictured by the heavy black dots in (c)-(f).}
\label{f Heegaardmoves}
\end{center}
\end{figure}

Begin by choosing a spanning tree $T$ of the black graph, let $x = x(T)$ denote the corresponding Kauffman state, and construct the diagram $\h_2(D,T)$ (of course, Proposition \ref{p h2indep} removes the dependence of $\h_2$ on $T$, but making this choice is useful for descibing what follows).  Orient $T$ out of the root $r_B$ ($\S$\ref{ss prelim_3}).  Now enumerate its vertices $r_B,v_1,\dots,v_k$ in such a way that for every directed edge $(v_i,v_j) \in E(T)$, we have $i < j$.  To each vertex $v_i$ there corresponds a curve $\beta_i$ which is associated to the unique edge $e_i$ in $T$ directed into $v_i$.  This is pictured in Figure~\ref{f Heegaardmoves}(c).

For $i=1,\dots,k$ in turn, handleslide the curve $\beta_i$ over all $\beta_j$ for which $(v_i,v_j) \in E(T)$.  We can then perform an isotopy to the transformed $\beta_i$ so that it is supported on the bottom of the Heegaard surface away from crossings in the diagram.  Having done so, each $\beta_i$ curve is supported nearby the curve $\alpha_i$ corresponding to $v_i$ (Figure~\ref{f Heegaardmoves}(d)).  Next, handleslide $\beta_i$ over every $\beta$ curve incident $\alpha_i$ which corresponds to an edge of the dual tree $T^*$.  Following an isotopy, the resulting $\beta_i$ is pictured as in Figure~\ref{f Heegaardmoves}(e).

Observe that at this stage, each $\beta_i$ only meets $\alpha$ curves which correspond to black regions, and meets them with sign and order exactly as in the diagram $\h_1$.  To actually obtain $\h_1$, we perform a sequence of destabilizations.  Thus, enumerate the vertices of $T^*$ by $w^1,\dots,w^\ell$ in such a way that $w^i$ is a leaf of the tree $T^*-w^1-\cdots-w^{i-1}$ for every $i$.  Adopt notation $\alpha^i$ and $\beta^i$ analogous to $\alpha_i$ and $\beta_i$ above.  Now observe that for every $i$, $\alpha^i$ and $\beta^i$ meet in a single intersection point, and that any other $\beta$ curve which meets $\alpha^i$ takes the form $\beta^j$ with $j > i$.  It follows that for $i=1,\dots,\ell$ in turn, we can perform handleslides and isotopies to $\beta^i$ so as to eliminate all its intersection points except the single one with $\alpha^i$, and then destabilize the pair $(\alpha^i,\beta^i)$.  In effect, this amounts to performing a $2$-surgery along each $\alpha^i$ in $\Sigma$ and erasing all the curves $\alpha^i$ and $\beta^i$ from the resulting surface.  The end result is the desired Heegaard diagram $\h_1$ (Figure~\ref{f Heegaardmoves}(f)).   (A rigorous argument as to why would go as follows: first, the $\alpha$ and $\beta$ curves in both this diagram and $\h_1$ give rise to isomorphic ribbon graphs; second, the regions in $\h_1$ are simply-connected; and third, the proof of Proposition \ref{p h2indep} goes through to show that the two Heegaard diagrams are diffeomorphic.  We leave further details to the reader.)

Having described the transformation from $\h_2$ into $\h_1$, we take a closer look at the generators of $\widehat{CF}(\h_1)$.  The $\alpha$ and $\beta$ curves in $\h_1$ come in pairs $(\alpha_w, \beta_w)$, $w \in V(\widetilde{W})$.  Intersection points in $\alpha_w \cap \beta_w$ correspond to edges incident the vertex $w$ in $W$.  Let $x$ denote a Kauffman state for $\D$, $T_W(x)$ the corresponding spanning tree of $W$, and orient $T_W(x)$ out of its root $r_W$.  For each vertex $w \in V(\widetilde{W})$, there is a unique edge in $T_W(x)$ directed into it, and we let $x_w$ denote the intersection point in $\alpha_w \cap \beta_w$ which corresponds to this edge.  The tuple $(x_w)$ is an intersection point in $\T_\alpha \cap \T_\beta$, and we call it the {\em Kauffman generator of $\widehat{CF}(\h_1)$} corresponding to $x$.  Here is an alternative description of this generator.  Regard $x$ as a choice of corners in $D$, so that exactly one appears at each crossing, and exactly one appears in each region except for the two regions incident the marked point $p$.  In the Heegaard diagram $\h_1$, at each corner that appears in a white region $w$, there is an intersection point between $\alpha_w$ and $\beta_w$.  The collection of these intersection points is the Kauffman generator corresponding to $x$ in $\widehat{CF}(\h_1)$.  For example, there are seven generators for the Heegaard diagram $\h_1(D)$ pictured in Figure~\ref{f HD1}, and five of these are Kauffman generators.

Now let $x = x(T)$ denote the Kauffman state corresponding to $T$, identified as an collection of intersection points between $\alpha$ and $\beta$ curves in $\h_2$. Every Heegaard move in the above transformation from $\h_2$ to $\h_1$ takes place in the complement of $x$ and the basepoint $z$, with the exception of the destabilizations, each of which deletes a point from $x$.  It follows that the absolute grading of $x \in \widehat{CF}(\h_2)$ equals the absolute grading of the resulting generator in $\widehat{CF}(\h_1)$, which in turn is precisely the Kauffman generator corresponding to $x$.  Since the choice of spanning tree $T$ and hence the Kauffman state $x$ here was arbitrary, we have established the following fact.

\begin{prop}\label{p gradingpreserving}

There is a sequence of Heegaard moves taking $\h_2$ to $\h_1$ which carries the generator of $\widehat{CF}(\h_2)$ corresponding to the Kauffman state $x$ onto the corresponding generator of $\widehat{CF}(\h_1)$.  The absolute gradings of these two generators are equal in their respective complexes.

\end{prop}

Notice that the sequence of moves from $\h_2$ to $\h_1$ given above is adapted to the choice of $x$.  Presumably it is the case that any one such sequence matches up the Kauffman generators between $\widehat{CF}(\h_1)$ and $\widehat{CF}(\h_2)$, and could be used to show that matched-up generators have the same absolute grading. Instead we have opted for a somewhat indirect approach, leaning on Proposition \ref{p h2indep}.


\vspace{.5in}

\section{Absolute Gradings and $Spin^c$ Structures.}\label{S absgradings}

Let $x$ denote a Kauffman state of the marked diagram $D$, identified as a genenerator of $\widehat{CF}(\h_2(D))$.  Our goal in this section is to give an algebraic-combinatorial formula for the absolute grading $\text{gr}(x)$ as well as its associated $spin^c$ structure $\t(x)$.  We first state the expression for $\text{gr}(x)$, then give an outline of its derivation, and finally provide a detailed proof in $\S\S$\ref{ss triangle}-\ref{ss chern}.  In $\S$\ref{ss spinc} we deduce the formula for $\t(x)$, and in $\S$\ref{ss alternating} we investigate the absolute gradings formula in the case of a non-split alternating link.

\subsection{Overview of the main result.}\label{ss absgradings1}

We begin by positing an important assumption that will remain in place in the remainder of the paper.

\begin{center}

{\em The link $K$ has non-zero determinant: $\det(K) \ne 0$.}

\end{center}

\noindent This assumption guarantees that $\Sigma(K)$ is a rational homology sphere, and that the Goeritz forms and coloring matrix ($\S$\ref{S prelim}) are invertible.

Given a Kauffman state $x$, we induce an orientation on the white graph $W$ in the following way.  Given an edge $e \in E(W)$, consider the crossing $c$ to which it corresponds, as well as the white region which abuts $c$ to the same side of the overstrand as $x(c)$.  We direct $e$ to point towards the endpoint corresponding to this white region.  At a vertex $w \in V(W)$, we compute the {\em signed degree} $d_x(w)$ as the number of edges directed out of $w$ minus the number of edges directed into $w$, with respect to this orientation on $W$.  Number the unmarked vertices $w_1,\dots,w_m$, and define the {\em degrees vector} $$v^W_x := (d_x(w_1),\dots,d_x(w_m))^T$$ as well as the Goeritz matrix $G_W$.  Let $q(v)$ denote the quadratic form $v^TG_W^{-1}v$. Finally, define $\delta(x)$ as the number of edges $e$ in $T_W(x)$ with $\mu(e) = -1$.
\begin{thm}\label{t grading}

For a Kauffman state $x$ in a marked diagram $D$, the absolute grading of the corresponding generator in $\widehat{CF}(\h_2(D))$ is given by the formula:

\begin{equation}\label{e absgrading}
\textup{gr}(x) = \delta(x) + { q(v^W_x) -2m - 3 \sigma(G_W) \over 4}.
\end{equation}

\end{thm}

The proof of Theorem \ref{t grading} goes as follows.  By Proposition \ref{p gradingpreserving}, we may calculate $\text{gr}(x)$ by regarding $x$ as a Kauffman generator of the Heegaard diagram $\h_1$.  This diagram fits into the Heegaard triple $\h_1(x)$ subordinate to the link $\L$ (defined at the end of $\S$\ref{ss h_1 2}), where $\Sigma(K) \cong S^3(\L)$.  Proposition \ref{p Htriple} identifies the $4$-manifold $X$ presented by $\h_1(x)$, and in particular $X$ is independent of $x$.  Nevertheless, in order to calculate $\text{gr}(x)$, we use the specific $\h_1(x)$ adapted to both $x$ and $X$.  We invoke the following formula (\cite{OS4mfld}, p. 385, Equation (12)):
\begin{equation}\label{e gradingformula}
\text{gr}(x) = - \mu(\psi_x) + 2n_z(\psi_x) + {c_1(\s_z(\psi_x))^2- 2\chi(X)- 3 \sigma(X) \over 4}.
\end{equation}  Here $\psi_x$ is a specific homotopy class of Whitney triangles whose domain in $\h_1(x)$ consists of a disjoint union of embedded triangles.  Now the task is to identify the terms in Equation (\ref{e gradingformula}) with the corresponding terms in Equation (\ref{e absgrading}).  We  easily compute the Maslov index and intersection number:
\begin{equation}\label{e Maslov}
\mu(\psi_x) = -\delta(x) \quad \mbox{and} \quad n_z(\psi_x)=0.
\end{equation} For the remaining terms, we must first understand the intersection pairing on $H_2(X)$.  A basis for $H_2(X)$ is given in the following standard way.  Each vertex $w \in V(\widetilde{W})$ corresponds to an unknotted link component $K_w \subset \L$. Take a disk spanning $K_w \subset S^3$, push its interior into the $4$-ball bound by $S^3$, and cap off with the core of the handle attachment along $K_w$.  The result is a sphere $\Sigma_w \subset X$, and the collection of the classes $[\Sigma_w]$ freely generate $H_2(X)$.
With respect to the Poincar\'e duals of these classes, the intersection form on $H^2(X)$ is given by the matrix $G_W^{-1}$. In particular, $$\chi(X) = \text{rk}(G_W) = m \quad \mbox{and} \quad \sigma(X) = \sigma(G_W).$$

\noindent It remains to calculate $c_1(\s_z(\psi_x))$.  This amounts to computing $ \langle c_1(\s_z(\psi_x)), [\Sigma_w] \rangle$ for each vertex $w$.  We write down a triply-periodic domain $\P_w$ in $\h_1(x)$ representing the homology class $[\Sigma_w]$ and apply the first Chern class formula (\cite{OS4mfld}, Prop. 6.3):
\begin{equation}\label{e chernformula}
\langle c_1(\s_z(\psi_x)),[\Sigma_w] \rangle = -2n_z(\P_w) + \#(\del \P_w) + \widehat{\chi}(\P_w) + 2 \sigma(\psi_x,\P_w).
\end{equation}  The bulk of the proof involves showing that the right-hand side of Equation (\ref{e chernformula}) reduces to $-d_x(w)$.  Thus we obtain the expression
\begin{equation}\label{e chernvector}
c_1(\s_z(\psi_x)) = v_x^W
\end{equation} with respect to the basis of Poincar\'e duals to the $[\Sigma_w]$, and so $c_1(\s_z(\psi_x))^2 = (v_x^W)^T G_W^{-1} v_x^W =q(v_x^W)$. This completes the identification of terms between Equations (\ref{e absgrading}) and (\ref{e gradingformula}) and hence the proof of Theorem \ref{t grading}.

We make two remarks regarding Theorem \ref{t grading}.  First, $\text{gr}(x)$ has an equivalent expression to that in Equation (\ref{e absgrading}) with reference to the black graph $B$ in place of $W$: one defines the degrees vector $v^B_x$ analogously to $v^W_x$, and sets $\delta(x)$ instead equal to the number of edges $e$ in $T_B(x)$ with $\mu(e) = +1$.  Second, it is possible to obtain an expression for $\text{gr}(x)$ by working solely with the diagram $\h_2$, proceeding along the lines sketched above, without ever making reference to $\h_1$.  However, in order to obtain the given expression involving the Goeritz form, the passage to $\h_1$ seems the easiest route to take.

\subsection{The triangle $\psi_x: \Delta \to \h_1(x)$.}\label{ss triangle}

Write $\h_1(x) = (\Sigma,\alpha,\beta(x),\gamma,z)$, where $Y_{\alpha \gamma} \cong \Sigma(K)$ and $\beta(x)$ is as in Definition~\ref{d edgecurve}.  Enumerate the vertices $w_0,\dots,w_n$ of $W$ so that $w_0 = r_W$ and each $w_j$ neighbors some $w_i$, $i< j$ in $T_W(x)$, for all $j > 0$.  For $j > 0$, let $\alpha_j$ and $\gamma_j$ denote the curves corresponding to vertex $w_j$.  Recall ($\S$\ref{ss h_1 2}) that each curve in $\beta(x)$ corresponds to an edge of $T_W(x)$; we let $\beta_j$ denote the one corresponding to the edge of $T_W(x)$ that directs into the vertex $w_j$, when $T_W(x)$ is oriented out of its root.  Since $T_W(x)$ is a tree, an easy induction shows that there is a unique intersection point between $\{ \beta_1,\dots,\beta_j \}$ and $\{ \alpha_1, \dots, \alpha_j \}$ for every $j$, and it lies in $(\beta_1 \cap \alpha_1) \times \cdots \times (\beta_j \cap \alpha_j)$.  The same reasoning applies to the $\gamma$ curves in place of the $\alpha$.  Let $x_0 \in \T_\alpha \cap \T_\beta$ and $x_1 \in \T_\beta \cap \T_\gamma$ denote the resulting unique intersection points.

\begin {figure}
\begin {center}
\includegraphics[width=2in]{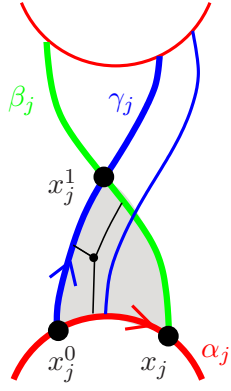}
\put(-60,5){$x_j$}
\put(-95,72){$x_j^1$}
\put(-95,5){$x_j^0$}
\put(-37,10){\textcolor{red}{$\alpha_j$}}
\put(-110,105){\textcolor{green}{$\beta_j$}}
\put(-72,105){\textcolor{blue}{$\gamma_j$}}
\caption{The triangular domain $\Delta_j$ (shaded) at a crossing $c$ with $\mu(c)$ = -1.  A dual spider is pictured, and $\alpha_j$ and $\gamma_j$ have been oriented as the boundary of $\P_{w_j}$ ($\S$\ref{ss periodic}).  The orientation on $\beta_j$ will depend on the relative values of $a_{e_j}$ and $b_{e_j}$.  This local picture appears on the bottom of the Heegaard surface, viewed from below. (Contrast this viewpoint with that of Figures~\ref{f periodic} and \ref{f periodicmodify}, in which the top half of the Heegaard surface has been sliced away.)}\label{f smalltriangle}
\end {center}
\end {figure}

Now fix a vertex $w_j$ and let $\{ x^0_j \} = \alpha_j \cap \beta_j$, $\{ x^1_j \} = \beta_j \cap \gamma_j$, and $x_j = x_{w_j}$ ($\S$\ref{ss heegaardmoves}).  There is a small triangular domain $\Delta_j$ connecting $x^0_j, x^1_j$, and $x_j$ in $\h_1(x)$ which will appear as shown in Figure~\ref{f smalltriangle} if $\mu(c) = -1$, and will appeared reflected from this picture if $\mu(c) = +1$.  We let $\psi_x \in \pi_2(x_0,x_1,x)$ denote the homotopy class whose domain is the union of these small triangles.  Clearly, $n_z(\psi_x) = 0$.

\begin {figure}
\begin {center}
\includegraphics[width=1in]{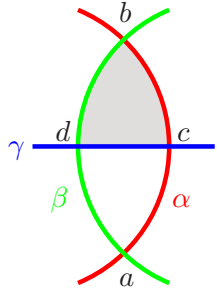}
\put(-38,0){$a$}
\put(-38,100){$b$}
\put(-16,55){$c$}
\put(-62,55){$d$}
\put(-18,30){\textcolor{red}{$\alpha$}}
\put(-64,30){\textcolor{green}{$\beta$}}
\put(-80,50){\textcolor{blue}{$\gamma$}}
\caption{The calculation $\mu(\Delta_j) = -1$.}\label{f maslovtriangle}
\end {center}
\end {figure}

Next we determine the Maslov index $\mu(\psi_x)$.  This is calculated as the sum of contributions from the individual $\Delta_j$.  For a domain $\Delta_j$ which appears reflected from Figure~\ref{f smalltriangle}, the Maslov index contribution is $\mu(\Delta_j) = 0$ (\cite{OS3mfld2}, Prop 9.5, or \cite{sucharit}, Section 2 and Theorem 4.1).  The contribution $\mu(\Delta_j)=-1$ for the kind pictured in Figure~\ref{f smalltriangle} follows as well from \cite{sucharit}.  Alternatively, its contribution can be gleaned from Figure~\ref{f maslovtriangle}.  If $\varphi_1 \in \pi_2(a,b)$ has domain equal to the bigon in this picture, and $\varphi_2 \in \pi_2(b,c,d)$, $\varphi_3 \in \pi_2(a,c,d)$ have domains equal to the two triangles, then $\mu(\varphi_1) = 1, \mu(\varphi_2) = \mu(\Delta_j)$, and $\mu(\varphi_3) = 0$ on the one hand, while additivity implies $\mu(\varphi_3) = \mu(\varphi_1 * \varphi_2) = \mu(\varphi_1) + \mu(\varphi_2)$.  The result $\mu(\Delta_j)=-1$ in this case follows.  (Thanks to Zolt\'an Szab\'o for providing this argument.) Now the sum of the $\mu(\Delta_j)$ results in $\mu(\psi_x) = - \delta(x)$, which is the expression (\ref{e Maslov}).

\subsection{The periodic class in $\h_1(x)$ representing $[\Sigma_w]$.}\label{ss periodic}

We revisit the description of the homology class $[\Sigma_w] \in H_2(X)$ ($\S$\ref{ss absgradings1}).  Express $X$ as the union $(Y_{\alpha \beta} \times I)  \cup_{U_\beta} (Y_{\beta \gamma} \times I)$, where we decompose $Y_{\alpha \beta} \times \{ 1 \} = U_\alpha \cup_\Sigma U_\beta$ and $Y_{\beta \gamma} \times \{ 1 \} = U_\gamma \cup_\Sigma U_\beta$.  This decomposition distinguishes the surface $\Sigma = \del U_\beta \subset X$.  Now fix a vertex $w \in V(\widetilde{W})$ and consider the curve $\gamma_w \subset \Sigma$.  It bounds a disk inside $Y_{\beta \gamma} \times I$ on the one hand, as well as one inside $Y_{\alpha \beta} \times I$ since $[\gamma_w] = [\alpha_w] + \sum_e c_e [\beta_e]$, for suitable $c_e$.  The union of these two disks represents the homology class $[\Sigma_w]$.  It follows that a periodic class $\P$ in $\h_1(x)$ will represent $[\Sigma_w]$ provided $\del \P = \gamma_w - \alpha_w - \sum_e c_e \beta_e$, for suitable $c_e$.

By way of this latter property, we can easily construct a periodic class $\P_w$ which represents $[\Sigma_w] \in H_2(X)$ and also satisfies $n_z(\P_w) = 0$.  Any point $y \in \Sigma - \beta(x)$ (cf. Definition~\ref{d edgecurve}) is connected to the basepoint $z$ via some path $P$ inside this subsurface which meets the $1$-cycle $\alpha_w \cup -\gamma_w$ with some algebraic multiplicity, having oriented the path from $z$ to $y$ and $\alpha_w$ and $\gamma_w$ counterclockwise.  This multiplicity is the coefficient of $\P_w$ on the region containing $y$.  It follows at once that $\del \P_w$ takes the required form and so $\P_w$ is the desired periodic class.  Notice, in particular, that $\P_w$ does not depend on the other $\alpha$ and $\gamma$ curves in its construction.

It is convenient to encode the periodic class $\P_w$ by a planar graph $H_w$ with labels on its faces.  The graph consists of a copy of the tree $T = T_W(x)$ along with the collection of edges $\Gamma(w)$ incident the vertex $w$.  There are two copies of an edge $e \in \Gamma(w) \cap E(T)$, and we distinguish the copy in $\Gamma(w)$ by pushing its interior clockwise off of the copy in $T$ when $\mu(e) = +1$, and counterclockwise when $\mu(e) = -1$.  We orient an edge in $\Gamma(w)$ out of $w$ when $\mu(e) = +1$ and into it when $\mu(e) = -1$.  The resulting planar graph has some faces, which we label by analogy to the construction of $\P_w$: for a point $y$ inside a face, we orient a path from the basepoint $z$ to it that avoids the edges of $T$, and record the oriented intersection number with the oriented edges in $\Gamma(w)$.

\begin {figure}
\begin {center}
\includegraphics[width=4in]{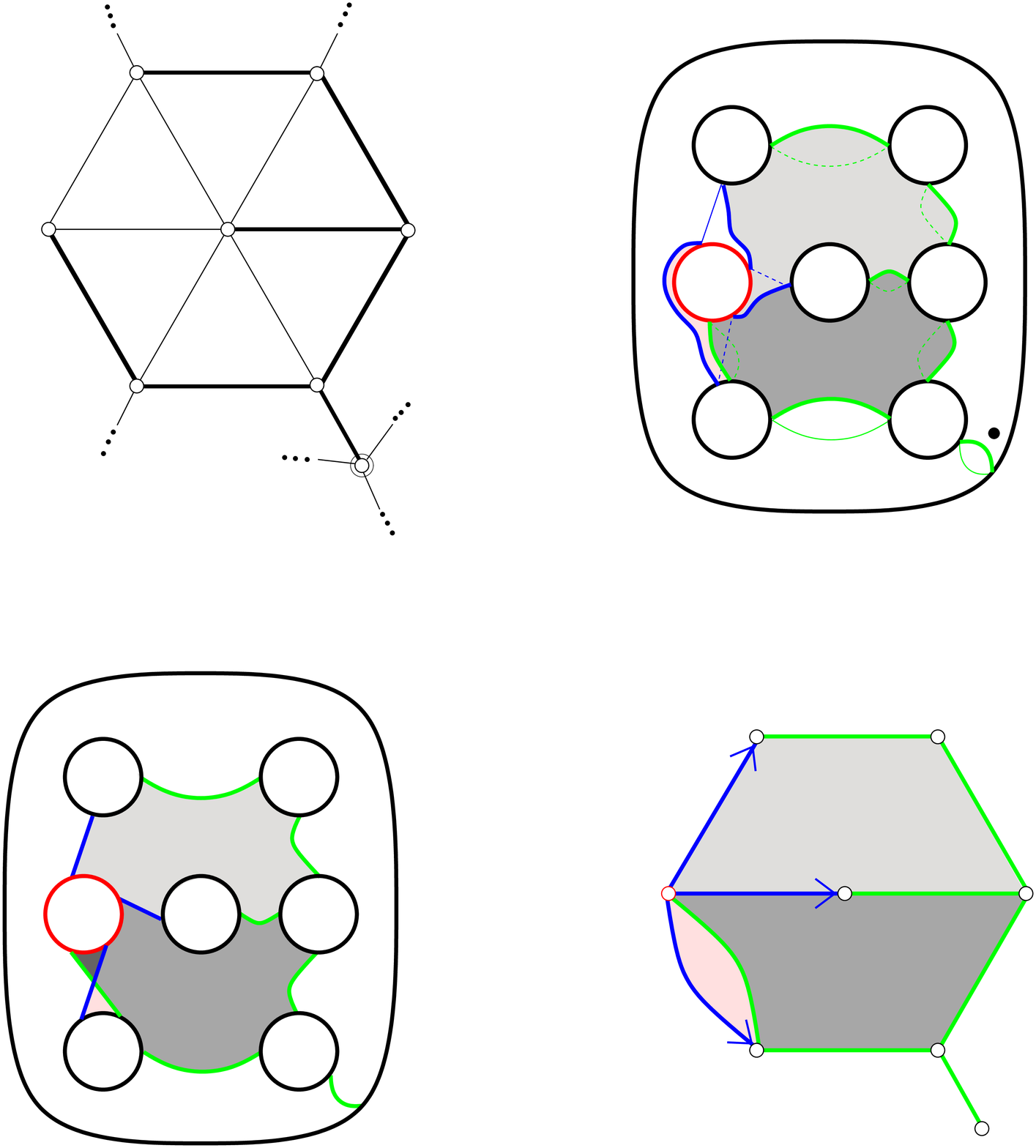}
\put(-320,320){(a)}
\put(-140,320){(b)}
\put(-320,130){(c)}
\put(-140,130){(d)}
\put(-290,260){$w$}
\put(-15,205){$z$}
\put(-105,95){\textcolor{blue}{$f_1$}}
\put(-80,75){\textcolor{blue}{$f_2$}}
\put(-110,35){\textcolor{blue}{$f_3$}}
\caption{(a) The white graph corresponding to a knot diagram, with a root marked, vertex $w$ identified, and spanning tree highlighted.  (b)  The Heegaard surface and the top of the periodic domain $\P_w$. (c)  The bottom of the periodic domain $\P_w$.  The top half of the Heegaard surface has been sliced off in this view.  (d) The affiliated face-labeled planar graph $H_w$.  In (b)-(d), a light gray region has coefficient $+1$; medium gray, $+2$; dark gray, $+3$; and pink, $-1$.}\label{f periodic}
\end {center}
\end {figure}

We extract a couple of pieces of information from the graph $H_w$.  For every edge $f \in \Gamma(w)$, consider the subgraph of $H_w$ obtained on deleting all the other edges in $\Gamma(w)$.  The resulting graph has a unique cycle $C(T,f)$ made up of $f$ and some edges of $T$, and it gets an orientation from the one on $f$.  We define a label $$\sigma(f) = \pm 1 \text{ according as } C(T,f) \text{ is oriented clockwise or counterclockwise.}$$   Observe that for every $f \in \Gamma(w) \cap E(T)$, we have $\sigma(f) = -1$.  As an example, consider the edge $f_3$ pictured in Figure~\ref{f periodic}(d), and also consider how that picture would change if the crossing corresponding to $f_3$ were reversed.  In general, all but one such $f$ is oriented into $w$ in the orientation on $W$ induced by $T$ ($\S$\ref{ss absgradings1}).  Similarly, for $f \in \Gamma(w) \setminus E(T)$, one checks that $\sigma(f) = \pm 1$ according as $f$ is directed out of or into $w$ in that orientation.  Consequently,

\begin{equation}\label{e dx}
d_x(w) = 2 + \sum_{f \in \Gamma(w)} \sigma(f) = 2 -|\Gamma(w) \cap E(T)| + \sum_{f \in \Gamma(w) \setminus E(T)} \sigma(f).
\end{equation}  Here the ``$2$" corrects for the unique edge in $\Gamma(w) \cap E(T)$ which is directed out of $w$ in $W$.  We also remark that the label on a face $F$ of $H_w$ is equal to the sum of $\sigma(f)$, over all edges $f \in \Gamma(w)$ for which $F$ is contained inside $C(T,f)$. Lastly, for every edge $e \in T \subset H_w$, let $(a_e,b_e)$ denote the pair of labels on the faces abutting along $e$.

\subsection{The first Chern class formula.}\label{ss chern}

In this section we calculate the terms appearing in the first Chern class formula (\ref{e chernformula}) to deduce the expression (\ref{e chernvector}).

\begin{prop}\label{p chernprop}

For the triangle $\psi_x$ and periodic domain $\P_w$ above, we compute
\begin{enumerate}

\item $n_z(\P_w) = 0$,
\item $\#(\del \P_w) = 2 + \sum_{e \in E(T)} |a_e - b_e|$,
\item $\widehat{\chi}(\P_w) = \sum_{f \in \Gamma(w) \setminus E(T)} \sigma(f) +|\Gamma(w) \cap E(T)| - \sum_{e \in E(T)} (a_e + b_e)$, and
\item $\sigma(\psi_x,\P_w) = -|\Gamma(w) \cap E(T)| + \sum_{e \in E(T)} \min \{a_e,b_e \}$.

\end{enumerate}

\end{prop}

Combining the terms in Proposition \ref{p chernprop} into (\ref{e chernformula}), and simplifying the result using (\ref{e dx}) and the ``miraculous cancelation" $|a-b|-(a+b)+2\min\{a,b\} = 0$, yields the desired identity (\ref{e chernvector}). Now we establish Proposition \ref{p chernprop} piece by piece.  Note that \ref{p chernprop}.1 is immediate from the construction of $\P_w$.
\\

\noindent {\em Proof of \ref{p chernprop}.2.}  We have $\del \P_w = \gamma_w - \alpha_w - \sum_e c_e \beta_e$, for suitable $c_e$, so $\#(\del \P_w) = 2 + \sum_e |c_e|$.  The value $|c_e|$ is computed as the absolute value in the difference between the coefficients on the pair of regions atop $\Sigma$ that abut along $\beta_e$, which is $|a_e-b_e|$.
\\

\begin {figure}
\begin {center}
\includegraphics[width=4in]{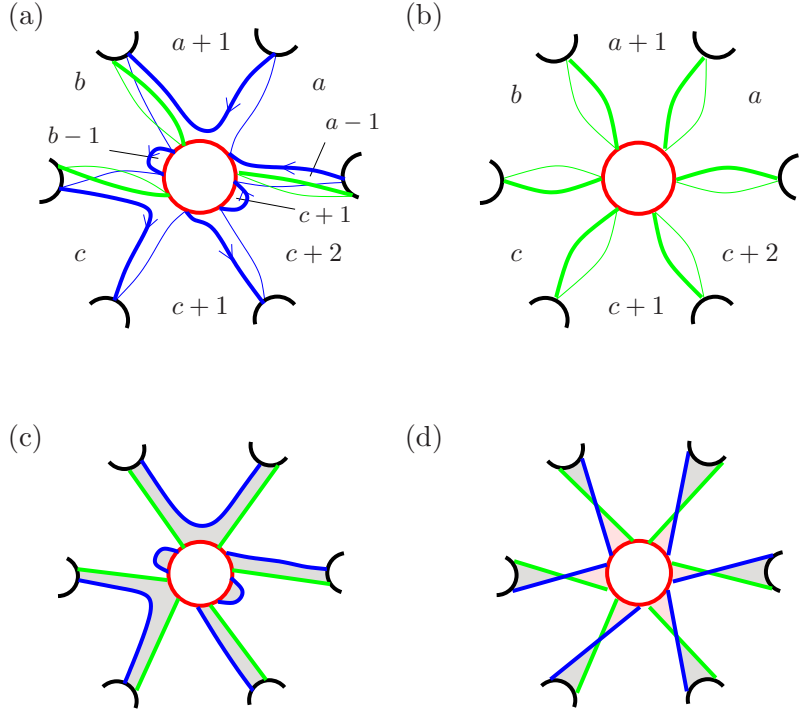}
\put(-300,260){(a)}
\put(-180,220){{\footnotesize $a-1$}}
\put(-185,235){{\small $a$}}
\put(-238,250){{\small $a+1$}}
\put(-285,215){{\footnotesize $b-1$}}
\put(-275,235){{\small $b$}}
\put(-275,170){{\small $c$}}
\put(-238,150){{\small $c+1$}}
\put(-190,185){{\footnotesize $c+1$}}
\put(-195,170){{\small $c+2$}}
\put(-150,260){(b)}
\put(-20,230){{\small $a$}}
\put(-73,250){{\small $a+1$}}
\put(-110,230){{\small $b$}}
\put(-110,170){{\small $c$}}
\put(-73,150){{\small $c+1$}}
\put(-30,170){{\small $c+2$}}
\put(-300,100){(c)}
\put(-150,100){(d)}
\caption{(a) The $2$-chain $\P_w$ nearby $\alpha_w$.  The coefficients on the top of $\P$ have been labeled.  (b) The $2$-chain $\P$ nearby $\alpha_w$.  (c) The top of the difference $2$-chain $\P-\P_w$, nearby $\alpha_w$.  (d) The bottom of $\P - \P_w$, nearby $\alpha_w$, with the top half of the Heegaard surface sliced away.  The gray regions have coefficient $+1$, and the pink ones have coefficient $-1$.}\label{f periodicmodify}
\end {center}
\end {figure}

\noindent {\em Proof of \ref{p chernprop}.3.}  We first modify $\P_w$ nearby the curve $\alpha_w$ to get a $2$-chain $\P$ with the same Euler measure $\widehat{\chi}$.  This is depicted in Figure~\ref{f periodicmodify}(b).  The boundary $\del \P$ is supported on edge curves $\beta_e$ ($\S$\ref{ss h_1 2}), $e \in E(T) \cup \Gamma(w)$.  Enumerate the edges in $\Gamma(w)$ in counterclockwise order.  The difference between the two-chains $\P - \P_w$ atop $\Sigma$ is a collection of bigons, squares, and hexagons, each with coefficient $+1$ (Figure~\ref{f periodicmodify}(c)).  There is a square between each consecutive pair of edges $e_1,e_2 \in \Gamma(w)$ with $\mu(e_1) = \mu(e_2)$, a bigon between each pair with $\mu(e_1) = +1$, $\mu(e_2) = -1$, and a hexagon between each pair with $\mu(e_1) = -1$, $\mu(e_2) = +1$.  The Euler measures of a bigon, square, and hexagon are $-1/2,0$, and $1/2$, respectively, and the number of bigons equals the number of hexagons, so the difference $2$-chain $\P - \P_w$ restricted to the top of $\Sigma$ has Euler measure zero.  The difference $2$-chain restricted to the bottom of $\Sigma$ consists of a pair of triangles of cancelling Euler measure for each edge $e \in \Gamma(w)$ (Figure~\ref{f periodicmodify}(d)), so it has Euler measure zero as well.  In total, $\widehat{\chi}(\P_w) =\widehat{\chi}(P)$, as desired.

Next, $\P$ decomposes into a sum of $2$-chains, one for each $f \in \Gamma(w) \setminus E(T)$, which we now describe.  Given such an edge $f$, the planar graph $T \cup \{ f \}$ contains a unique cycle $C(T,f)$, which may have some vertices in its interior.  There are as many of these vertices as there are edges interior to $C(T,f)$, since $T$ is a tree.  The cycle $C(T,f)$ gives rise to a subsurface $S \subset \Sigma$ whose boundary consists of edge curves $\beta_e$, $e \in C(T,f)$, and its genus is the number of edges interior to $C(T,f)$.  (An instructive example is the case of the edge $f_1$ in Figure~\ref{f periodic}(d).) Therefore, $\widehat{\chi}(S) = 2-2\# \{f \in E(T) \text{ interior to } C(T,f) \}-\#E(C(T,f))$.  Lastly, we sign $S$ by $\sigma(f)$ and denote by $S_f$ the resulting $2$-chain.  Then the sum of $S_f$ over all $f \in \Gamma(w) \setminus E(T)$ equals $\P_w$ (compare the remark following (\ref{e dx})), and we can calculate $$\widehat{\chi}(\P) = \sum_{f \in \Gamma(w) \setminus E(T)} \widehat{\chi}(S_f).$$

In order to recover the expression in \ref{p chernprop}.3, we reconsider the calculation of $\widehat{\chi}(S_f)$.  Write the value $\sigma(f)$ to each side of a curve in $\beta(x)$ that appears in $S_f$.  So $\sigma(f)$ appears once next to each boundary component of $S_f$ and twice next to each $\beta_e$, $e$ interior to $C(T,f)$.  The sum of these labels is $-\widehat{\chi}(S_f) - \sigma(f)$, hence the sum over all $f$ is $-\widehat{\chi}(\P) - \sum_{f \in \Gamma(w) \setminus E(T)} \sigma(f)$.  On the other hand, summing up the labels $\sigma(f)$ that appear next to the curve $\beta_e$ gives $a_e + b_e$, for each $e \in E(T) \setminus \Gamma(w)$, and $a_e+b_e-1$, for each $e \in E(T) \cap \Gamma(w)$.  It follows that $$-\widehat{\chi}(\P) - \sum_{f \in \Gamma(w) \setminus E(T)} \sigma(f) = \sum_{e \in E(T)} (a_e+b_e) - |\Gamma(w) \cap E(T)|.$$ Substituting $\widehat{\chi}(\P) = \widehat{\chi}(\P_w)$ into this equation results in \ref{p chernprop}.3.
\\

\noindent {\em Proof of \ref{p chernprop}.4.}  For the definition and use of the {\em dual spider number} $\sigma(u,\P)$, see \cite{OS4mfld}, Section 6, especially the discussion between Lemma 6.2 and the proof of Lemma 6.1. In the case at hand, $\sigma(\psi_x,\P_w)$ decomposes as a sum of terms $\sigma(\Delta_e,\P_w)$.  For each, we choose an interior point $p$ and arcs $a$, $b$, $c$ to $\del \Delta_e$ as shown in Figure~\ref{f smalltriangle}.  We note that the mirror image of this picture will appear in the case of a crossing $c$ with $\mu(c)=+1$, but this will have no effect on what follows.

If $e \in E(T) \setminus \Gamma(w)$, then $\Delta_e$ is disjoint from $\del_\alpha \P_w \cup \del_\gamma \P_w$, and so $\#(a \cap \del'_\alpha \P_w) = \#(c \cap \del'_\gamma \P_w) = 0$.  The value of $n_p(\Delta_e,\P_w)$ is given by $a_e$ or $b_e$ depending on the labeling, so we may assume it is $a_e$. Examining the cases $a_e = b_e, a_e > b_e$, and $a_e < b_e$ in turn, we deduce that
$\#(b \cap \del'_\beta \P_w) = \min \{ 0, b_e - a_e \}$.  In total,
\begin{eqnarray*}
\sigma(\Delta_e,\P_w) &=& n_p(\Delta_e,\P_w) + \#(a \cap \del'_\alpha \P_w) + \#(b \cap \del'_\beta \P_w) + \#(c \cap \del'_\gamma \P_w) \\ &=& a_e + \min \{ 0, b_e - a_e \} \\ &=& \min \{ a_e , b_e \}.
\end{eqnarray*}  If $e \in \Gamma(w) \cap E(T)$, it is still true that $\Delta_e$ is disjoint from the shifted boundary components $\del'_\alpha \P_w \cup \del'_\gamma \P_w$, giving $\#(a \cap \del'_\alpha \P_w) = \#(c \cap \del'_\gamma \P_w) = 0$.  The value of $n_p(\Delta_e,\P_w)$ in this case is given by either $a_e-1$ or $b_e-1$ depending on the labeling, so we may assume it is $a_e-1$.  Exactly as before, $\#(b \cap \del'_\beta \P_w) = \min \{ 0, b_e - a_e \}$, and we obtain $$ \sigma(\Delta_e,\P_w) = \min \{ a_e , b_e \}-1.$$  Summing the value of $\sigma(\Delta_e,\P_w)$ over all $e \in E(T)$ results in \ref{p chernprop}.4.

\subsection{$Spin^c$ structures.}\label{ss spinc}

The complex $\widehat{CF}(\h_2)$ decomposes into subcomplexes $\widehat{CF}(\h_2,\t)$ indexed by $\t \in Spin^c(\Sigma(K))$, and each Kauffman generator $x$ is contained in $\widehat{CF}(\h_2,\t(x))$ for some $\t(x) = \t_z(x)$.  We describe how to compute $\t(x)$.

In general, if $X^4$ is the result of attaching $2$-handles to $D^4$, and $Y = \del X$ is a rational homology sphere, then the long exact sequence in cohomology of the pair $(X, Y)$ splits off a short exact sequence $$0 \to H^2(X, Y) \to H^2(X) \to H^2(Y) \to 0.$$  The first Chern class $c_1: Spin^c(X) \to H^2(X)$ sets up a 1-1 correspondence between $Spin^c$ structures on $X$ and characteristic covectors in $H^2(X)$ for the intersection pairing on $X$, and $Spin^c(Y)$ is identified as the set of $2 \cdot H^2(X,Y)$-orbits of characteristic covectors.

In the case at hand, we have $X = \Sigma(F)$, an identification $H^2(X) \cong \Z^m$, and the intersection form given by the matrix $G_W$.  A characteristic covector $v \in \Z^m$ is defined by the condition that $v_i \equiv (G_W)_{ii} \, (\mod 2)$ for all $i$, and $Spin^c(\Sigma(K))$ is the set of $2 \cdot \text{im}(G_W)$-orbits of such vectors in $\Z^m$.  The $Spin^c$ structure $\t(x)$ lifts to $\s_z(\psi_x) \in Spin^c(X)$, and Equation (\ref{e chernvector}) expresses $c_1(\s_z(\psi_x)) = v_x^W \in \Z^m$.  Thus, $\t(x)$ corresponds to the orbit of $v_x^W$.  When $\det(K)$ is odd, the story simplifies somewhat.  In this case, there is a canonical identification $$c_1 : Spin^c(\Sigma(K)) \to H^2(\Sigma(K)) \cong \textup{coker}(G_W),$$ and under this identification, we have

\begin{equation}\label{e spinceqn}
\t(x) = [v_x^W] \in \textup{coker}(G_W).
\end{equation}

\subsection{Correction terms and alternating links.}\label{ss alternating}

We specialize now to the case that $K$ is a non-split alternating link, and let $D$ denote an alternating diagram of $K$.  Corollary \ref{c Lspace} implies that each Kauffman generator $x$ of $\widehat{CF}(\h_2(D))$ is unique in the $Spin^c$ structure $\t$ it occupies; consequently, the correction term in a given $Spin^c$ structure $\t = \t(x)$ is calculated as the absolute grading $\text{gr}(x)$.  Thus Theorem \ref{t grading} can be used to calculate the correction terms $d(\Sigma(K),\t)$.  On the other hand, Theorem 3.4 of \cite{OSdoublecover} provides such a formula as well.  Let us compare the two.

We begin by recalling the result of \cite{OSdoublecover} (cf. \cite{OSplumbed}, Corollary 1.5).  Let $G = G_W$ denote the Goeritz matrix corresponding to the white graph of $D$, $q$ its associated quadratic form, and $m$ its rank.  As in $\S$\ref{ss spinc}, we interpret $\t \in Spin^c(\Sigma(K))$ as a $2\cdot \textup{im}(G)$-orbit of characteristic covectors
for $G$.  With this notation, Theorem 3.4 of \cite{OSdoublecover} reads as follows:

\begin{equation} \label{e d-1}
d(\Sigma(K),\t) = \max _{v \in \t}\frac{ q(v) + m}{4}.
\end{equation}  On the other hand, the matrix $G$ is negative definite and $\delta(x) = 0$ for all $x$, so the expression (\ref{e absgrading}) reduces to

\begin{equation}\label{e d-2}
d(\Sigma(K),\t(x)) = \frac{ q(v_x^W) + m}{4}.
\end{equation}  Comparing Equations (\ref{e d-1}) and (\ref{e d-2}) gives the following result.

\begin{prop}\label{p d-3}

Let $x$ denote a Kauffman state in a connected alternating diagram.  The vector $v_x^W$ attains the maximum value of $q$ in the $2 \cdot \textup{im}(G)$-orbit $\t(x)$.

\end{prop}

It is interesting to give a direct algebraic proof of Proposition \ref{p d-3}.  Thus, we must identify $q(v_x^W)$ as the maximum value of $q$ over the $2 \cdot \text{im}(G)$-orbit of $v_x^W$.  In other words, we seek the inequality $q(v_x^W) \geq q(v_x^W-2Gy)$ for all integer vectors $y$; or, what is the same, that

\begin{equation}\label{e ineq}
y^T G y \leq v_x^W \cdot y, \text{ for all integer vectors } y.
\end{equation}  In fact, we claim that Inequality (\ref{e ineq}) holds if we replace $v_x^W$ by any vector $v$ induced by an orientation $\overrightarrow{W}$ of $W$.  Recall that this means that the $w_i$-entry of $v$ is equal to the number of edges directed out of $w_i$, minus the number directed into it, in $\overrightarrow{W}$.  Thus, given an orientation $\overrightarrow{W}$, form a matrix $M$ whose rows are indexed by $E(W)$, whose columns are indexed by $V(\widetilde{W})$, and whose $(e,w)$ entry is $0, 1,$ or $-1$ according as $e$ is not incident $w$, is directed away from it, or is directed into it in $\overrightarrow{W}$.  One checks that $$M^T M = -G \quad \mbox{ and } \quad v = M^T \bf{1},$$ where $\bf{1}$ denotes the all $1$'s vector of length $m$.  Therefore, $$y^T G y = -(My) \cdot (My) = - (z_1^2+\cdots+z_m^2)$$ and $$v \cdot y = {\bf 1} \cdot (My)=z_1+\cdots+z_m,$$ where $My = (z_1,\dots,z_m)^T$.  As the $z_i$ are integers, the inequality $y^T G y \leq v \cdot y$ follows and completes the proof of Proposition \ref{p d-3}.


\vspace{.5in}

\section{The domain of a homotopy class.}\label{S domain}

\subsection{Statement of the result.}\label{ss domain}

We briefly review how to compute the domain of a homotopy class connecting two generators in $\widehat{CF}(\h)$ in the case that $\h = (\Sigma,\alpha,\beta,z)$ presents a rational homology sphere $Y$.  Orient all the $\alpha$ and $\beta$ curves and $\Sigma$, let $\{ \eta_i \}$ denote the collection of oriented segments of $\alpha$ and $\beta$ curves that run between two consecutive intersection points, and fix generators $x,y \in \widehat{CF}(\h)$.  We write the 1-chain $$\epsilon(x,y)  = \sum_i d_i \eta_i, \quad d_i \in \Z,$$ which in turn bounds a {\em rational} 2-chain $$\sum_j q_j R_j, \quad q_j \in \Q.$$  Fix a segment $\eta_i$ and let $R_j$ and $R_k$ denote the regions that abut along it to the left and right, respectively.  We obtain the relation $$q_j-q_k = d_i,$$ and the collection of these relations for all $i$ determines the 2-chain up to a rational multiple of $\sum_j R_j$.  Enforcing $q_j=0$ for the region $R_j$ containing the basepoint $z$ removes this indeterminacy, and we denote the resulting 2-chain by $\D$.

\begin{defin}

The 2-chain $\D$ constructed in this way is {\em the domain connecting $x$ to $y$}.

\end{defin}

\noindent Notice that we make no assumption on $[\epsilon(x,y)] \in H_1(Y)$ in making this definition.  However, in the case that $[\epsilon(x,y)]=0$, the domain $\D$ equals $\D(\phi)$ for the unique homotopy class $\phi \in \pi_2(x,y)$ satisfying $n_z(\phi)=0$.

In Proposition \ref{p h2regions} we identified the regions in the Heegaard diagram $\h_2$.  In effect, they come in two types: those atop the diagram, which are in 1-1 correspondence with arcs of $D$, and those on the bottom, which are in 1-1 correspondence with regions of $D$.  For an arc $\alpha$ of $D$, we let $R_\alpha$ denote the corresponding region of $\h_2$.  Similarly, for a region $v$ of $D$, we let $R_v$ denote the corresponding region of $\h_2$.  The choice of marked point $p$ on $D$ will identify some of these regions, so we understand that $R_\alpha$ and $R_v$ may coincide for some choices of $\alpha$ and $v$.

Under these identifications, we split the determination of $\D$ into three parts.  To begin with, the assumption that $n_z(\D)=0$ allows us to restrict attention to unmarked arcs and regions of the diagram $D$.  Enumerate these arcs $\alpha_1,\dots,\alpha_k$, black vertices $b_1,\dots,b_\ell$, and white vertices $w_1,\dots,w_m$.  We define $$ \D_\alpha = \sum_{i=1}^k n_{\alpha_i} R_{\alpha_i}, \quad \D_B = \sum_{i=1}^\ell n_{b_i} R_{b_i}, \quad \mbox{and} \quad \D_W = \sum_{i=1}^m n_{w_i} R_{w_i}.$$  Thus, may we regard $\D_\alpha$ as a rational vector $(n_{\alpha_i})$ indexed by unmarked arcs of $D$, and similar statements hold for $\D_B$ and $\D_W$.  Next, mark a crossing in $D$ and enumerate the others $c_1,\dots,c_k$.  Given a Kauffman state $x$
and an unmarked crossing $c$, define

$$ \xi(x,c) =
\begin{cases}
1 , \quad \mbox{if $x(c)$ lies in a black region,} \\
0 , \quad \mbox{if $x(c)$ lies in a white region,}
\end{cases}
$$  and collect these values into the vector $$\xi_x = (\xi(x,c)).$$  Let $A$ denote the corresponding coloring matrix ($\S$\ref{ss prelim_2}), $G_B$ and $G_W$ the Goeritz matrices ($\S$\ref{ss prelim_1}), and $v^B_x$ and $v^W_x$ the degrees vectors ($\S$\ref{ss absgradings1}).

\begin{thm}\label{t domain}

For a pair of Kauffman states $x$ and $y$ in a marked diagram $D$, the domain $\D$ connecting $x$ to $y$ in the Heegaard diagram $\h_2(D)$ is given by the formulas

$$ \D_\alpha = A^{-1} (\xi_x - \xi_y), \quad \D_B = -{1 \over 2} G_B^{-1} (v^B_x - v^B_y), \quad \mbox{and} \quad \D_W = {1 \over 2} G_W^{-1} (v^W_x - v^W_y).$$

\end{thm}

The reason for the difference in signs in the formulas for $\D_B$ and $\D_W$ is due in effect to the asymmetric dependence of the incidence number $\mu$ on the black-white coloration: switching the roles of black and white reverses the value of $\mu$ at a crossing.

We remark on one peculiar feature of Theorem \ref{t domain}.  Recall that Equation (\ref{e chernvector}) in $\S$\ref{S absgradings} identified the degrees vector $v^W_x$ with minus the first Chern class of a specific $Spin^c$ structure on $\Sigma(F)$.  Its reincarnation in this context seems more than accidental, but a good explanation is lacking.  Regarding $G_W$ as the intersection pairing on $\Sigma(F)$ and $v^W_x, v^W_y \in H^2(\Sigma(F);\Z)$, the expression for $\D_W$ enables us to identify it as a class in $H_2(\Sigma(F);\Q)$; when $[\epsilon(x,y)]=0$, this class is integral.  Still, the significance of this class remains mystery. Of course, these same remarks apply with the roles of the black and white graphs switched.

In $\S$\ref{ss Dalpha} we settle the formula for $\D_\alpha$.  Its proof is a straightforward application of the procedure described at the beginning of this section.  In this situtation, we need only consider the segments $\eta_i$ which occur as portions of the $\beta$ curves atop $\h_2$.  In $\S$\ref{ss DB&W} we settle the formulas for $\D_B$ and $\D_W$.  We apply the same method, using only segments of the $\beta$ curves on the bottom of $\h_2$.  However, in order to put the resulting answer into the form given in Theorem \ref{t domain}, more work is needed, and in particular we must revisit the passage from $\h_2$ to $\h_1$.

\subsection{Calculation of $\D_\alpha$.}\label{ss Dalpha}

\begin{figure}[htb!]
\centering
\includegraphics[height=2in]{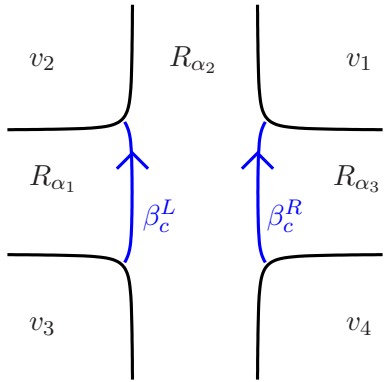}
\put(-82,120){$R_{\alpha_2}$}
\put(-135,75){$R_{\alpha_1}$}
\put(-20,75){$R_{\alpha_3}$}
\put(-92,60){\textcolor{blue}{$\beta^L_c$}}
\put(-45,60){\textcolor{blue}{$\beta^R_c$}}
\put(-15,120){$v_1$}
\put(-135,120){$v_2$}
\put(-135,20){$v_3$}
\put(-15,20){$v_4$}
\caption{The local picture atop $\h_2$ at a crossing $c$.}  \label{f Dtoplocal}
\end{figure}

Choose a crossing $c$ and let $\alpha_2$ denote the overstrand and $\alpha_1$, $\alpha_3$ the other two arcs that meet at $c$.  Orient the curve $\beta_c$ in the Heegaard diagram $\h_2$.  It has two segments $\beta^L_c$ and $\beta^R_c$ on top of $\h_2$ which appear oriented in the same direction, as in Figure \ref{f Dtoplocal}.  Notice that this is opposite the corresponding local picture in knot Floer homology \cite{OSaltknots}.  One of these segments will continue onto the bottom of $\Sigma$ in case the marked region is incident $c$, but this will not affect what follows.  Also shown in the picture are the regions $v_1,\dots,v_4$ that appear at $c$ ($\S$\ref{ss DB&W}).

Now fix a pair of Kauffman states $x$ and $y$.  Let $a_i$ denote the coefficient on $R_{\alpha_i}$ in the domain connecting $x$ to $y$, and $d_L,d_R$ the coefficients on segments $\beta^L_c,\beta^R_c$ in $\epsilon(x,y)$.  Thus $$a_1-a_2 = d_L \quad \mbox{and} \quad a_2-a_3 = d_R,$$ implying the single relation $$2a_2-a_1-a_3 = d_R-d_L.$$ Depending on the states $x(c)$ and $y(c)$, the difference $d_R - d_L$ takes on one of the values $0$ or $\pm 1$.  There are a few cases to consider.  In considering each, it is helpful to imagine a choice of spanning tree $T$ in constructing $\h_2$ and continuing the curve $\beta_c$ to the bottom of $\Sigma$ accordingly.  Ultimately, of course, this choice is immaterial.

\begin{itemize}

\item  {\em The states $x(c)$ and $y(c)$ agree or are diagonally opposite at $c$.}  In either case $d_L = d_R$, which implies the relation $2a_2 - a_1 - a_3 = 0$.

\item  {\em The states $x(c)$ and $y(c)$ lie to the same side of the overstrand, but differ.}  Assume that $x(c)$ lies in the bottom-right quadrant of Figure \ref{f Dtoplocal} and $y(c)$ lies in the upper-right.  Then $d_R - d_L= 1$ and we obtain the relation $2a_2-a_1-a_3 = 1$.  A different placement of $x(c)$ and $y(c)$ can be transformed into this one by (i) rotating the local picture, (ii) swapping the roles of $x$ and $y$, or (iii) both.  The effect of (i) is to reverse the roles of $d_R$ and $d_L$ as well as $a_1$ and $a_3$, giving the relation $2a_2-a_1-a_3=-1$ instead.  The effect of (ii) is to negate $a_1,a_2,$ and $a_3$, which also results in the relation $2a_2-a_1-a_3=-1$.  These effects are canceled in (iii) and result in the original relation.

\item  {\em The states $x(c)$ and $y(c)$ lie to the same side of the understrand, but differ.}  Assume that $x(c)$ lies in the bottom-right quadrant of Figure \ref{f Dtoplocal} and $y(c)$ lies in the bottom-left.  We obtain $d_R - d_L = 1$ and hence the relation $2a_2-a_1-a_3=1$.  The three other possible placements of $x(c)$ and $y(c)$ can be treated now exactly as in the previous case.

\end{itemize}

Now consider the black-white coloring of the regions of $D$.  It is easy to check that the expressions for $2a_2-a_1-a_3$ in each of the preceding cases condense into the following unified form: $$-\mu(c) \cdot (2a_2-a_1-a_3) = \xi(x,c)-\xi(y,c).$$  The collection of these relations, one for each crossing $c$, is expressed in the matrix equation $$A \cdot \D_\alpha = \xi_x -\xi_y,$$ where $A$ denotes the coloring matrix.  This gives the formula for $\D_\alpha$ expressed in Theorem \ref{t domain}.

\subsection{Calculation of $\D_B$ and $\D_W$.}\label{ss DB&W}

We proceed as in the previous subsection, obtaining a relation on the coefficients of $\D$ on the bottom of $\h_2$, one for each crossing $c$ of the diagram.  Let $\beta^{13}_c$ and $\beta^{24}_c$ denote the two segments of $\beta_c$ appearing on the bottom of $\h_2$, where $\beta^{13}_c$ runs between regions $v_1$ and $v_3$ at $c$ and $\beta^{24}_c$ runs between $v_2$ and $v_4$.  Let $d_{13}$ and $d_{24}$ denote the coefficients on these segments in $\epsilon(x,y)$ and $a_i$ the coefficient on region $R_{v_i}$, $i = 1,\dots,4$.  We obtain the relations $$a_1-a_3 = d_{13} \quad \mbox{and} \quad a_2-a_4 = d_{24},$$ whence the single relation $$-a_1+a_2+a_3-a_4 = d_{24}-d_{13}.$$  Again, this is most easily seen by making a choice of spanning tree in constructing $\h_2$ and extending the curve $\beta_c$ onto the bottom of $\h_2$.  As in the case of $\D_\alpha$, depending on the states $x(c)$ and $y(c)$, the difference $d_{24} - d_{13}$ takes on one of the values $0$ or $\pm 1$.  Treating each case in turn, we obtain the following concise result.  Define $$\psi(x,c) = \pm 1, \mbox{ according as $x(c)$ lies to the right / left of the oriented overstrand at $c$,}$$ and collect these values into a vector

\begin{equation}\label{e psivector}
\psi_x = (\psi(x,c)).
\end{equation}  Then

\begin{equation}\label{e localrelation}
-a_1+a_2+a_3-a_4 = {1 \over 2}(\psi(x,c)-\psi(y,c)).
\end{equation}  Verification is routine, and the collection of these relations, one for each crossing $c$, is expressed in a matrix equation

\begin{equation}\label{e DB&W1}
M_0 \cdot (\D_B \oplus \D_W) = {1 \over 2}(\psi_x-\psi_y).
\end{equation}

Equation (\ref{e DB&W1}) is the most direct way of expressing the defining relation for $\D_B$ and $\D_W$.  However, in order to convert (\ref{e DB&W1}) into the form given in Theorem \ref{t domain}, we first modify it into the intermediate form (\ref{e DB&W2}).  First, make a choice of spanning tree $T \subset B$.  This induces orientations of $B$ and $W$ and gives rise to a Kauffman state $x = x_T$ ($\S$\ref{ss absgradings1}).  Every edge $e \in E(T) \cup E(T^*)$ gives rise to a segment of a $\beta$ curve on the bottom of $\h_2$ running between the same regions as $e$, and we orient this segment according to the direction on $e$.  Notice that at a given crossing $c$, the orientations on the two segments of $\beta_c$ on the bottom of $\h_2$ due to $e$ and $e^*$ are compatible in that they induce the same orientation of $\beta_c$.  So the choice of $T$ (or equivalently, the Kauffman state $x$) defines an orientation of the curves $\beta_c$.  Now define $$d(x,c) = \pm 1, \mbox{ according as the segment of $\beta_c$ is directed out of / into the region containing $x(c)$.}$$  In the same way we define $d(v,c)$ for a region $v$ incident $c$, and set $d(v,c) = 0$ in case $v$ is not incident $c$ at all or it meets it twice.  This assignment is depicted in Figure \ref{f d(v,c)}.

\begin{figure}[htb!]
\centering
\includegraphics[height=1.5in]{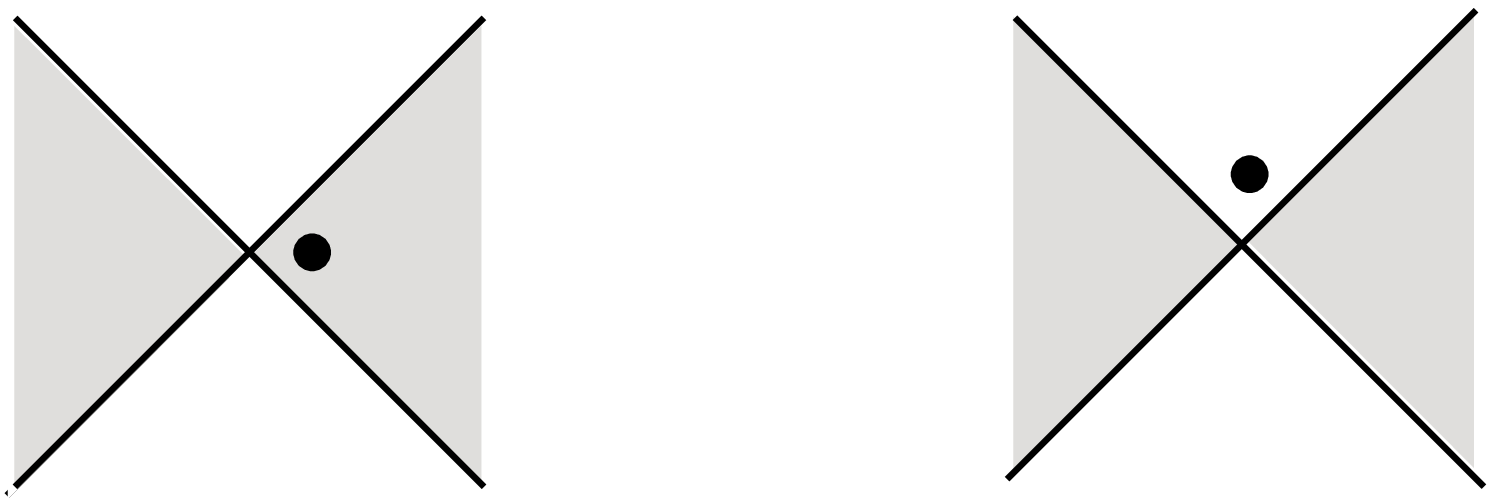}
\put(-235,50){$-1$}
\put(-315,50){$1$}
\put(-275,90){$\mu(c)$}
\put(-285,10){$-\mu(c)$}
\put(-65,90){$-1$}
\put(-55,10){$1$}
\put(-25,50){$\mu(c)$}
\put(-102,50){$-\mu(c)$}
\caption{The incidence number $d(v,c)$.  Pictured is the crossing $c$, where the black dot represents the position of $x_T(c)$.  Assuming the region $v$ meets the crossing $c$ once, the value of $d(v,c)$ is given by the number written in it; otherwise it takes the value $0$.}  \label{f d(v,c)}
\end{figure}

\noindent Also, set $$\kappa(x,c) = \pm 1, \mbox{ according as $x(c)$ is in a white / black region,}$$ and $$\alpha(x,c) = \mu(c) \cdot d(x,c) \cdot \kappa(x,c).$$  With these definitions in place, a simple check shows that (\ref{e localrelation}) can be put into the following form: $$\sum_{j=1}^4 d(v_j,c) a_j = {1 \over 2}(\alpha(x,c)  - \alpha(y,c));$$  and the collection of these relations, one for each crossing $c$, condense into a matrix equation

\begin{equation}\label{e DB&W2}
M \cdot (\D_B \oplus \D_W) = {1 \over 2}(\alpha_x-\alpha_y).
\end{equation}  Here $\alpha_x$ denotes the vector of values $\alpha(x,c)$, while $M$ is the matrix of values $(d(v,c))$, its rows indexed by all crossings $c$ and its columns by unmarked regions $v$.

In spite of the additional data required to express (\ref{e DB&W2}) as compared with the equivalent (\ref{e DB&W1}), it is easier to convert it into the desired form.  We claim that

\begin{equation}\label{e alphad}
M^{-1} \cdot \alpha_x = (- G_B^{-1} d^B_x) \oplus (G_W^{-1} d^W_x).
\end{equation}  Together, (\ref{e DB&W2}) and (\ref{e alphad}) imply the formulas for $\D_B$ and $\D_W$ in Theorem \ref{t domain}.  Observe that while both $M$ and $\alpha$ depend on the choice of spanning tree $T \subset B$, the left-hand side of (\ref{e alphad}) is independent of this choice: indeed, a different choice of $T$ results in multiplying some columns of $M$ and the corresponding entries of $\alpha_x$ by $-1$, which get canceled in the product $M^{-1} \cdot \alpha_x$.  Consequently, it suffices to prove (\ref{e alphad}) in the special case that $x = x_T$.  Notice that $d(x_T,c) = -1$ for all crossings $c$, which simplifies the expression for $\alpha(x_T,c)$, a convenient fact we use below.

We realize (\ref{e alphad}) by examining the passage from $\h_2$ to $\h_1$ ($\S$\ref{ss heegaardmoves}) at the homological level.  Let $\h_1'$ denote the intermediate Heegaard diagram obtained prior to the destabilizations.  The matrix $M^T$ is the intersection matrix $I(\h_2)$ (Definition \ref{d intersectionmatrix}), having oriented the $\alpha$ curves of $\h_2$ to run counterclockwise and the $\beta$ curves as above.  The matrix $(G_B | 0)^T$ is the submatrix of $I(\h_1')$ induced on the columns $\beta_i$ corresponding to edges in the tree $T \subset B$.  Now, if two Heegaard diagram $H_1$ and $H_2$ are related by a handleslide of a $\beta$ curve, then $I(H_2)$ is related to $I(H_1)$ by an elementary column operation.  If they are related instead by an isotopy, then $I(H_1) = I(H_2)$. Thus, the passage from $\h_2$ to $\h_1'$ entails a sequence of {\em row} operations which transform $M$ into $I(\h_1')^T$. In short, there exists a matrix $E_B$ with the property that

\begin{equation}\label{e E}
E_B \cdot M = (G_B \; | \; 0).
\end{equation}

\noindent Similarly there is a matrix $E_W$ so that $E_W \cdot M = (0 \; | \; G_W)$.  It corresponds to the passage from $\h_2$ to the Heegaard diagram $\h_1(D')$, where $D'$ agrees with $D$ but with the opposite coloration of its regions.  In what follows we concentrate on the case of $E_B$, as the case of $E_W$ follows analogously.

\begin{figure}[htb!]
\centering
\includegraphics[height=1.5in]{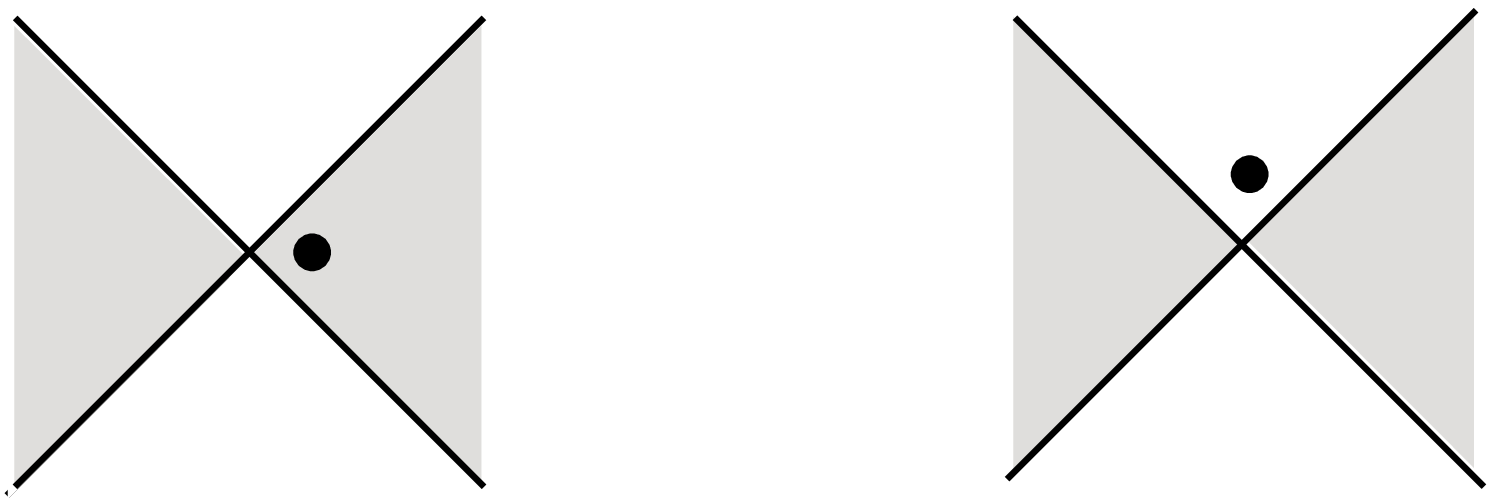}
\put(-247,50){$-\mu(c)$}
\put(-315,50){$\mu(c)$}
\put(-20,50){$+1$}
\put(-100,50){$-1$}
\caption{The incidence number $\sigma(v,c)$.  Pictured is the crossing $c$, where the black dot represents the position of $x_T(c)$.  Assuming the region $v$ meets the crossing $c$ once, the value of $\sigma(v,c)$ is given by the number written in it; otherwise it takes the value $0$.}  \label{f sigma}
\end{figure}

It stands to actually write down an expression for $E_B$ obeying (\ref{e E}) and also establish

\begin{equation}\label{e E0}
E_B \cdot \alpha_x = - d^B_x;
\end{equation}

\noindent together, (\ref{e E}) and (\ref{e E0}) prove the part of (\ref{e alphad}) concerning $B$, and the part concerning $W$ follows in exactly the same way.  In order to write down $E_B$, we define one final type of incidence number $\sigma(v,c)$ between an unmarked black vertex $v$ and crossing $c$, which depends on the Kauffman state $x_T$.  Its value is given as in Figure \ref{f sigma}.  We set $E_B$ equal to the matrix of values $\sigma(v,c)$, indexed by unmarked black vertices $v$ and all crossings $c$.  Now we proceed to establish (\ref{e E}).  Fix $v \in V(B)$ and $u \in V(B) \cup V(W)$.  The $(v,u)$-entry of the product $E_B \cdot M$ is $$(E_B \cdot M)_{vu} = \sum_c \sigma(v,c) d(u,c).$$  The general term $\sigma(v,c) d(u,c)$ in the summation is non-zero unless the crossing $c$ touches both regions $v$ and $u$.  Assuming this occurs, there are a few cases to consider.
\\

\noindent 1. $u \in V(B), u \ne v$.  If $x(c)$ lies in $v$, then $\sigma(v,c) = \mu(c)$ and $d(u,c) = 1$, so their product is $\mu(c)$.  If $x(c)$ instead lies in $u$, then the signs on $\sigma(v,c)$ and $d(u,c)$ negate and their product remains the same: $\mu(c)$.  If $x(c)$ instead lies in one of the white regions, there are four possibilities to consider, depending on which white region and the type of crossing at $c$.  In any event, we always obtain $$\sigma(v,c)d(u,c) = \mu(c).$$  It follows then that  $$(E_B \cdot M)_{vu} = \sum_{c \text{ incident } u \text{ and } v} \mu(c).$$
\\

\noindent 2. $u = v$.  Note that $d(c,v) = -d(c,u')$ for $u'$ the black region which touches $c$ opposite to the corner of $v$ under consideration.  In particular, it is possible that $u' = v$, in which case both sides of this equation are $0$.  The crossings $c$ incident $v$ are partitioned according to the other black region to which they are incident.  Using the formula worked out in the preceding case, we conclude that $$(E_B \cdot M)_{vv} = - \sum_{c \text{ incident } v} \mu(c).$$
\\

\noindent 3. $u \in V(W)$.  First suppose that no crossing $c$ meets the same region $u$ twice.  The value of the product $\sigma(v,c) d(u,c)$ can be gotten from Figure \ref{f d(v,c)sigma}.  The regions $u$ and $v$ abut along disjoint edges in the knot diagram; here we view the knot diagram as a planar graph whose vertices are the crossings.  Each edge ends at a pair of crossings $c_1$ and $c_2$.  A simple check shows that the values $\sigma(v,c_1) d(u,c_1)$ and $\sigma(v,c_2) d(u,c_2)$ cancel one another regardless of the placement of $x(c_1)$ and $x(c_2)$, using only the property that $x(c_1)$ and $x(c_2)$ lie in different regions.  It follows that the $(v,u)$-entry of the product $E_B \cdot M$ is equal to $0$ in this case.  In the event that some crossing meets the same region $u$ twice, essentially the same argument carries through.  In this case the regions $u$ and $v$ abut along disjoint paths in the knot diagram, where vertices internal to a path correspond to crossings which the region $u$ meets twice.  For each internal crossing we have $d(u,c) = 0$, and for the endpoints of the paths we get canceling contributions just as in the previous case.  Consequently, $$(E_B \cdot M)_{vu} = 0$$ in any event.
\\

\begin{figure}[htb!]
\centering
\includegraphics[height=1.5in]{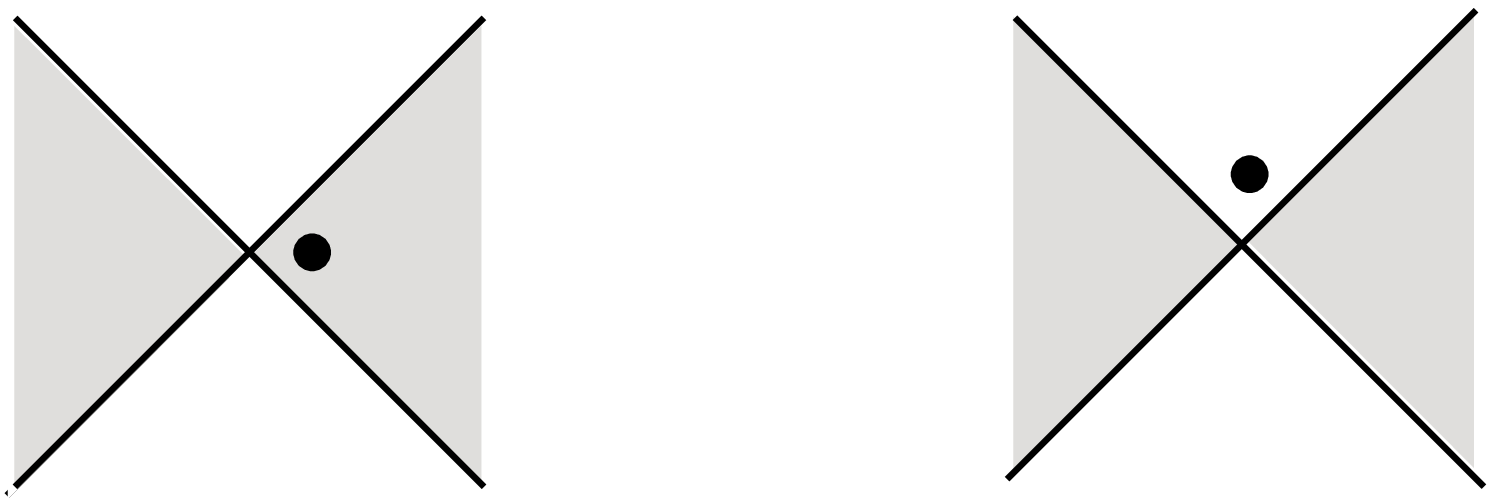}
\put(-235,50){$-1$}
\put(-315,50){$+1$}
\put(-275,90){$+1$}
\put(-275,10){$-1$}
\put(-60,90){$-1$}
\put(-60,10){$+1$}
\put(-20,50){$+1$}
\put(-100,50){$-1$}
\caption{The incidence number $\sigma(v,c)d(u,c)$.  Pictured is the crossing $c$, where the black dot represents the position of $x(c)$.  Assuming the region $u$ meets $c$ once, the value of $\sigma(v,c)d(u,c)$ is gotten by multiplying together the values in the regions $u$ and $v$.}  \label{f d(v,c)sigma}
\end{figure}

Cases 1-3 together establish (\ref{e E}).  Now we settle (\ref{e E0}).  We seek the identity
\begin{equation}\label{e E2}
\sum_c \sigma(v,c) \alpha(x,c) = - d_x(v) ,\text{ for all }v \in V(B).
\end{equation}

\noindent  Recall that for $x := x_T$, $d(x,c) = -1$ for all crossings $c$, so the expression for $\alpha(x,c)$ reduces to $-\mu(c) \kappa(x,c)$.  Suppose that $x(c)$ is black.  Then $\sigma(v,c) = \pm \mu(c)$ according as $x(c)$ is not or is contained in the region $v$, while $\alpha(x,c) = \mu(c)$ regardless.  Hence the product $\sigma(v,c) \alpha(x,c) = \pm 1$ according as $x(c)$ is not or is contained in the region $v$.  Recalling the orientation on $B$ induced by $T$ ($\S$\ref{ss absgradings1}), this is to say that the edge $e_c$ is oriented into or out of the vertex $v$.  Next suppose that $x(c)$ is white; thus $\alpha(x,c) = - \mu(c)$.  Examining the few possibilities, we see that the value of the product $\sigma(v,c) \alpha(x,c)$ is also $\pm 1$ according as the edge $e_c \in E(B)$ is directed into or out of the vertex $v$.  Adding up all these contributions results in the identity (\ref{e E2}) and hence (\ref{e E0}).  Tracing backwards, we have at last settled (\ref{e alphad}), which completes the proof of Theorem~\ref{t domain}.


\vspace{.5in}

\section{A spectral sequence.}\label{S specseq}

\subsection{Construction of the spectral sequence.}\label{ss specseq}

As pointed out in $\S$\ref{ss h_2}, there is a strong similarity between the standard doubly-pointed Heegaard diagram $\h(D) = (\Sigma, \alpha, \beta, z, z')$ presenting the knot $K \subset S^3$ and the Heegaard diagram $\h_2(D) = (\Sigma, \alpha, \beta', z)$ presenting $\Sigma(K)$, starting with a marked diagram $D$ of $K$.  In this section, we exploit this similarity to construct a pair of spectral sequences $\{(E_r,d_r)\}$ and $\{(E'_r,d'_r)\}$ which converge to $E_\infty \cong \widehat{HFK}(K)$ and $E'_\infty \cong \widehat{HF}(\Sigma(K))$, respectively.  Both terms $E_0$ and $E'_0$ are free abelian, generated by Kauffman states of the diagram $D$, and moreover the complexes $(E_0,d_0)$ and $(E'_0,d'_0)$ are isomorphic.  It follows that $E_1 \cong E'_1$.  Although the term $E_1$ is not a knot invariant, it admits a simple combinatorial description (Theorem \ref{t E1}), it is independent of the choice of complex structure on $\Sigma$, and it is useful in making some of the calculations in $\S$\ref{S calculations}.  The two spectral sequences arise out of a standard construction in Heegaard Floer homology that we describe in the next paragraph.  After setting this up, the remainder of the section is directed towards establishing Theorem \ref{t E1}.

To begin with, arrange so that the restriction of each pair $\beta_i$, $\beta'_i$ to the top of $\Sigma$ is the same; then the tops of both $\h$ and $\h_2$ are the same, and the intersection points in $\T_\alpha \cap \T_\beta$ are in 1-1 correspondence with Kauffman states of $D$ in both cases.  Choose a complex structure $J$ on $\Sigma$ and a generic complex structure on $\text{Sym}^g(\Sigma)$ close to the one induced by $J$.  The complex $\widehat{CF}(\h)$ decomposes into a direct sum of subcomplexes, where a pair of generators $x,y \in \widehat{CF}(\h)$ belong to the same subcomplex iff there is a domain $\D$ connecting them which satisfies $n_z(\D) = n_{z'}(\D) = 0$.  Place a basepoint in each region that meets the bottom of $\h$.  We define a relative filtration on each subcomplex by setting $$\F(x,y) = \sum_w n_w(\D),$$ the sum over all basepoints $w$ distinct from $z$ and $z'$.  The domain $\D$ is uniquely determined by the condition $n_z(\D) = 0$, so $\F$ is well-defined. In the same way, we obtain a relative filtration $\F_2$ on each subcomplex $\widehat{CF}(\h_2,\s)$, $\s \in Spin^c(\Sigma(K))$.  The relatively filtered complexes $(\widehat{CF}(\h),\f)$ and $(\widehat{CF}(\h_2),\f_2)$ give rise to a pair of affiliated spectral sequences $\{(E_r,d_r)\}$ and $\{(E'_r,d'_r)\}$, respectively.  Both terms $E_0$ and $E'_0$ are free abelian, generated by Kauffman states of $D$, and the sequences converge to $E_\infty \cong \widehat{HFK}(K)$ and $E'_\infty \cong \widehat{HF}(\Sigma(K))$, respectively.  Both differentials $d_0$ and $d'_0$ count pseudo-holomorphic disks $\phi$ with $\sum_w n_w(\D(\phi))=0$; as $n_w(\D(\phi)) \geq 0$ for any pseudo-holomorphic disk $\phi$ connecting $x$ to $y$, it follows that these differentials count disks $\phi$ for which the support of $\D(\phi)$ is contained in the top of its respective Heegaard diagram.  Moreover, the closure of the union of regions not containing a basepoint $w$, $z$, or $z'$ is the same for both $\h$ and $\h_2$.  It follows that $(E_0,d_0) \cong (E'_0,d'_0)$, and consequently $E_1 \cong E'_1$.

\subsection{Subcomplexes of $(E_0,d_0)$.}\label{ss subcomplexes}

In order to describe $E_1$, we first examine this question: given a pair of generators $x$ and $y$, when does there exist a domain $\D$ connecting $x$ to $y$ which is supported on top of $\Sigma$?  We can ask this question with reference to either $\h$ or $\h_2$; of course, the answer is the same in either case.  Recall that the Kauffman state $x$ gives rise to a spanning tree $T_W(x) \subset W$, which in turn induces an orientation of $W$ ($\S$\ref{ss absgradings1}); we refer to this as the orientation on $W$ {\em induced} by $x$.

\begin{lem}\label{l support}
The pair of generators $x$ and $y$ are connected by a domain $\D$ supported on top of $\Sigma$ iff $x(c)$ and $y(c)$ lie to the same side of the overstrand at $c$ for every crossing $c$, which happens iff $x$ and $y$ induce the same orientation on $W$ (equivalently, on $B$).
\end{lem}

\begin{proof}

It is easy to check that the generators $x$ and $y$ induce the same orientation on $W$ (equivalently, on $B$) iff $x(c)$ and $y(c)$ lie to the same side of the overstrand at $c$, for every crossing $c$.  This is to say that $\psi_x = \psi_y$ ($\S$\ref{ss DB&W}, Equation (\ref{e psivector})). By Equation (\ref{e DB&W1}), this is the case iff $\D$ is supported on top of $\Sigma$, taking note that the intersection matrix $M_0$ for $\Sigma(K)$ is invertible ($\det(K) \ne 0$).

\end{proof}

Following Lemma \ref{l support}, the complex $(E_0,d_0)$ decomposes into a direct sum of subcomplexes $\oplus_\psi \; \c_\psi$, where $\c_\psi$ is generated by Kauffman states $x$ satisfying $\psi_x = \psi$.  We now describe how the generators of a given subcomplex $\c_\psi$ are related.  Thus, fix a pair of Kauffman generators $x,y$ in $\c_\psi$.  Let $H$ denote the bipartite graph of incidences between regions and crossings in the diagram $D$.  The Kauffman state $x$ has a natural interpretation as a matching $m_x$ in $H$ between the {\em unmarked} regions and crossings, and conversely every such matching corresponds to a Kauffman state.  We construct another maximal matching $m'_x$ in $H$ by pairing each crossing $c$ with the other region that meets $c$ to the same side of the overstrand as $x(c)$.  A Kauffman state $y$ satisfies $\psi_y = \psi_x$ exactly when $m_y$ is a subset of the edges in $m_x \cup m'_x$. Notice that the subgraph on the edges in $m_x \cup m'_x$ consists of some number of (even) cycles $C_1,\dots,C_n$, as well as two star-shaped components with star vertex at either marked region.  Therefore, there are $2^n$ such Kauffman states $y$.

\begin{defin}\label{d solitary}
The Kauffman state $x$ is {\em solitary} if there are no cycles in $m_x \cup m'_x$.
\end{defin}

\noindent Thus a solitary Kauffman state $x$ is uniquely determined by the vector $\psi = \psi_x$, and it alone generates $\c_\psi$.  Examining $\c_\psi$ more closely in the general case, two-color the edges of $C_1,\dots,C_n$ red and blue. The edges of $H$ can be identified with the intersection points between the $\alpha$ and $\beta$ curves in $\h$ (aside from the single intersection point on $\beta_p$).  Let $r_i$ denote the tuple of intersection points corresponding to the red edges and $s_i$ the tuple corresponding to the blue edges of $C_i$ for every $i$, and let $t$ denote the tuple of intersection points corresponding to the edges of $m_x$ not in any $C_i$.  Then the Kauffman generators of $\c_\psi$ are the elements of the set $\G_\psi = \{r_1,s_1\} \times \cdots \times \{r_n,s_n\} \times \{ t \}$.

\begin{figure}
\begin{center}
\includegraphics[width=5in]{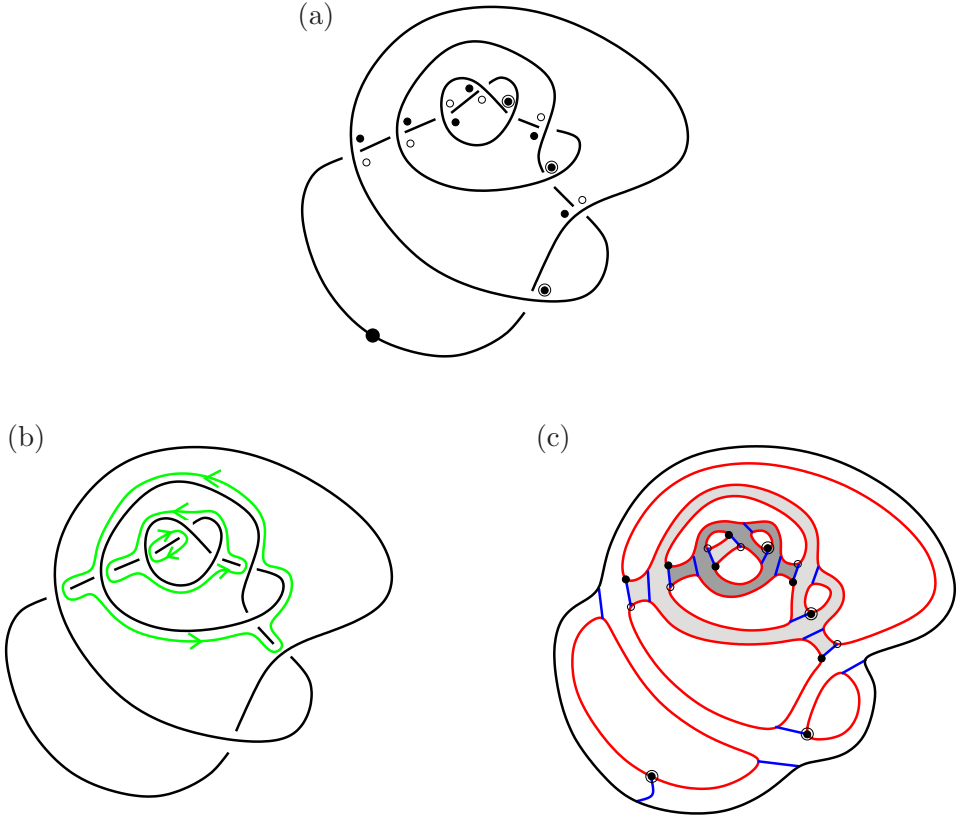}
\put(-250,300){(a)}
\put(-360,140){(b)}
\put(-160,140){(c)}
\caption{(a) A marked knot diagram with a pair of Kauffman states $x$ (corresponding to the open circles) and $y$ (corresponding to the filled-in circles).  (b) The three oriented curves $\gamma_i$ connecting $x$ to $y$, pictured in green.  (c) The domain $\D$ connecting $x$ to $y$, interpreted as Kauffman generators.  The lightly shaded regions have coefficient $+1$ and the darkly shaded region has coefficient $+2$.}
\label{f E1domain} 
\end{center}
\end{figure}

Next, we describe the domain connecting a pair of generators $x,y \in \G_\psi$ (cf. \cite{OSmutation}, Section 2.2 for an alternative approach). We begin by constructing a $1$-cycle $\gamma$ in $\Sigma$ connecting the components of $x$ and $y$, as in \cite{OS3mfld1}, Section 2.4.  For every bounded region $R$ with $x(R) \ne y(R)$, traverse the portion of the curve $\alpha_R$ that runs clockwise from $x(R)$ to $y(R)$; if $R$ is unbounded, traverse it in the counterclockwise direction.  For every crossing $c$ for which $x(c) \ne y(c)$, traverse the portion of the curve $\beta_c$ that runs between $y(c)$ and $x(c)$ atop $\Sigma$.  The union of these oriented portions of curves is the $1$-cycle $\gamma$.  This $1$-cycle decomposes into a union of disjoint, oriented curves $\gamma_i$ which connect the points of $r_i$ and $s_i$, one for each $i$ for which $C_i \subset m_x \cup m_y$.  We regard the $\gamma_i$ as a collection of disjoint, oriented curves in the plane of the knot diagram $D$.  Every arc $a$ in the diagram $D$ is separated from the unbounded region by some subset of the $\gamma_i$, and we set $c(a)$ equal to the signed count of the curves in this subset, where $\gamma_i$ gets weight $\pm 1$ according to whether it is oriented counterclockwise or clockwise.  Identify an arc $a$ with the corresponding region in the Heegaard diagram $\h$ (equivalently, $\h_2$).  See Figure~\ref{f E1domain}.

\begin{lem}\label{l E1domain}
The domain connecting $x$ to $y$ is $$\D = \sum_a c(a) \cdot a.$$
\end{lem}

\noindent We remark that $\del \D$ does not equal $\gamma$ on the nose in general, but the following proof implies that the two always differ by some number of copies of full $\alpha$ curves.

\begin{proof}

Let $\D' = \sum_a c(a) \cdot a$ and let $\D$ denote the domain connecting $x$ to $y$.  The domain $\D$ is supported on top of $\Sigma$, hence so is $\del \D$.  It follows that $\del_{\beta_c} \D$ consists of the oriented arc of $\beta_c$ from $y(c)$ to $x(c)$ on top of $\Sigma$ for every crossing $c$ (which is void in case $x(c) = y(c)$).  On the other hand, this is precisely $\del_{\beta_c} \D'$ restricted to the top of $\Sigma$.  According to $\S$\ref{ss Dalpha}, this is sufficient to conclude that the restrictions of $\D$ and $\D'$ to the top of $\Sigma$ agree.  Since $\D$ is supported on top of $\Sigma$, every region which meets both the top and bottom of $\Sigma$ appears with $0$ coefficient in $\D$ and hence in $\D'$.  Lastly, every region supported on the bottom of $\Sigma$ appears with coefficient $0$ in $\D'$.  It follows that $\D = \D'$, as desired.

\end{proof}

\subsection{Punctured polygons.}\label{ss polygon}

In the case that $m_x \cup m_y$ contains a single cycle $C_i$, the $1$-cycle $\gamma$ connecting $x$ to $y$ is the single curve $\gamma_i$, which we may take to be oriented counter-clockwise by switching the roles of $x$ and $y$ if necessary.  In this case, the domain $\D$ in this case is a {\em punctured polygon} according to the following definition.  Recall that a domain has an {\em interior corner} at an intersection point $q \in \alpha_i \cap \beta_j$ if it has coefficient $1$ on one of the four corners of regions incident $q$ and $0$ on the three others.

\begin{defin}\label{d polygon}

A domain $\D$ from $x$ to $y$ is a {\em punctured polygon} if (a) it is a planar region with every coefficient $0$ or $1$, (b) there are no components of $x$ and $y$ interior to $\D$, (c) all but one component of $\del \D$ is a full $\alpha$ curve, and (d) all corners of $\D$ are interior.

\end{defin}

Given a punctured polygon $\D$ from $x$ to $y$, fix a component $y_i \in \alpha_i \cap \beta_i$ of $y$.  If $y_i$ lies on $\D$, let $\widehat{\beta}_i$ denote the maximal arc of $\beta_i$ which lies in $\D$ and has both its endpoints on $\del \D$; by Definition \ref{d polygon} (b), $y_i$ is one of these endpoints. If $y_i$ does not lie on $\D$, let $\widehat{\beta}_i = y_i$.  Now define a directed graph $\Gamma(\D)$ on the $\alpha$ curves by declaring $(\alpha_i,\alpha_j)$ to be an edge if $\widehat{\beta}_i$ runs from $y_i$ to a point on $\alpha_j$.  The graph $\Gamma(\D)$ contains a directed cycle corresponding to its unique polygonal boundary component, as well as a directed loop at each $\alpha_i$ for which $\widehat{\beta}_i = y_i$.

\begin{defin}\label{d arborescent}

The punctured polygon $\D$ is {\em arborescent} if there are no other cycles in the graph $\Gamma(\D)$.

\end{defin}

\noindent It is easy to see that the property of arborescence is equivalent to the property that $\D - (\widehat{\beta}_1 \cup \cdots \cup \widehat{\beta}_n)$ is connected.  The importance of Definition \ref{d arborescent} is the following fact, which is in essence Lemma 3.11 of \cite{OScube}.

\begin{lem}\label{l polygonlemma}

If the domain of a homotopy class $\phi \in \pi_2(x,y)$ is a punctured polygon, then $\mu(\phi) = 1$; if in addition it is arborescent, then $\#\widehat{\M}(\phi) = \pm 1$, for any choice of complex structure on $\text{Sym}^g(\Sigma)$ sufficiently close to one induced by a complex structure on $\Sigma$.

\end{lem}

\begin{figure}
\begin{center}
\includegraphics[width=2in]{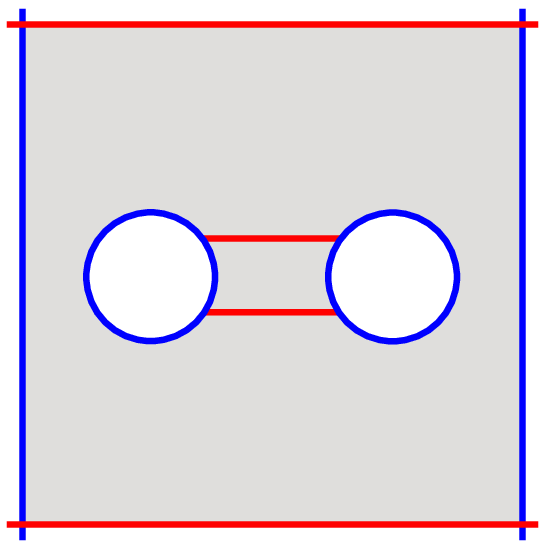} 
\put(-77,100){$\D_1$}
\put(-77,70){$\D_2$}
\put(-150,0){$a$}
\put(0,140){$a$}
\put(0,0){$b$}
\put(-150,140){$b$}
\put(-100,62){$x$}
\put(-50,77){$x$}
\put(-100,77){$y$}
\put(-50,62){$y$}
\caption{The domain $\D = \D_1 + \D_2$ of a map $\phi$ for which $\#\widehat{\M}(\phi)$ will depend on the choice of complex structure.}\label{f poly}
\end{center}
\end{figure}

We point out that the condition (b) in Definition \ref{d polygon} is absent from the original statement of this Lemma in \cite{OScube}, which is necessary for its conclusion to hold ($\mu(\phi) > 1$ otherwise).  Reassuringly, this additional condition is always met in applications of that Lemma in \cite{OScube}.  Observe also that the assumption on arborescence is essential in Lemma \ref{l polygonlemma}.  For consider the punctured polygon $\D$ depicted in Figure~\ref{f poly}.  If $\D = \D(\phi)$ for some homotopy class $\phi$, then $\mu(\phi) = 1$, but the value of $\#\widehat{\M}(\phi)$ will depend on the choice of complex structure.  Thanks to Zolt\'an Szab\'o for indicating this example.  

\begin{proof}

The assertion that $\mu(\phi) = 1$ follows from Lipshitz's formula for the Maslov index (\cite{robert}, Section 4, esp. Corollary 4.3).

Now suppose that $\D$ is arborescent.  Endow $\Sigma$ with an arbitrary complex structure, and let $\D'$ denote the decorated Riemann surface obtained by erasing $\beta - (\widehat{\beta_1} \cup \cdots \cup \widehat{\beta_n})$ from $\D$, retaining its complex structure.  In the event that $\beta - (\widehat{\beta_1} \cup \cdots \cup \widehat{\beta_n})$ is disjoint from $\D$, we have $\D = \D'$, and $\#\widehat{\M}(\phi) = \pm 1$ according to \cite{OScube}, Lemma 3.11.  Moreover, the proof of \cite{OScube}, Lemma 3.11 shows how to construct a Heegaard diagram for which $\D' = \D(\phi')$ (as Riemann surfaces) for some homotopy class $\phi' \in \pi_2(x',y')$, and that $\#\widehat{\M}(\phi') = \pm 1$ in this case.

Let us compare $\widehat{\M}(\phi)$ and $\widehat{\M}(\phi')$.  Thus, let $u \in \M(\phi)$ and $u' \in \M(\phi')$.  The maximum principle implies that the image of $u$ is contained in $\text{Sym}^\ell(\D) \times \{ y_{\ell+1} \} \times \cdots \times \{ y_n \}$, where $y_{\ell+1}, \dots, y_n$ are those components for which $\widehat{\beta}_j = y_j$.  Similarly, the image of $u$ is contained in $\text{Sym}^\ell(\D') \times \{ y'_{\ell+1} \} \times \cdots \times \{ y'_m \}$.  Since $\D$ and $\D'$ are isomorphic as Riemann surfaces, it follows that we have an identification of moduli spaces $\M(\phi) \cong \M(\phi')$.  Now $\#\widehat{\M}(\phi) = \pm 1$ follows.

\end{proof}

\subsection{Determination of $E_1$.}\label{ss E_1}

Although the complex $(E_0,d_0)$ will in general depend on the choice of complex structure on $\Sigma$, Lemmas \ref{l E1domain} and \ref{l polygonlemma} provide a sufficient understanding to enable the computation of its homology $E_1$.  As a motivation, let us return to the example described in Figure~\ref{f poly}.  Suppose the domain $\D$ pictured there appears in a particular Heegaard diagram, place a basepoint in each region besides $\D_1$ and $\D_2$, and consider the Floer chain complex whose differential counts holomorphic disks which avoid these basepoints.  Choosing the Heegaard diagram suitably, we can find a quadruple of generators which coincide away from the pictured region and meet it in one of four ways: $(a,x),(a,y),(b,x),(b,y)$.  This quadruple induces a subcomplex which looks as shown:

\[ \xymatrix @R=2pc @C=-.5pc
{&& (a,x) \ar[lld]_{\D_2} \ar[rrd]^{\D} && \\
(a,y) \ar[rrd]_{\D} &&&& (b,x) \ar[lld]^{\D_2} \\
&& (b,y) && }
\]  Both of the maps labeled by $\D_2$ support a unique holomorphic representative for every choice of complex structure, thanks to the Riemann mapping theorem (a special case of Lemma \ref{l polygonlemma}).  From this alone we can conclude that the subcomplex pictured is acyclic, regardless of the value $\# \widehat{\M}(\phi)$, $\D(\phi)=\D$.

The proof of the following Lemma is in essence a generalization of this particular example.

\begin{lem}\label{l acyclic}

The subcomplex $\c_\psi$ is acyclic if it has rank $ > 1$.

\end{lem}

\begin{proof}

Suppose that $\c_\psi$ has rank $2^n >1$.  Number the (unoriented) curves $\gamma_1,\dots,\gamma_n$ in such a way that if $\gamma_i$ is contained interior to $\gamma_j$, then $i < j$.  Select a pair $x,y \in \G_\psi$ so that $C_i$ is the unique cycle in $m_x \cup m_y$ and $\gamma_i$ is oriented counterclockwise.  We pin down $r_i$ and $s_i$ (which previously depended on an arbitrary two-coloring of $C_i$) by declaring that $r_i \subset x$ and $s_i \subset y$.  For $n>1$ and $j=1,\dots,n-1$, let $\G_j$ denote those generators $x \in \G_\psi$ for which $s_{j+1} \cup \cdots \cup s_n \subset x$, let $\G_n = \G_\psi$, and for $j=1,\dots,n$, let $\c_j$ denote the subgroup of $\c_\psi$ induced on $\G_j$.  We claim that $\c_j$ is a subcomplex of $\c_\psi$.  For suppose that $x \in \G_j$ and $y \in \G_k - \G_{k-1}$ for some $k > j$.  In the $1$-cycle connecting $x$ to $y$, $\gamma_k$ is oriented clockwise and is not contained interior to any other curve.  It follows by Lemma \ref{l E1domain} that the domain $\D(\phi)$ connecting $x$ to $y$ has a negative coefficient, so that $\M(\phi) = \emptyset$ by the non-negativity principle.  Hence $d_0(\c_j) \subset \c_j$, as desired.

Consequently we obtain a filtration $\c_1 \subset \cdots \subset \c_n = \c_\psi$.  We prove by induction on $j$ that $\c_j$ is acyclic; the case $j=n$ is the statement of the Lemma.  Thus suppose $j=1$, and let $x$ and $y$ denote the Kauffman generators of $\c_1$ with $r_1 \subset x$ and $s_1 \subset y$.  We claim that the domain $\D$ connecting $x$ to $y$ is arborescent.  For if it were not, let $\alpha_1,\dots,\alpha_\ell$ denote the vertices of a non-trivial directed cycle in $\Gamma(\D)$.  We can decompose $\D = \D_1 \cup \D_2$ so that $\D_1 \cap \D_2 = \widehat{\beta}_1 \cup \cdots \cup \widehat{\beta}_\ell$ and $\D_1$ is a punctured polygon.  Let $r$ denote the corner points of $\D_1$ that are in $x$ and $s$ the other corner points of $\D_1$.  Then $y' = x \cup s - r$ is another Kauffman generator, and $\D_1$ connects the pair $x, y'$.  It follows that $y' \in \c_\psi$ as well, but the curve $\gamma = \del \D_1$ connecting the pair $x,y'$ is contained interior to $\gamma_1$, a contradiction.  Hence $\D$ is an arborescent punctured polygon; now an application of Lemma \ref{l polygonlemma} shows that $\c_1$ is acyclic.

For the induction step, suppose that $\c_j$ is acyclic for some $1 \leq j < n$, and consider the short exact sequence $$0 \to \c_j \to \c_{j+1} \to \c_{j+1} / \c_j \to 0.$$  The bijection $\G_j \to \G_{j+1} \setminus \G_j$ via $x \mapsto x \cup r_{j+1} - s_{j+1}$ induces a group isomorphism $f: \c_j \to \c_{j+1} / \c_j$.  The domain connecting a pair of generators $x,y \in \G_j$ is the same as the domain connecting the pair $f(x), f(y)$.  But if a pair of homotopy classes $\phi$ and $\phi'$ have the same domain, then $\M(\phi) \cong \M(\phi')$.  It follows that $f$ defines a chain isomorphism, so both $\c_j$ and $\c_{j+1} / \c_j$ are acyclic.  Now the long exact sequence in homology implies that $\c_{j+1}$ is acyclic. This completes the induction step and the proof of the Lemma.

\end{proof}

\begin{thm}\label{t E1}
The group $E_1$ is free abelian, generated by solitary Kauffman states.
\end{thm}

\begin{proof}

The complex $(E_0,d_0)$ decomposes into a direct sum of subcomplexes $\oplus_\psi \; \c_\psi$.  Those of rank $1$ are generated by solitary Kauffman states, and those of rank $>1$ are acyclic.

\end{proof}

\noindent As a passing remark, Theorem \ref{t E1} and the fact that the differential on $\widehat{CF}(\h(D))$ vanishes for an alternating diagram together imply that every Kauffman state in an alternating diagram is solitary.  More substantially, we have the following result.

\begin{cor}\label{c E1collapse}
If the knot $K$ admits a diagram $D$ for which the number of solitary Kauffman states is equal to $\det(K)$, then $\text{rk} \; \widehat{HFK}(K) = \det(K)$ and $\Sigma(K)$ is an L-space.
\end{cor}

\noindent Note that it is unknown whether $\text{rk} \; \widehat{HFK}(K) = \det(K)$ implies that $K$ is $\widehat{HFK}$-thin.

\begin{proof}

This is immediate from the construction of the spectral sequences, Theorem \ref{t E1}, and the inequalities $\text{rk} \; \widehat{HFK}(K) \geq \det(K)$ and $\text{rk} \; \widehat{HF}(\Sigma(K)) \geq \det(K)$.

\end{proof}

In the effort to find knots whose branched double-covers are L-spaces, it is an interesting problem to find examples of non-alternating knots satisfying the hypothesis of Corollary \ref{c E1collapse}.  Certainly there are non-alternating diagrams whose number of solitary Kauffman states equals $\det(K)$. For example, changing one crossing in the alternating diagram of the trefoil knot gives a diagram of the unknot, and for a suitable placement of the marker $p$ there will be a unique solitary Kauffman state.  However, a different placement of $p$ can result in a marked diagram with three solitary Kauffman states.


\vspace{.5in}

\section{Calculations.}\label{S calculations}

In this section we explain how to use the foregoing results to make calculations of $\widehat{HF}(\Sigma(K))$, where $K$ is a link with non-zero determinant, and in particular carry this out for some knots with $\leq 10$ crossings.

\subsection{Making a calculation.}\label{ss calculation}

In general, suppose that $\h$ is a Heegaard diagram presenting a rational homology sphere $Y$, and we wish to compute $\widehat{HF}(Y)$ as the homology of the chain complex $\widehat{CF}(\h)$.  One proceeds by enumerating the intersection points $\T_\alpha \cap \T_\beta$, calculating their associated $spin^c$ structures, and (in the case that $\h$ fits into a Heegaard triple subordinate to a link) writing down their absolute gradings.  The hard work is in determining the differential on $\widehat{CF}(\h)$.  In general this depends on the choice of complex structure on the Heegaard surface $\Sigma$ and a perturbation of the induced complex structure on $\text{Sym}^g(\Sigma)$.  However, in small examples, one can often circumvent the difficult analysis of moduli spaces en route to the answer by applying a few simple ideas which we discuss below.  For now, simply note that the holomorphicity condition restricts our attention to pairs of generators connected by a domain with non-negative coefficients.

For the case of interest, the above data can be determined in a completely algorithmic way.  We start with a marked diagram $D$ of $K$ with its regions two-colored.  From this we write down its coloring matrix and Goeritz matrices, and enumerate its Kauffman states.  Each Kauffman state corresponds to a generator of $\widehat{CF}(\h_2(D))$, whose absolute grading and associated $spin^c$ structure are determined by Theorem \ref{t grading} and $\S$\ref{ss spinc}, respectively.  Theorem \ref{t domain} enables us to determine the domain connecting a pair of generators, and in particular enables us to judge for which pairs the non-negativity condition is met.

We wrote code in Mathematica to do this (cf. a related program by Ciprian Manolescu \cite{Mprogram}).  To input a given example into the program, one begins with a marked diagram of a knot and enters in by hand the white and black graphs, the crossing incidence numbers $\mu$, and the coloring matrix.  In practice, this process is somewhat laborious, and (at least in the hands of the author) susceptible to careless errors which may take a while to detect and fix.  Therefore, it is sometimes preferable to make the calculation by another method.  For example, for the case of an alternating knot, one knows that $\Sigma(K)$ is an L-space (cf. Corollary \ref{c Lspace}), so what remains to calculate are its correction terms.  This can be done as explained in $\S \S3.5-3.6$ (indeed, calculations of this sort were carried out in \cite{OSdoublecover} and \cite{OSunknotting}).  For the case of a Montesinos knot, $\Sigma(K)$ can be realized as the boundary of a negative definite plumbing on a tree with one bad vertex, and the algorithm of \cite{OSplumbed} enables a calculation of $HF^+(\Sigma(K))$.  All but twenty knots with $\leq 10$ crossings belong to one of these two types (\cite{kawauchi}, Table F.2).  For the remainder, we calculated $\widehat{HF}(\Sigma(K))$ using our program, making use of the minimum crossing diagram for $K$ from Rolfsen's table (\cite{rolfsen}, Appendix C).

When performing a calculation, it is extremely helpful to know in advance something about the structure of the answer.  For example, if $K$ is {\em quasi-alternating}, then $\Sigma(K)$ is an L-space (\cite{OSdoublecover}, Definition 3.1 and Proposition 3.3), and if $K$ has braid index $3$, then this is nearly so: $\widehat{HF}(\Sigma(K))$ is torsion-free, and it is monic in all but at most one $spin^c$ structure \cite{quasi3braids}.  The following exceptional knots are quasi-alternating: $9_{47}$, $9_{49}$, and $10_n$, $n \in \{$148, 149, 150, 151, 155, 156, 157, 158, 159, 162, 163, 164, 165$\}$.  For $9_{47}$ and $9_{49}$, pictures are shown in (\cite{Mquasi}, Figure 10); for those of braid index 3, this follows from (\cite{quasi3braids}, Theorem 5.6). For the remainder, it is possible to find a crossing in the Rolfsen diagram of $K$ for which the resolution $K_0$ is quasi-alternating, the resolution $K_1$ is a two-component alternating link, and $\det(K_0)+\det(K_1)=\det(K)$. Furthermore, we can argue that $\Sigma(10_{160})$ is an L-space, although we could not determine whether this knot is quasi-alternating.  To see this, there is a crossing in its Rolfsen diagram for which one resolution is $9_{42}$ and the other is $7_1^2$.\footnote{In Rolfsen notation.  Also known as L7a6 in Thistlethwaite notation.}  The knot $9_{42}$ is not quasi-alternating, but $\Sigma(9_{42})$ is an L-space (one can show this either using the spanning tree model or the algorithm of \cite{OSplumbed}), $7_1^2$ is a non-split alternating link, and $\det(9_{42}) + \det(7_1^2) = \det(10_{160})$.  Now the argument of \cite{OSdoublecover}, Proposition 3.3 applies to show that $\Sigma(10_{160})$ is an L-space.  In addition to the knots listed here with L-space branched double-covers, the knots $10_{152}$ and $10_{161}$ have braid index $3$.

Thus, it is typically the case that we know in advance that $\widehat{HF}(\Sigma(K),\t)$ is free and monic, and what remains to determine is the correction term $d(\Sigma(K),\t)$ in which this group is supported.  Suppose we are in this position in a particular example.  If $\widehat{CF}(\h_2(D),\t)$ is supported in a pair of gradings for some particular $spin^c$ structure $\t$, then $\widehat{CF}(\h_2(D),\t) \cong \Z_{(d)}^k \oplus \Z_{(d \pm 1)}^{k-1}$ for some $d \in \Q$ and $k \in \Z_+$. In this case we can conclude from the Euler characteristic that $d(\Sigma(K),\t) = d$.  Alternatively, suppose that there is a unique solitary Kauffman generator in a particular $spin^c$ structure $\t$. In this case Theorem \ref{t E1} implies that $d(\Sigma(K),\t)$ equals the grading of this generator.  If $\widehat{CF}(\h_2(D),\t)$ does not have one of the two foregoing properties, it is possible that $\widehat{CF}(\h_2(D),\overline{\t})$ does, in which case we can invoke the conjugation symmetry $d(\Sigma(K),\t) =  d(\Sigma(K),\overline{\t})$.  Occasionally, neither $\widehat{CF}(\h_2(D),\t)$ nor $\widehat{CF}(\h_2(D),\overline{\t})$ has either property. In this case, we can still read off $d(\Sigma(K),\t) \; (\mod 1)$ from the grading of any generator of one of these two groups.  If $d(\Sigma(K),\s) \equiv d(\Sigma(K),\t) \; (\mod 1)$ implies that $\s = \t$ or $\overline{\t}$, then we can try to change the position of the marked point $p$ to get a new marked diagram $D'$.  Having done so, we can locate the pair of subcomplexes $\widehat{CF}(\h_2(D'),\t)$ and $\widehat{CF}(\h_2(D'),\overline{\t})$ by the gradings of their generators $(\mod 1)$, and see whether we can apply one of the preceding observations to one of the new complexes to determine $d(\Sigma(K),\t)$.

For the eighteen knots enumerated above, these techniques suffice to determine $d(\Sigma(K),\t)$ every time we know that $\widehat{HF}(\Sigma(K),\t)$ is monic.  Indeed, for many of these knots, $\det(K)$ is relatively close to the number of generators of $\widehat{CF}(\h_2(D))$, so that frequently $\widehat{CF}(\h_2(D),\t)$ consists of very few generators.  In addition, for a suitable choice of marked point $p$ in the diagrams of $10_{152}$ and $10_{160}$, and $\t_0$ the $spin$ structure on $\Sigma(K)$, the spectral sequence of $\S$\ref{S specseq} which computes $\widehat{HF}(\Sigma(K),\t_0)$ collapses at the $E_1$ term.  Thus we can easily calculate $\widehat{HF}(\Sigma(K),\t_0)$ in these two cases as well, even though these groups are not monic.

The two remaining knots $10_{153}$ and $10_{154}$ are more difficult to handle.  For one thing, $\det(10_{153})=1$, so in this case we wind up considering a chain complex with $71$ generators in a single $spin^c$ structure.  In addition, we have no a priori result asserting that $\widehat{HF}(\Sigma(K))$ is torsion-free for either of these knots, although this turns out to be the case.  Instead, we must analyze $\widehat{CF}(\h_2(D))$ more closely, applying the foregoing techniques with some additional observations.  For example, it is commonly the case that we have a pair of generators $x$ and $y$ for which $d_0(x) = \pm y$ and $x$ is not connected to any other generator by a non-negative domain.  Then $d(x) = \pm y$, and $d^2=0$ implies that $d(y) = 0$.  This latter fact may not be so clear judging from domains of maps, as $y$ could be connected to some other generator $z$ by a non-negative domain $\D$.  In any event, we learn that $\# \widehat{\M}(\phi) = 0$ for any $\phi$ with $\D = \D(\phi)$.  That is, if there is some other pair of generators $y'$ and $z'$ connected by this same domain $\D$, then we can conclude that $d(y')$ does not involve any contribution from $z'$.  One can run this sort of argument and close in on the answer by comparing what the different complexes $\widehat{CF}(\h_2(D))$ reveal by moving the marked point around.  As a final resort, it is necessary on occasion to draw the domain connecting a pair of generators which is not counted by $d_0$, recognize it as an arborescent punctured polygon, and apply Lemma \ref{l polygonlemma}.  With a good deal of perseverance, these techniques enable the determination of $\widehat{HF}(\Sigma(K))$, $K = 10_{153}$ and $10_{154}$.

\subsection{Results.}\label{ss results}

The calculation of $\widehat{HF}(\Sigma(K))$ for the twenty non-alternating, non-Montesinos knots $K$ with $\leq 10$ crossings is reported as follows.  We identify $Spin^c(\Sigma(K)) \cong H^2(\Sigma(K))$ via $\t \mapsto c_1(\t)$, and when $\Sigma(K)$ is an L-space, we record the element $$\sum_\t d(\Sigma(K),\t) \cdot \t \in \Q[H^2(\Sigma(K))].$$  We express an element of $H^2(\Sigma(K))$ as a monomial $x^j$ or $x^j y^k$, depending on the number of generators of this group. In the four instances in which $\Sigma(K)$ is not an L-space, $\widehat{HF}(\Sigma(K),\t)$ has rank one except in the $spin$ structure $\t_0$, and $\widehat{HF}(\Sigma(K),\t_0)$ takes the form $\Z_{(d)}^k \oplus \Z_{(d \pm 1)}^{k-1}$, where $d = d(\Sigma(K),\t_0)$.  In this case we write $$\underline{d}^k + (d \pm 1)^{k-1} + \sum_{\t \ne t_0} d(\Sigma(K),\t) \cdot \t.$$  Thus, for example, $\widehat{HF}(\Sigma(10_{153})) = \widehat{HF}(\Sigma(10_{153},\t_0)) \cong \Z_{(0)}^3 \oplus \Z_{(1)}^2$.  We briefly point out that these calculations can be used to determine $HF^+(\Sigma(K))$ from the algebraic fact that $\widehat{HF}(Y,\t) \cong \Z_{(d)}^k \oplus \Z_{(e)}^{k-1}$ with $e = d \pm 1$ implies that $HF^+(Y,\t) \cong \mathcal{T}^+_d \oplus \Z_{(\min \{d, e \})}^{k-1}$.  When $H^2(\Sigma(K))$ is cyclic, the conjugation symmetry manifests itself as the equality between the coefficients on terms $x^j$ and $x^{\det(K)-j}$.  When $H^2(\Sigma(K)) \cong \Z / a \Z \oplus \Z / b \Z$, we compare the coefficients on terms $x^j y^k$ and $x^{a-j} y^{b-k}$ instead.  The rational number appearing in front of each expression is $1 / \det(K)$ or $1 / ( 2 \det(K) )$ when $H^2(\Sigma(K))$ is cyclic, and $1 / b$ or $1 / (2b)$ when $H^2(\Sigma(K)) \cong \Z / a \Z \oplus \Z / b \Z$ with $a \; | \; b$.

\vspace{.5in}

\begin{enumerate}

\item[\fbox{$9_{47}$}] ${1 \over 18} \cdot (-9 + 3 x + 3 x^2 + 7 y - 17 x y - 17 x^2 y - 17 y^2 -
    5 x y^2 - 5 x^2 y^2 - 9 y^3 + 3 x y^3 + 3 x^2 y^3 -
    5 y^4 + 7 x y^4 + 7 x^2 y^4 - 5 y^5 + 7 x y^5 +
    7 x^2 y^5 - 9 y^6 + 3 x y^6 + 3 x^2 y^6 - 17 y^7 -
    5 x y^7 - 5 x^2 y^7 + 7 y^8 - 17 x y^8 - 17 x^2 y^8)$

\item[\fbox{$9_{49}$}] ${1 \over 5} \cdot (5 + 3 x - 3 x^2 - 3 x^3 + 3 x^4 - y - 3 x y + x^2 y +
    x^3 y - 3 x^4 y + y^2 - x y^2 + 3 x^2 y^2 + 3 x^3 y^2 -
    x^4 y^2 + y^3 - x y^3 + 3 x^2 y^3 + 3 x^3 y^3 -
    x^4 y^3 - y^4 - 3 x y^4 + x^2 y^4 + x^3 y^4 - 3 x^4 y^4)$

\item[\fbox{$10_{148}$}] ${1 \over 62} \cdot (31 - 33 x + 23 x^{2} - 49 x^{3} - x^{4} + 43 x^{5} - 41 x^{6} - 5 x^{7} +
    27 x^{8} - 69 x^{9} - 45 x^{10} - 25 x^{11} - 9 x^{12} + 3 x^{13} +
    11 x^{14} + 15 x^{15} + 15 x^{16} + 11 x^{17} + 3 x^{18} - 9 x^{19} -
    25 x^{20} - 45 x^{21} - 69 x^{22} + 27 x^{23} - 5 x^{24} - 41 x^{25} +
    43 x^{26} - x^{27} - 49 x^{28} + 23 x^{29} - 33 x^{30})$

\item[\fbox{$10_{149}$}] ${1 \over 41} \cdot (41 + 35 x + 17 x^{2} - 13 x^{3} + 27 x^{4} - 27 x^{5} - 11 x^{6} -
    7 x^{7} - 15 x^{8} - 35 x^{9} + 15 x^{10} - 29 x^{11} - 3 x^{12} +
    11 x^{13} + 13 x^{14} + 3 x^{15} - 19 x^{16} + 29 x^{17} - 17 x^{18} +
    7 x^{19} + 19 x^{20} + 19 x^{21} + 7 x^{22} - 17 x^{23} + 29 x^{24} -
    19 x^{25} + 3 x^{26} + 13 x^{27} + 11 x^{28} - 3 x^{29} - 29 x^{30} +
    15 x^{31} - 35 x^{32} - 15 x^{33} - 7 x^{34} - 11 x^{35} - 27 x^{36} +
    27 x^{37} - 13 x^{38} + 17 x^{39} + 35 x^{40})$

\item[\fbox{$10_{150}$}] ${1 \over 29} \cdot (-29 + 5 x - 9 x^{2} - 13 x^{3} - 7 x^{4} + 9 x^{5} - 23 x^{6} +
    13 x^{7} + x^{8} - x^{9} + 7 x^{10} - 33 x^{11} - 5 x^{12} - 25 x^{13} -
    35 x^{14} - 35 x^{15} - 25 x^{16} - 5 x^{17} - 33 x^{18} + 7 x^{19} -
    x^{20} + x^{21} + 13 x^{22} - 23 x^{23} + 9 x^{24} - 7 x^{25} - 13 x^{26} -
    9 x^{27} + 5 x^{28})$

\item[\fbox{$10_{151}$}] ${1 \over 86} \cdot (-43 + 25 x + 57 x^{2} + 53 x^{3} + 13 x^{4} - 63 x^{5} - 3 x^{6} +
    21 x^{7} + 9 x^{8} - 39 x^{9} + 49 x^{10} - 71 x^{11} - 55 x^{12} -
    75 x^{13} + 41 x^{14} - 51 x^{15} - 7 x^{16} + x^{17} - 27 x^{18} -
    91 x^{19} - 19 x^{20} + 17 x^{21} + 17 x^{22} - 19 x^{23} - 91 x^{24} -
    27 x^{25} + x^{26} - 7 x^{27} - 51 x^{28} + 41 x^{29} - 75 x^{30} -
    55 x^{31} - 71 x^{32} + 49 x^{33} - 39 x^{34} + 9 x^{35} + 21 x^{36} -
    3 x^{37} - 63 x^{38} + 13 x^{39} + 53 x^{40} + 57 x^{41} + 25 x^{42})$

\item[\fbox{$10_{152}$}] ${1 \over 22} \cdot (\underline{- 11}^2 + 11 + 9x -19x^2 -7 x^3 + x^4 + 5x^5 + 5 x^6 + x^7 - 7 x^8 - 19x^9 + 9x^{10})$

\item[\fbox{$10_{153}$}] $\underline{0}^3 + 1^2$

\item[\fbox{$10_{154}$}] ${1 \over 13} \cdot (\underline{13}^2 + 0 +9x^1 +3x^2 +3x^3 + x^4+17x^5 -x^6 -x^7 +17x^8 +x^9 +3x^{10} +3x^{11} +9x^{12})$

\item[\fbox{$10_{155}$}] ${1 \over 5} \cdot (0 -4 x - 6 x^2 - 6 x^3 - 4 x^4 - 4 y + 2 x y + 2 x^4 y -
    6 y^2 - 2 x^2 y^2 - 2 x^3 y^2 - 6 y^3 - 2 x^2 y^3 -
    2 x^3 y^3 - 4 y^4 + 2 x y^4 + 2 x^4 y^4)$

\item[\fbox{$10_{156}$}] ${1 \over 35} \cdot (35 + 19 x - 29 x^{2} + 31 x^{3} + 59 x^{4} + 55 x^{5} + 19 x^{6} -
    49 x^{7} - 9 x^{8} - x^{9} - 25 x^{10} + 59 x^{11} - 29 x^{12} -
    9 x^{13} - 21 x^{14} + 75 x^{15} - x^{16} + 31 x^{17} + 31 x^{18} -
    x^{19} + 75 x^{20} - 21 x^{21} - 9 x^{22} - 29 x^{23} + 59 x^{24} -
    25 x^{25} - x^{26} - 9 x^{27} - 49 x^{28} + 19 x^{29} + 55 x^{30} +
    59 x^{31} + 31 x^{32} - 29 x^{33} + 19 x^{34})$

\item[\fbox{$10_{157}$}] ${1 \over 7} \cdot (-7 + x - 3 x^2 - 5 x^3 - 5 x^4 - 3 x^5 + x^6 - y - x y +
    x^2 y + 5 x^3 y - 3 x^4 y + 5 x^5 y + x^6 y + 3 y^2 -
    5 x y^2 + 3 x^2 y^2 - x^3 y^2 - 3 x^4 y^2 - 3 x^5 y^2 -
    x^6 y^2 + 5 y^3 + 3 x y^3 + 3 x^2 y^3 + 5 x^3 y^3 -
    5 x^4 y^3 + x^5 y^3 - 5 x^6 y^3 + 5 y^4 - 5 x y^4 +
    x^2 y^4 - 5 x^3 y^4 + 5 x^4 y^4 + 3 x^5 y^4 +
    3 x^6 y^4 + 3 y^5 - x y^5 - 3 x^2 y^5 - 3 x^3 y^5 -
    x^4 y^5 + 3 x^5 y^5 - 5 x^6 y^5 - y^6 + x y^6 +
    5 x^2 y^6 - 3 x^3 y^6 + 5 x^4 y^6 + x^5 y^6 - x^6 y^6)$

\item[\fbox{$10_{158}$}] ${1 \over 45} \cdot (0 + 32 x + 38 x^{2} + 18 x^{3} - 28 x^{4} - 10 x^{5} - 18 x^{6} + 38 x^{7} -
    22 x^{8} - 18 x^{9} - 40 x^{10} + 2 x^{11} + 18 x^{12} + 8 x^{13} -
    28 x^{14} + 2 x^{16} - 22 x^{17} + 18 x^{18} + 32 x^{19} + 20 x^{20} -
    18 x^{21} + 8 x^{22} + 8 x^{23} - 18 x^{24} + 20 x^{25} + 32 x^{26} +
    18 x^{27} - 22 x^{28} + 2 x^{29} - 28 x^{31} + 8 x^{32} + 18 x^{33} +
    2 x^{34} - 40 x^{35} - 18 x^{36} - 22 x^{37} + 38 x^{38} - 18 x^{39} -
    10 x^{40} - 28 x^{41} + 18 x^{42} + 38 x^{43} + 32 x^{44})$

\item[\fbox{$10_{159}$}] ${1 \over 78} \cdot (39 + 23 x - 25 x^{2} + 51 x^{3} - 61 x^{4} - 49 x^{5} - 69 x^{6} +
    35 x^{7} - 49 x^{8} - 9 x^{9} - x^{10} - 25 x^{11} - 81 x^{12} -
    13 x^{13} + 23 x^{14} + 27 x^{15} - x^{16} - 61 x^{17} + 3 x^{18} +
    35 x^{19} + 35 x^{20} + 3 x^{21} - 61 x^{22} - x^{23} + 27 x^{24} +
    23 x^{25} - 13 x^{26} - 81 x^{27} - 25 x^{28} - x^{29} - 9 x^{30} -
    49 x^{31} + 35 x^{32} - 69 x^{33} - 49 x^{34} - 61 x^{35} + 51 x^{36} -
    25 x^{37} + 23 x^{38})$

\item[\fbox{$10_{160}$}] ${1 \over 21} \cdot (-21 + 5 x - x^{2} + 3 x^{3} - 25 x^{4} - x^{5} - 9 x^{6} - 7 x^{7} +
    5 x^{8} - 15 x^{9} - 25 x^{10} - 25 x^{11} - 15 x^{12} + 5 x^{13} -
    7 x^{14} - 9 x^{15} - x^{16} - 25 x^{17} + 3 x^{18} - x^{19} + 5 x^{20})$

\item[\fbox{$10_{161}$}] ${1 \over 5} \cdot (\underline{-5}^2+ 0 -x +x^2 + x^3 -x^4) $

\item[\fbox{$10_{162}$}] ${1 \over 70} \cdot (35 + 59 x - 9 x^{2} - 29 x^{3} - x^{4} - 65 x^{5} + 59 x^{6} - 49 x^{7} +
    31 x^{8} + 19 x^{9} + 55 x^{10} - x^{11} - 9 x^{12} + 31 x^{13} -
    21 x^{14} - 25 x^{15} + 19 x^{16} - 29 x^{17} - 29 x^{18} + 19 x^{19} -
    25 x^{20} - 21 x^{21} + 31 x^{22} - 9 x^{23} - x^{24} + 55 x^{25} +
    19 x^{26} + 31 x^{27} - 49 x^{28} + 59 x^{29} - 65 x^{30} - x^{31} -
    29 x^{32} - 9 x^{33} + 59 x^{34})$

\item[\fbox{$10_{163}$}] ${1 \over 102} \cdot (-51 - 91 x - 7 x^{2} - 3 x^{3} - 79 x^{4} - 31 x^{5} - 63 x^{6} +
    29 x^{7} + 41 x^{8} - 27 x^{9} + 29 x^{10} + 5 x^{11} - 99 x^{12} -
    79 x^{13} + 65 x^{14} - 75 x^{15} - 91 x^{16} + 17 x^{17} + 45 x^{18} -
    7 x^{19} + 65 x^{20} + 57 x^{21} - 31 x^{22} + 5 x^{23} - 39 x^{24} +
    41 x^{25} + 41 x^{26} - 39 x^{27} + 5 x^{28} - 31 x^{29} + 57 x^{30} +
    65 x^{31} - 7 x^{32} + 45 x^{33} + 17 x^{34} - 91 x^{35} - 75 x^{36} +
    65 x^{37} - 79 x^{38} - 99 x^{39} + 5 x^{40} + 29 x^{41} - 27 x^{42} +
    41 x^{43} + 29 x^{44} - 63 x^{45} - 31 x^{46} - 79 x^{47} - 3 x^{48} -
    7 x^{49} - 91 x^{50})$

\item[\fbox{$10_{164}$}] ${1 \over 45} \cdot (0 -32 x - 38 x^{2} - 18 x^{3} + 28 x^{4} + 10 x^{5} + 18 x^{6} -
    38 x^{7} + 22 x^{8} + 18 x^{9} - 50 x^{10} - 2 x^{11} - 18 x^{12} -
    8 x^{13} + 28 x^{14} - 2 x^{16} + 22 x^{17} - 18 x^{18} - 32 x^{19} -
    20 x^{20} + 18 x^{21} - 8 x^{22} - 8 x^{23} + 18 x^{24} - 20 x^{25} -
    32 x^{26} - 18 x^{27} + 22 x^{28} - 2 x^{29} + 28 x^{31} - 8 x^{32} -
    18 x^{33} - 2 x^{34} - 50 x^{35} + 18 x^{36} + 22 x^{37} - 38 x^{38} +
    18 x^{39} + 10 x^{40} + 28 x^{41} - 18 x^{42} - 38 x^{43} - 32 x^{44})$

\item[\fbox{$10_{165}$}] ${1 \over 78} \cdot (-39 + 5 x - 19 x^{2} + 45 x^{3} + 41 x^{4} - 31 x^{5} - 15 x^{6} +
    89 x^{7} - 31 x^{8} + 93 x^{9} - 7 x^{10} - 19 x^{11} + 57 x^{12} +
    65 x^{13} + 5 x^{14} + 33 x^{15} - 7 x^{16} + 41 x^{17} + 21 x^{18} +
    89 x^{19} + 89 x^{20} + 21 x^{21} + 41 x^{22} - 7 x^{23} + 33 x^{24} +
    5 x^{25} + 65 x^{26} + 57 x^{27} - 19 x^{28} - 7 x^{29} + 93 x^{30} -
    31 x^{31} + 89 x^{32} - 15 x^{33} - 31 x^{34} + 41 x^{35} + 45 x^{36} -
    19 x^{37} + 5 x^{38})$

\end{enumerate}

\vspace{.5in}

\section{Speculation.}\label{S speculation}

Theorem \ref{t grading} expresses the Maslov grading of a Kauffman generator $x$ in $\widehat{CF}(\h_2)$: $$\textup{gr}(x) = \delta(x) + { q(v^W_x) -2m - 3 \sigma(G_W) \over 4}.$$  Notice that for a pair of generators $x$ and $y$ connected by the differential $d_0$, we have $\textup{gr}(x)-\textup{gr}(y) = 1$, and the domain connecting $x$ to $y$ is supported on top of $\h_2$, implying that $v_x^W = v_y^W$.  It follows in this case that $\delta(x)-\delta(y)=1$, i.e. {\em $d_0$ lowers $\delta$ by $1$}.  In fact, in small examples we found that it is always the case that the full differential $d$ lowers $\delta$ by $1$.  

This phenomenon is particularly striking in examples when the reduced white graph $\widetilde{W}$ is a tree and the Goeritz matrix $G_W$ is negative definite.  In other words, $\Sigma(K)$ is the boundary of a negative definite plumbing on a tree.  As a representative case, consider the well-known plumbing on the $E_8$ diagram, in which every vertex receives weight $-2$. This is the reduced white graph of a particular diagram of the $(3,5)$ torus knot, and the manifold it describes is the Poincar\'e dodecahedral space $\Sigma(2,3,5)$.  In this case, the spanning tree model has $61$ generators, each supported in the unique $spin^c$ structure and (notably) in different gradings.  We represent each generator by a dot $\bullet$ and draw an arrow between two dots $\bullet \to \bullet$ if the generator to the left is connected to the one on the right by a non-negative domain.  Each arrow therefore represents a potential differential connecting a pair of generators.  In this way, the complex $\widehat{CF}(E_8)$ decomposes as a union of linear chains $\bullet \to \bullet \to \cdots \to \bullet$, one of which terminates on the unique generator in grading $-2$.  Now, it is known that $\widehat{HF}(\Sigma(2,3,5)) \cong \Z_{(-2)}$ (cf. \cite{OSgrading}, Section 8 or \cite{OSplumbed}, Section 3.2).  This enables us to conclude that in each linear chain, the maps alternate between isomorphisms and trivial maps.  Inspecting the complex more closely, we notice that the trivial maps connect pairs of generators for which $\delta(x)-\delta(y) = -1$, while the isomorphisms connect pairs for which $\delta(x)-\delta(y) = 1$.  These facts are not at all obvious by considering the domains of maps.  For example, one of the domains corresponding to a non-trivial map has a region with coefficient $11$, which is remarkably large for such a small example.

In spite of the very limited evidence, it is tempting to consider the following conjecture.

\begin{conj}\label{conj deltagrading}

In the complex $\widehat{CF}(\h_2(D))$, the differential $d$ lowers $\delta$ by $1$.  The induced $\delta$-grading on the homology group $\widehat{HF}(\Sigma(K))$ is independent of $D$ up to an overall shift, and can be pinned down to a well-defined integer absolute grading.

\end{conj}

In case Conjecture~\ref{conj deltagrading} were true, we would naturally like to compare the $\delta$-grading on $\widehat{HF}(\Sigma(K))$ with that on the Khovanov homology $Kh(K)$ and knot Floer homology $\widehat{HFK}(K)$.  In particular, we might ask whether the inequality of ranks $\text{rk} \; Kh_{red}(K, \Z / 2) \geq \text{rk} \; \widehat{HF}(\Sigma(K),\Z/2)$ implied by the Ozsv\'ath-Szab\'o spectral sequence persists at the level of $\delta$-graded groups.  As it stands, that spectral sequence is not known to reflect any of the finer features of the two theories involved: the integer bigrading on $Kh_{red}(K,\Z/2)$, or the rational grading on and decomposition into $spin^c$ structures of $\widehat{HF}(\Sigma(K))$.  Note that the quantity $\delta(x)$ is the same as the $v$-grading in the spanning tree complex for Khovanov homology in \cite{CK}.

We also note that in many small examples, such as knots on $\leq 12$ crossings and several pretzel and torus knots, we observed an inequality of ranks $\text{rk} \; \widehat{HFK}(K) \geq \text{rk} \; \widehat{HF}(\Sigma(K))$.  For knots on $\leq 12$ crossings this is no surprise, since $\text{rk} \; \widehat{HFK}(K) = \text{rk} \; Kh(K)$ in this range.  It is tempting to speculate on the relationship between $\widehat{HF}(\Sigma(K))$ and $\widehat{HFK}(K)$, due to the strong resemblance between the Heegaard diagrams $\h(D)$ presenting $K \subset S^3$ and $\h_2(D)$ presenting $\Sigma(K)$.  The fact that the $d_0$ differential in the two spectral sequences of Section~\ref{S specseq} is the same is particularly curious, since we know a priori that it respects the Alexander filtration grading on $\widehat{HFK}$ (absent in $\widehat{HF}(\Sigma(K))$) and the decomposition into $spin^c$ structures of $\widehat{HF}(\Sigma(K))$ (absent in $\widehat{HFK}$).

In a different direction, Andr\'as N\'emethi has proposed a very interesting model for the Floer homology $HF^+$ of a negative-definite plumbed rational homology sphere (\cite{nemethi}, esp. Section 3 and Conjecture 5.2.4).  In this model, $HF^+$ is enhanced by an additional non-negative integer grading.  Any plumbed rational homology sphere takes the form $\Sigma(K)$ for some link $K$ with $\det(K) \ne 0$.  It is very interesting to compare N\'emethi's model to our own, especially with a view towards our Conjecture \ref{conj deltagrading} and his Conjecture 5.2.4.

\end{document}